\documentclass[12pt,a4paper]{article}
\usepackage[utf8]{inputenc}
\usepackage{amsmath}
\usepackage{relsize}
\usepackage{amsthm} 
\usepackage{amsfonts}
\usepackage{amssymb}
\usepackage{url}
\usepackage{graphicx}
\usepackage{listings}
\usepackage{color}
\usepackage{mathtools}
\usepackage{physics}
\usepackage[bottom=2.5cm, right=2.5cm, left=2.5cm, top=2.5cm]{geometry}

\ProvideTextCommand{\DJ}{OT1}{\leavevmode\raisebox{-.5ex}{\makebox[0pt][l]{\hskip-.07em\accent"16\hss}}D}
\newcommand{\tens}[1]{%
	\mathbin{\mathop{\otimes}\displaylimits_{#1}}%
}

\newtheorem{remark}{Remark}[section]

\newtheorem*{conclusion}{Conclusion}
\newtheorem*{conjecture}{Conjecture}
\newtheorem{definition}{Definition}[section]

\newtheorem{lemma}{Lemma}[section]

\newtheorem*{lemma5}{Lemma 5.5'}

\newtheorem{corollary}{Corollary}[section]
\newtheorem{theorem}{Theorem}[section]

\newtheorem{proposition}{Proposition}[section]
\newtheorem{example}{Example}[section]
\newtheorem*{question}{Question}

\newtheorem*{theorema}{Theorem}
\newtheorem*{thA}{Theorem A}

\newtheorem*{thover}{Theorem 1.1$^\prime$}

\author{Alexander Kushkuley \\ akushkuley@gmail.com}
\title{Identities  of  Irreducible Representations and Gassmann Equivalence}
\begin{document}
	\maketitle
\;\;\; \;\;\;  \;\;\; \;\;\; \;\;\; \;\;\;  \;\;\; \;\;\;  \;\;\; \;\;\;  \;\;\; \;\;\; \;\;\; \;\;\;  \;\;\; \;\;\;                     dedicated to Boris Plotkin
\newline\newline
	\begin{abstract}
\noindent Identities of complex irreducible representations of finite groups can be explicitly constructed from character value sets. Among other things, these identities determine   representations up to Gassmann equivalency. 
Some examples of identities related to spherical space forms and to representations of finite  $p$-groups are presented. 
Some old results on irreducible representations with the same identities are revisited 
	\end{abstract}
	
\tableofcontents	
	\section{Introduction}
	   \numberwithin{equation}{section}
	The notion of \textit{identical relations} (or simply \textit{identities}) of/in group representations  was introduced by Plotkin in a context of universal algebra (cf. \cite{P} and references therein).  
	Let $ \Gamma $ be a  group and let $ \rho : \Gamma \rightarrow \textnormal{GL}(V) $  be a representation of $\Gamma$ in a vector space $V $ over a field $k$. 
	Let $ kF $ be a group algebra of a free group $ F = F(Y) $ with a countable set of free generators $ Y $. An element $ u = u(y_1,\cdots , y_t) \in kF, \; y_i \in Y, \;   i = 1,\cdots , t $ is said to be an identity of/in the  representation  $\rho$ (cf. \cite{P}) if 
	$  u( \rho( \mathfrak{h}(y_1)),  \cdots  \rho( \mathfrak{h}(y_t)) ) $ vanishes in $ \textnormal{End}(V)$  for  any homomophism 
	$ \mathfrak{h} : F(Y) \rightarrow G  $.
	
	\begin{example} (cf. \cite{P}-\cite{PK}, \cite{V}).
		Let $ f_i(y_1, \cdots, y_t ) \in F(Y), \; i = 1, \cdots , k $.  If a disjunctive formula 
		\begin{align}
		w = 	f_1(y_1, \cdots , y_t ) = 1) \vee ( f_2(y_1, \cdots , y_t ) = 1) 
			\vee  (f_k(y_1, \cdots , y_t ) = 1   
		\end{align}
		holds identically in a group $G$ then \textbf{disjunctive identity}  
		\begin{align}
		 u(w) = u_0( f_1 - 1) u_1  \cdots  u_{k-1} (f_k - 1) u_k , \; u_i \in Y
		\end{align} 	
		holds in any $G$-representation. 

The identity (1.2) is satisfied by a faithful   irreducible representation  of a group $ G $ if and only if disjunction (1.1)  identically holds in $ G $  (cf. \cite{P}-\cite{PK}, \cite{V}).
As a concrete example of (1.1-2), set 

\begin{align}
	\!\!\!\!		\mathcal{C}_m \equiv 	\mathcal{C}_m(y_1, \cdots, y_m) = u_0\prod_{1 \leq i<j \leq m}^m ( y_i - y_j) u_{ij}
\end{align}
\noindent where $  y_i, \; u_0, u_{ij} $ are pairwise distinc free variables in $Y$.
	The identity (1.3) holds in any representation of any group of order less than   $ m  $.  The identity (1.3) holds in a  faithful   irreducible representation of  $ G $ if and only if  $ |G| < m  $. 
\end{example}

\begin{example} Any polynomial identity that holds in $ \textnormal{End}(V) $  
is  an identity of the representation $\rho : G \rightarrow \textnormal{GL}(V)$. In particular,	any non-commutative polynomial that is multilinear and skew-symmetric in $ n^2 +1$ of its variables is an identity of any group representation of dimension $\leq n$. 
Another example of this kind is 
a standard  identity of order $2n$
\begin{align}
	s_{2n}( y_1, \cdots y_{2n}) \equiv    \sum_{\sigma \in S_{2n}} \textnormal{sign}(\sigma)y_{\sigma(1)} \cdots y_{\sigma(2n)} \nonumber
\end{align}
 It is well known (cf. e.g. \cite{Artin}) that identity $ s_{2n} $ holds in a full matrix algebra $M_r(k) $ and hence in any group representation of dimension  $ r \leq n $.
\end{example}
\noindent 
It was shown by Plotkin (\cite{PK}, see also \cite{P} and a well written exposition in \cite{V}) that  all identities of a non-modular regular representation of a finite group follow from a finite set of identities of a kind described in examples above. 

Turning to irreducible representations (irreps), it is natural to try to use trace (character) as a blueprint for a representation identity.  
 The main goal of this paper (cf. sections 1, 3 and 4) is  to present a few character-related recipes of writing identities  
 of irreps of finite groups.  
 It turns out, for example, that \textit{Gassmann equivalence} of exact irreps  corresponds to  a well defined set of identities (see Proposition 3.1). In particular, we are interested in finding explicit specifications of  identities that distinguish between irreps. This can be done only up to  \textit{similarity} of representations - the following definition and subsequent discussion will make it clear. 

\begin{definition}(\cite{P}).
Two representations    $ \rho_i : G_i \rightarrow \textnormal{Aut}(V_i), \; i =1,2 $ 
	are called \textbf{similar}
	\footnote[1]{ This terminology is not standard. For example,   similar representations are called \textbf{isomorphic} in \cite{P}, \textbf{automorphically equivalent} in \cite{Sutton} and
	\textbf{equivalent modulo automorphisms} in \cite{Ikeda}} 
	  if there is an isomorphism  
	$ \alpha: G_1 \rightarrow G_2$ 
	such that representations $ \rho_2  \alpha $ and   
	$ \rho_1 $ are equivalent. Two characters  $ \chi_1, \chi_2 $ of a finite group $G$  are called similar if there is $ \alpha \in \textnormal{Aut}(G) $ such that 
	$ \chi_2 \alpha = \chi_1  $. 

\end{definition}
\begin{remark}
	Exact representations  $ \rho_i : G_i \rightarrow \textnormal{GL}(V), \; i =1,2 $ are similar if and only if the images $ \rho_1(G_1) $ and $ \rho_2(G_2) $ are conjugate in $ \textnormal{GL}(V) $. It should be noted that the notion of similarity of representation is quite different from that of representation equivalence. Some examples emphasizing the difference can be found in Appendix 2   
\end{remark}
\noindent It is obvious  (cf. \cite{P}) that similar representations have  same identities and an important fact established by  Plotkin is that the reverse of this statement is also true for faithful nonmodular irreps of  finite groups  (see \cite{P} and sections 5, 6 below)  
\begin{theorem} (Plotkin \cite{P} (Prop. 2.2.4.1)). 
Exact nonmodular irreps of finite groups  that have the same identities are similar
\end{theorem}
\noindent A new proof of this theorem that allows generalization to $p$-\textit{defect zero} modular irreps  is presented in section 5.
Similar results for irreps of algeraic groups are considered in section 6. A concrete  example of an identity of the natural representation of $ \textnormal{SL}_2(k) $ can be found in Appendix 1.   
   
 We wind up this introduction by quoting  a simple fact (first observed by Plotkin) that demonstrates usefulness  of "disjunctive identities"  (1.2) (cf. \cite{P}-\cite{PK}, \cite{V}). 
  Let $ G,\; H  $ be finite groups and let $ \rho, \; \sigma $ be   faithful irreducible   representations of $G$ and $H$ respectfully.

\begin{lemma} (cf. \cite{P}-\cite{PK}, \cite{V}). Let  $A \subset \textnormal{End}(V)$ be a linear span of $ \rho(G) $.
	\begin{enumerate}
		\item[(1)] Identity $u(w)$ (1.2) is satisfied by the  irrep $\rho$ if and only if the group $G$ identically satisfies the disjunction $w$ (1.1)
		\item[(2)] Suppose that  disjunction (1.1) does not hold in $G$.
		Then for any $ a \in A$ there are  $ g_1, \cdots, g_t \in \rho(G) $ and
		$ a_0, a_1, \cdots , a_k \in A $    such that 
\begin{align}
	a = a_0 (f_1(g_1, \cdots , g_t ) -1 )a_1 \cdots  a_{k-1}(f_k(g_1, \cdots , g_t ) - 1) a_k  
\end{align}	
		
			\item[(3)] If in addition to condition (2) the ground field $k$ is algebraically closed, then (1,4) holds for any $ a \in \textnormal{End}(V) $ and in particular, there are   
	 $ h_0, \cdots , h_k \in \rho(G) $ such that 
			$$ \textnormal{tr}(\; h_0 (f_1(g_1, \cdots , g_t ) -1 )h_1 \cdots  h_{k-1}(f_k(g_1, \cdots , g_t ) - 1) h_k  \;  ) \neq 0 $$ 
				
			\item[(4)] There is a finite set   $ \mathfrak{D}(G) \subset kF(Y) $ of expressions (1.2) such that

		\begin{enumerate} 
			\item[(a)] identities $ \mathfrak{D}(G)$  hold in any representation of $G$ 
			\item[(b)] if identities $ \mathfrak{D}(G)$ hold in  irrep  $\sigma$  then $ H $ is a section of $ G$ and  in particular, if $ |H| =  |G| $ then $ G \approx H $
		\end{enumerate}
		
	\end{enumerate}
\end{lemma}

\noindent Let's briefly mention some details (cf.  \cite{P}-\cite{PK}). It is clear that  disjunction (1.1) implies identity (1.2). On the other hand, if (1.1) does not hold in $ G $   then by exactness of the irrep $\rho$ there are  $ g_1, \cdots, g_t \in \rho(G) $ such that $  f_i(g_1, \cdots , g_t ) -1  \neq 0 $ for all $ i = 1, \cdots , k $ and  statement  (2) follows from the simplicity of  $A$. Over algebraically closed field $ A = \textnormal{End}(V)$ and this yields (3).
The statement (4) is a standard fact of model theory:  
a  group  $\Gamma$ is a section (factor group of a subgroup) of a group  $\Sigma$ if and only if $\Gamma$ satisfies all disjunctions (1.1) satisfied by $ \Sigma $ (cf. e.g.\cite{Cohn}).

\begin{remark}
	A simple (albeit inefficient) algorithm that determines the (finite) set of identities $ \mathfrak{D}(\Gamma)$ for a finite group $\Gamma$ can be found in \cite{P}-\cite{PK} (see also \cite{V})
\end{remark}

\subsection{ Basic Setup }	
From now on, unless explicitly stated otherwise, we will assume that the ground  filed $k$ is the field of  complex numbers $\mathbb{C} $. Representations of finite (compact) groups  will be assumed to be unitary if necessary. 
\begin{remark}
	It is well known that any continuous representation of a compact group over the field of complex (real) numbers is equivalent to a unitary (orthogonal) representation. 
	It is easy to see that equivalent unitary (orthogonal) representations are unitary (orthogonally) equivalent  and essentially the same is true for similar unitary (orthogonal) representations (cf. e.g \cite{Sutherland} and section 6.3)
\end{remark}

\noindent 

  Let $ G $ be a finite group of order $m>1$  and let $ \rho : G \rightarrow \textbf{U}(V) $ be its faithful irreducible unitary representation (unirrep) in a vector space $ V$ of dimension $  n$.  
Denote by $  \chi_{\rho}$   the character of representation $\rho $.
Let  $C_1, \; C_2, \cdots , \;  C_s \subset G $  be all  classes of conjugate elements  in $G$ so    
that  $ \textnormal{range}(\chi_{\rho}) $ is a finite set of complex numbers $  \{ \chi_1, \cdots, \chi_k \} \subset  \mathbb{C}$ where  $  k= |\textnormal{range}(\chi_{\rho})| \leq s $. 
\newline\newline
\noindent
\textbf{ Let us fix this otherwise arbitrary  unirrep  $ \rho : G \rightarrow \textbf{U}(V) $ as a notational convenience. For the same  reason let's keep around another  exact  unirrep 
$ \sigma : H \rightarrow \textbf{U}(W) $ of a finite group $H$ in a finite dimensional space $ W $.   
}

\section{Character Identities}

Take free variables  $ x \in Y, \; Y_m= \{ y_1, \cdots , y_m \} \subset Y $ and  
set 
\begin{align}
	\Psi_m \equiv \Psi_m (x, Y_m) = 
	y_1 x y^{-1}_1  + 	y_2 x y^{-1}_2  + \cdots
	+ y_m x y^{-1}_m 
\end{align}

\begin{remark} For any element 
	$$ u = u(y_1,\cdots, y_q) = \sum_{i=1}^{p} \alpha_i f_i(y_1,\cdots, y_q) \in kF(Y) $$ 
define its "conjugate" as
$$
u^{*} = u(y_1,\cdots, y_q) = \sum_{i=1}^{p} \bar{\alpha}_i f_i^{-1}(y_1,\cdots, y_q)
$$
For example,
	\begin{align}
	\Psi^{*}_m \equiv  \Psi_m (x^{-1}, Y_m)= 
y_1 x^{-1} y^{-1}_1  + 	y_2 x^{-1} y^{-1}_2  + \cdots
+ y_m x^{-1} y^{-1}_m      \tag{2.1*}
	\end{align}
Clearly, the value of $  u^{*} \in kF $  will be adjoint to the value of $ u $  if both are evaluated on the same group elements in a unitary representation.  
 
\end{remark}
\noindent
\begin{lemma} 
If $ g_1, \cdots g_m \in G $ are pairwise distinct then for any 
	$ a \in  \textnormal{End}(V) $
	\begin{align}
		\Psi_m(a, \rho(g_1), \cdots \rho(g_m) ) = (m/n) \textnormal{tr} (a) I_V  \nonumber
	\end{align} 
	where $ I_V $ denotes a unity matrix in $ \textnormal{End}(V) $. In particular,	
 if $g_0 \in G $ then 
	\begin{align}
		\Psi_m(\rho(g_0), \rho(g_1), \cdots \rho(g_m) ) = (m/n) \chi_{\rho}(g_0) I_V 
	\end{align} 
\end{lemma}
\noindent The proof is obvious. $  \Psi_m(a, \rho(g_1), \cdots \rho(g_m) ) $  commutes with $\rho(G)$ and therefore is a scalar matrix  $ \lambda I_V $. 
 To find $\lambda$, note that
\begin{align}
\textnormal{tr}(\lambda I_V ) = 	n \lambda = \textnormal{tr}(\; 	\Psi_m(a, \rho(g_1), \cdots \rho(g_m) ) \; ) = m\; \textnormal{tr}( a )   \nonumber  
\end{align}
and therefore $ \lambda = (m/n)\textnormal{tr}(a) $ as stated.
\newline
\newline 
Set
\begin{align}
	 \Psi_{m,i} \equiv \Psi_m - (m/n) \chi_i, \; i = 1, \cdots, k 
\end{align}
 and using $ \mathcal{C}_m(Y_m) $ (1.3) as a "guard term"  set  
\begin{align}
\!\!\!\!\!\!\!\!\!\!\!\!				
			\Psi_m(\rho)  \equiv	\Psi_m( \chi_{\rho})  = \;
	 \mathcal{C}_m(Y_m)  
	\Psi_{m,1} v_1 \cdots v_{k-1} \Psi_{m,k} ;
\;	v_1, \cdots, v_{k-1} \in Y
\end{align}

\begin{lemma} \;
	\begin{enumerate}
		\item[(i)] Any irrep of a group of order $ \leq m $   satisfies the identity
		 $$ \Theta_m \equiv \Theta_m( y, Y_m )= \mathcal{C}_m(Y_m) ( \Psi_m y - y 	\Psi_m), \;y \in Y $$  
		\item[(ii)]  The identity $ \Psi_m(\rho) $ (2.4) holds in $ \rho $
		\item[(iii)] 	If identity $ \Psi_m(\rho) $ holds in   $\sigma$ then   
		\begin{enumerate}
			\item[(iii).1] $|H| \leq m $
				\item[(iii).2]
	 if $  	|H| = m$ then
	$	
\textnormal{range} ( \chi_{\sigma}) \subset 
	(\dim\sigma/\dim\rho) \textnormal{range} ( \chi_{\rho})
$
		\end{enumerate}

\item[(iv)] If the order of a group $\Gamma$  is $m$, then identity $ \Psi_m(\rho) $ (2.4)
holds in any irrep $ \tau$ of $ \Gamma$ that satisfies the condition
$	
\textnormal{range} ( \chi_{\tau}) \subset 
(\dim\tau/\dim\rho) \textnormal{range} ( \chi_{\rho})
$ 		
	\end{enumerate}
\end{lemma}
\noindent  The following general remark will be routinely used  below 
\begin{remark}
Representations $ \rho $ and $ \sigma$ are faithful, hence,  as a shortcut we will  identify elements of $ G$ with their $\rho$-images in $ \textnormal{GL}(V() \subset \textnormal{End}(V)  $ and elements of $ H $ with their $\sigma$-images in $ \textnormal{GL}(W) \subset \textnormal{End}(W)  $. In other words,  free variables $ x, y_i, u_{ij}  \cdots \in Y $   take values in  $ \rho(G) \subset \textnormal{GL}(V)  $  or in  $ \sigma(H) \subset \textnormal{GL}(W) $ depending on a context
\end{remark}
\paragraph{Proof of Lemma 2.2}
Evaluating $\Psi_m(\rho) $ (2.4) on a variable value assignment in $\rho(G)  $, we have $\mathcal{C}_m(g_1, \cdots g_m ) = 0 $ unless  $ g_1, \cdots g_m \in G  $ are  pairwise distinct. Assuming the latter, and applying  Lemma 2.1, we see that $ \Psi_{m}(g_0, g_1, \cdots g_m)  $ is a scalar matrix for any $ g_0 \in g $ and that
$ \Psi_{m,i}(g_0,\cdots , g_m) = 0 $ for some $i, \;  1 \leq i \leq k $.  This proves (i) and (ii). 

Suppose that $\sigma$ affords  (2.4).
If $ |H| > m $,  fix $(m+1)$-element subset  $ S =\{ g_1, \cdots , g_{m+1} \} 
$ of $H$. By (2.4),  $ \Psi_m( g_0, h_1, \cdots , h_m )$  (evaluated  in $\sigma(H)$) is  a scalar matrix in $\textnormal{End}(W)  $ 
for any $m$-element subset  $\{ h_1, \cdots , h_m \} $ of $ S$. It is easy to see then, that   
\begin{align}
g^{-1}_i g_0 g_i -  g^{-1}_j g_0 g_j  = \lambda I_W 
\end{align}
 for any  $ g_i, g_j \in S, \; 1 \leq i,\;j \leq m  $. Taking traces on both sides of (2.5) we find that $ \lambda = 0 $ and  
therefore  $ g_i g^{-1}_j $  commutes with $g_0$. Since  subset $S$ is arbitrary, $  g_0 $ must belong to the center of $H$. Since the choice of $ g_0$ is also arbitrary, $ H $ must be commutative and hence cyclic.  Therefore,  $  \dim \sigma =1 $ and 
$ \textnormal{range}(\chi_{\sigma}) $ is the set of all  roots of unity of order $|H|$. 
On the other hand, by (2.3)-(2.4) we must have
$ (m/1)\textnormal{range}(\chi_{\sigma }) \subset (m/n)\textnormal{range}(\chi_{\rho}) $ where in turn 
$  \textnormal{range}(\chi_{\rho}) \subset 
\mathbb{Q}(\xi) $ for some primitive $m$-th root of unity  $\xi$ (cf. e.g. \cite{Serr}, \cite{Serre2}). This, however,  contradicts  $ |H| >  m$ assumption. 
In the remaining case of $ |H| = m, $ it follows from Lemma 2.1 that
$$
 \Psi_m( g_0, h_1, \cdots , h_m ) = (m/\dim \sigma )\chi_{\sigma}(g_0) I_W
 \textnormal{ for any } g_0, h_1, \cdots , h_m \in H 
$$ 
 yielding (iii) when juxtaposed with (2.4). And, finally, the statement (iv) should be now  quite obvious   
\begin{example}
	Looking at  character table of the alternating group $ A_5 $ (cf. e.g. \cite{Serre2}) we see that
	\begin{enumerate}
		\item[(1)] four-dimensional irrep of $ A_5 $ satisfies the identity 
		$$  \mathcal{C}_{ 60} 
		\Psi_{60}v_1
		\left( \Psi_{60} - 60  \right)v_2   (\Psi_{60} - 15) v_{3}  \left( \Psi_{60} +  15 \right)  $$ 
		\item[(2)] five-dimensional irrep of $ A_5 $  satisfies the identity 
		$$  \mathcal{C}_{ 60} 
		\Psi_{60}v_1
		\left( \Psi_{60} - 60  \right)v_2   (\Psi_{60} - 12) v_{3}  \left( \Psi_{60} +  12 \right)  $$ 
		\item[(3)] both  three-dimensional irreps of $ A_5 $ satisfy the identity 
		$$  \mathcal{C}_{ 60} 
		\Psi_{60}v_1
		\left( \Psi_{60} - 60  \right)v_2   \left(\Psi_{60} + 10(1+\sqrt{5}) \right) v_{3}  \left(\Psi_{60} - 10(1+\sqrt{5}) \right)  $$ 
		where $ v_1, v_2, v_3 $ are additional free variables in $ Y $. Note that three-dimwnsional irreps of $ A_5 $ are similar and that non-similar irreps of $A_5 $ are separated by the identities (1)-(3)
	\end{enumerate}   
\end{example}
\begin{example}
	Character values of the standard $n-1$-dimensional  irrep of the symmetric group $ S_n $ are $ n-1, \; n-3, \; \cdots , \; 0, \; -1$ and by Lemma 2.1 it satisfies  the identity (cf. (2.1), (2.4))
	$$  \mathcal{C}_{n!} 
	( \Psi_{n!} - n!  )v_1
		\left( \Psi_{n!} - \frac{n!}{n-1} (n - 3)  \right)v_2 \cdots    (\Psi_{n!} - 0) v_{n-1}  \left( \Psi_{n!} +  \frac{n!}{n-1} \right)  $$ 

\end{example}

\begin{corollary}
	If  $G$ is cyclic and identity $\Psi_m(\rho) $ (2.4) holds in $\sigma$ then either 
	$ |H| < |G| $ or   $ H \approx G $
\end{corollary}
\noindent Proof. By Lemma 2.2 (iii).1 we have $ |H| \leq |G| $. Assuming  that $ |H| = |G| $ and applying Lemma 2.2 (iii).2 we get 
 $\; \chi_{\sigma}(h)  \chi_{\sigma}(h)^{*} = (\dim \sigma)^2$ 
for all $ h \in H $. Since the character 
$ \chi_{\sigma} $ is irreducible this leads to
\begin{align}
	|H| = \sum_{h \in H }  \chi_{\sigma}(h)\chi_{\sigma}(h)^{*} \; =
	 \;  |H|(\dim\sigma)^2 \nonumber
\end{align}
Thus $ \dim\sigma = 1$, $ H $ is cyclic, $ \textnormal{range}(\chi_{\sigma })\subset \textnormal{range}(\chi_{\rho})$ and therefore  $ H \approx G $  
\begin{remark} The irrep  $ \rho$ satisfies identity $	
\mathcal{C}_m(Y_m)  
\Psi_{m,1} \cdots  \Psi_{m,k} \equiv \Psi'_m(\rho)  	$ that is obtained from $ 	\Psi_m(\rho)  $ by substitution $ v_i \rightarrow 1, \; i = 1, \cdots k - 1 $. It is easy to see that $ \Psi_m(\rho) $ and $ \Psi'_m(\rho) $
are equivalent modulo  identity $ \Theta_m $  and the conclusion of  the statement (ii) of Lemma 2 remains valid if one requires that $ \sigma $ satisfies  
both identities $  \Theta_m $ and   $	\Psi'_m(\rho) $ instead of just one identity $ \Psi_m(\rho) $. Note that the identity $ \Theta_m $ guarantees that $ \Psi_m$ commutes with the representation. This remark is applicable to some other statements below (sections 2.2, 3)      
\end{remark}
\noindent It is easy to see that the identity $\Psi_m(\rho)$ is  invariant under the action of the Galois group  $ \textnormal{Gal}(\mathbb{Q}(\sqrt[m]{1})/\mathbb{Q}) $ on character values.       
We interrupt the main line of this discussion in order to review some basic facts related to a notion of \textit{Galois conjugacy}  of representations (cf. e.g. \cite{Serr}, \; \cite{Serre2}, \; \cite{conj}) as it turns out to be highly relevant to the subject at hand
\subsection{Galois Conjugate Representations}
\noindent As was mentioned already, similar representations have the same identities. In particular, since the identity (2.4) depends only on the set of values of the character $ \chi_{\rho} $, one has 
$ \Psi_m(\chi_{\rho}) \equiv \Psi_m( \chi_{\rho \alpha} ) $ for any automorphism $ \alpha $ of $G$ (cf. (2.4)). 
Another operation on characters that preserves the set of character values is Galois conjugation.

Let $ \Gamma $ be any group   and let $ \theta  : \Gamma \rightarrow \textnormal{GL}_n(V)$ be a representation of $\Gamma$ in a complex vector space of dimension $ n $. Let $ \epsilon $ be an automorphism of the field of complex numbers. Field automorphism $ \epsilon $ acts on a group of invertible matrices and composing $ \epsilon $ with $ \theta$ 
one gets a representation $\epsilon \theta : \Gamma \rightarrow \textnormal{GL}(V) $ (cf. e.g. \cite{Serre2}). Hence there is a (left) action of field automorphisms on representations.  On the other hand (see Introduction), there is another  (right) action  of the automorphism group $ \textnormal{Aut}(\Gamma) $ on representations  of $ \Gamma $. By definition these two actions commute: if  $\alpha \in  \textnormal{Aut}(\Gamma) $ then $ \epsilon(\theta \alpha) = (\epsilon \theta) \alpha $.

If the group $ \Gamma $ is finite of order  $|\Gamma| = m $ then $\rho $ is realizable (conjugate in \textnormal{GL}(V)) to a representation over $ \mathbb{Q}(\omega_m) $ where $\omega_m $ is a  primitive root of unity of order $m$ (cf. e.g. \cite{Serr}).
The cyclotomic field $ \mathbb{Q}(\omega_m) $ is invariant under automorphisms of $ \mathbb{C}$,
 hence one can
assume that $ \epsilon \in \mathfrak{G}(\Gamma)  \equiv  \textnormal{Gal}(\mathbb{Q}(\omega_m)/\mathbb{Q}) $.

\begin{definition} (cf. e.g. \cite{Serre2}).
	A representation $ \theta' $ is said to be Galois-conjugate to a finite-dimensional representation $\theta$ of a finite group $ \Gamma $ if there is $ \epsilon \in  \mathfrak{G}(\Gamma) $ such that $ \theta'$ is equivalent to $ \epsilon \theta$  
\end{definition}
\noindent  Both representations $ \theta$ and its Galois-conjugate representation $ \epsilon \theta $ are realizable over $ \mathbb{Q}(\omega_m) $ and since traces of Galois-conjugate matrices are Galois-conjugate the function $ \epsilon \chi_{\theta} $ is a character of the representation $ \epsilon \theta $. 
\begin{remark} Galois group $\mathfrak{G}(\Gamma) $ is abelian.
Hence (see \cite{conj}),	 if a character $ \chi $ is Galois-conjugate to $ \chi_{\theta} $ then, in fact, there is $ \epsilon' \in  \textnormal{Gal}(\mathbb{Q}(\chi_{\theta} )/\mathbb{Q}) $ such that  $ \chi = \epsilon' \chi_{\theta} $  , where  
 $ \mathbb{Q}(\phi) \subset \mathbb{Q}(\omega_m) $ denotes the field generated by values of a $\Gamma$-character $\phi $. We see that automorphisms of the field of complex numbers act on  representation $ \theta$ as elements of $  \textnormal{Gal}(\mathbb{Q}(\chi_{\theta})/ \mathbb{Q}) $  	
\end{remark}
\noindent   As was mentioned above we have by definition 
\begin{lemma} Let $ \theta $ be a representation of a finite group $\Gamma$. If $ \epsilon \in \mathfrak{G}(\Gamma) $ and $ \alpha \in \textnormal{Aut}(\Gamma) $ then
	representations $ (\epsilon \theta ) \alpha $ and  
	$ \epsilon (\theta \alpha) $ are equivalent.
	In other words,  we have
	\begin{equation}
		\epsilon (\chi_{\theta} \alpha) = 	(\epsilon \chi_{\theta}) \alpha
		\nonumber
	\end{equation} 
\end{lemma}
\noindent  to check this directly, write
$$ 	\epsilon (\chi_{\theta} \alpha)(h) = \epsilon (\; Tr( \; \theta(\alpha(h) \; )\;)  =   \; Tr( \epsilon \; \theta(\alpha(h) \; ) =  \chi_{\epsilon \theta}( \alpha(h)) = (\epsilon \chi_{\theta}) \alpha(h)
 $$
 for any $ h \in \Gamma $ and $ \epsilon \in \mathfrak{G}(\Gamma) $
\begin{corollary}
	Let 
	$  \rho $ be a finite dimensional representation of a group $ \Gamma $ and let 
	$ \alpha  $ be an automorphism of $\Gamma$. If $ \rho $ is Galois conjugate to another representation $ \sigma$ then 
	representations $ \rho \alpha $ and $ \sigma \alpha $ are also Galois conjugate. 
	
\end{corollary}
\noindent

\begin{example}
	Suppose that representation $\theta$ is exact and that there is $ g \in  \textnormal{GL}(V) $ such that  
	$  \epsilon \theta(\Gamma) = g^{-1} \theta(\Gamma) g   $. Then by definition (cf. Remark 2.1) representations $ \theta$ and its Galois-conjugate representation $ \epsilon \theta $ are similar. Two  $3$-dimensional representations of alternating group $A_5$ (Example 2.2) are Gallois conjugate and similar (cf.  e.g. \cite{Serre2}) 
\end{example}
\begin{example}
	The symmetric group $S_6$ has two similar non-equivalent  irreps of dimension $5$ (cf. e.g. \cite{Wildon_L}, \cite{S6}). These irreps are not Galois-conjugate simply because all representations of a symmetric group are defined over the field of rational numbers. 
\end{example}
\begin{example}
	The sporadic simple Mathieu group $M_{11}$ has two irreps of 
	order $16$ that are Galois conjugate but not similar (cf. \cite{M11}, \cite{Hughes}).  Some other examples of Galois-conjugate non-similar irreps can be found in \cite{Marin} (see also section 3.1.2). By Theorem 1.1, in all these cases there are identities that distinguish between irreps in question. Writing these identities explicitly seems like an interesting problem (cf. Proposition 3.3)
\end{example}
\noindent Any automorphism $\phi$ of 
the field $ k \equiv \mathbb{C} $ over the field $ \mathbb{Q} $ can be extended to an automorphism $  kF(Y) \rightarrow kF(Y) $ by a rule
$$ \sum_{i=1}^k a_i f_i  \rightarrow  \sum_{i=1}^k \phi(a_i) f_i, \; \alpha_i\in k, \; f_i \in F(Y) $$
and again by definition, we have. 

\begin{lemma}
	Let  $ u \in kF $ be an identity of representation $ \theta$.
	\begin{enumerate}
		\item [(i)] $ \epsilon(u) $ is an identity of the representation $ \epsilon \theta$ for any 
		$ \epsilon \in  \mathfrak{G}(\Gamma) $
		\item[(ii)] in particular any identity of $\theta$ with rational coefficients holds in $\epsilon \theta$
		\item[(iii)] Suppose that a representation 
		$\theta$ is defined over a subfield $ K \subset \mathbb{C}$ Then any identity of $\theta$ is a linear combination of identities with coefficients in $K$
	\end{enumerate}
\end{lemma}
\begin{corollary}
	All identities of any representation of a symmetric group follow from identities with rational coefficients 
\end{corollary}
\noindent Finally, as was mentioned above, we have 
\begin{lemma} Let $ \theta  : \Gamma \rightarrow \textnormal{GL}_n(V)$ be an irrep of a finite group $\Gamma, \; (|\Gamma| = m ) $. Then $ \Psi_m(\theta) \equiv \Psi_m( \epsilon \theta ) $ for any  $ \epsilon \in \mathfrak{G}(\Gamma) $. In other words, identity (2.4) is preserved by the action of  $\mathfrak{G}(\Gamma)$ on $kF(Y)$.  	  
\end{lemma}
\noindent Proof. It is well known that $\textnormal{range}(\chi_{\theta}) = \textnormal{range}(\epsilon\chi_{\theta}) $. The proof of this last statement boils down to a combination of the following well known facts (cf. e.g. \cite{Serre2}) that we summarize as follows.
\noindent 
\begin{lemma}\
\begin{enumerate}
	\item [1)]  the group  $ \mathfrak{G}(\Gamma) $
	is isomorphic to $ \mathbb{Z}_m^{*} $
	\item[2)] let $ \epsilon_t \in \mathfrak{G}(\Gamma) \approx \mathbb{Z}_m^{\ast} $ be a Galois automorphism that corresponds to an integer $t $ coprime to $ m $.
	Then $ \epsilon_t(\theta)(x) $ and $ \theta(x^t) $ have the same spectrum for all $ x \in \Gamma $ 

	\item [3)] therefore, for an integer $t $ coprime to $ m $, the action of a corresponding  Galois automorphism $ \epsilon_t \in \mathfrak{G}(\Gamma) \approx \mathbb{Z}_m^{\ast} $ on a character $ \chi $ of $  \Gamma $ is given by $ \epsilon_t(\chi)(x) = \chi(x^t) , \; x \in \Gamma $  
	\item [4)]for any integer $t $ coprime to $ m $ the map $ x \rightarrow x^t , \; x \in \Gamma $ is a bijection 

\end{enumerate}
\end{lemma}
\begin{conclusion}
 Identity (2.4) cannot be used to distinguish between Galois conjugate irreps
\end{conclusion} 
\noindent We end this section with another simple property of Galois-conjugate representations  

\begin{lemma}
	Let $ \Gamma $ be a finite group that contains a subgroup $ \Sigma $. If characters $ \phi, \psi$ of  $\Sigma $ are Galois conjugate then
	induced characters  $ \textnormal{Ind}^{\Gamma}_{\Sigma}\phi , \; \textnormal{Ind}^{\Gamma}_{\Sigma}\psi $  are also Galois conjugate	
\end{lemma}

\subsection{Dimension Identity}

Let  $  Y_m = \{y_1, y_2 , \cdots ,y_{m} \}  \subset Y $ and  $ X_m = \{ x_1, x_2 , \cdots ,x_{m} \} \subset Y  $ be two sets of independent variables.
Let
\begin{align}
	\mathcal{D}_{n=\dim\rho}(X_m,Y_m) = \mathcal{C}_m (X_m) \mathcal{C}_m (Y_m) 
\left(\sum_{i=1}^{m} \Psi_m( x_i, Y_m ) x_i^{-1} - (m/n)^2 \right)   \nonumber \\
\equiv \; \mathcal{C}_m (X_m) \mathcal{C}_m (Y_m)  	
\left( \sum_{i,j=1}^{m}[y_j,x_i^{-1}] - (m/n)^2 \right)
\end{align}

\begin{lemma} 
	The identity (2.6)
	holds in $ \rho$
\end{lemma}  
\noindent Proof.
Let $ x_i \rightarrow h_i \in G $ and $ y_i \rightarrow g_i \in G ,\; i = 1, \cdots, m  $ be any assignment of pairwise distinct variables. It is well known and easy to check (cf. Lemma 2.1) that $\sum_{i=1}^{m} \chi_{\rho}(h_i)h_i^{-1}  $ is a scalar matrix $  (m/n)I_V $. Hence, by Lemma 2.1
\begin{align}
 \sum_{i=1}^{m} \Psi_m( h_i, g_1, \cdots, g_m  ) h^{-i}_i  
	 = (m/n) \sum_{i=1}^{m} \chi_{\rho}(h_i)h_i^{-1}   
	= (m^2/n^2)I_V
\end{align}
The following corollary is  well known (cf. \cite{Serre2}, 8.13, exercise 27)   
\begin{corollary} (cf. \cite{Serre2}).
	The average of pairwise group commutators in any faithful  $n$-dimensional irrep of a
	finite group is equal to $(1/n)^2I$.  
\end{corollary}
\noindent In this form the statement is true for not necessarily finite compact groups.
\begin{corollary} The expectation of a commutator in a faithful $n$-dimensional irrep of a compact group is equal to $	(1/n^2) I	$ and therefore, expectation of the trace of  a commutator is  $	1/n 	$    
\end{corollary}
\noindent Indeed, 	
	let $\chi_{\tau}$ be a character of a faithful unitary irrep  $ \tau $ of a
	compact group $G$ and let $\mu$ be a normalized Haar (probability)  measure on $G$. 
	Using the same arguments as above, we get
\begin{align}
	\underset{y \in G}{\int} \tau(y)^{-1}\tau(x)^{-1} \tau(y) d\mu(y)  =  \frac{1}{\dim \tau }\chi_{\tau}(x^{-1}) \nonumber
\end{align}	
and, therefore,	
	\begin{align}
		\underset{x,y \in G}{\iint} \;	
		 [\tau(x),\tau(y)] d\mu(x)d\mu(y) =
		\underset{x \in G}{\int} \left(	\underset{y \in G}{\int} \tau(y)^{-1}\tau(x)^{-1} \tau(y) d\mu(y) \right) \tau(x) d\mu(x) =     \nonumber \\
		\!\!\!\!\!\!\!\!\!\!\!\! =	\; \frac{1}{\dim \tau }\underset{x \in G}{\int} \chi_{\tau}(x^{-1})\tau(x) d\mu(x) = 
		\frac{1}{(\dim\tau)^2} I	 \nonumber
	\end{align}  

\begin{remark}
	 Define an expectation of $ u = u (y_1, , \cdots y_t ) \in kF(Y) $ over a representation $\tau$ of  a compact group $G$  as 
	\begin{align}
	\mathbb{E}(\tau)(u)	\;\; = \underset{(g_1, \cdots, g_t)  \in G \times \cdots \times G }\int u( \tau(g_1), \cdots, \tau(g_t) ) d\mu	\nonumber
	\end{align} 
where $\mu$ is a normalized  Haar measure on a  Cartesian product $ G ^t $ of $t$ copies of $ G $. It is clear that the expectation 
$	\mathbb{E}(\tau)(u)$ is a scalar matrix when $\tau$ is irreducible. Note that $ u $ is an identity of $ \tau $ if and only if $	\mathbb{E}(\tau)(uu^{*}) = 0 $. Indeed, suppose that $	\mathbb{E}(\tau)(uu^{*}) = 0 $.  Then $ \mathbb{E}(\tau)(tr(uu^{*})) = 0 $ and since the matrix $uu^{*} $ is symmetric positive semi-definite, $ tr( u( \tau(g_1), \cdots, \tau(g_t)  u( \tau(g_1), \cdots, \tau(g_t)^{*}) $ must be zero for any $ g_1, \cdots g_t \in G $
\end{remark}


\begin{lemma}
	If dimension  identity (2.6) is satisfied by representation  $ \sigma$ then either
	\begin{itemize}	
		\item[(i)] $ \dim \rho  = \dim \sigma  = 1 $ or
		\item[(ii)] 	  $ |H| \leq m $ and  if  $ H = m $ then $\dim \sigma = n = \dim \rho $
	\end{itemize} 
\end{lemma}
\noindent Proof. We work with variable value assignments in $ \sigma(H)$. By (2.7)
\begin{align}
	\sum_{i=1}^{m} \Psi_m( h_i, g_1, \cdots, g_m  ) h^{-i}_i   \; = \; (m^2/n^2)I_W
\end{align}
for pairwise distinct $ g_i \in H $ and  $ h_i \in H, \; i = 1, \cdots, m $. If $ |H| >m $, fix 
the subset  
 $ \{ g_1, \cdots  g_m \}  \subset H $ and  
 evaluate left hand side of (2.8) on subsets
  $ \{ 1, \; h_1, \cdots, h_{m-1} \} \subset H  $ and 
$ \{ h_1,  \cdots , h_{m-1},  h_m \neq 1 \} $. Comparing the results we get
$
mh_m = \Psi_m( h_m, g_1, \cdots, g_m  )  
$ and since $h_m$ is arbitrary, 
$
\Psi_m( h, g_1, \cdots, g_m) = m h 
$
for any $ h \in H $. Moreover, since
$ \Psi_m( h, g_1, \cdots, g_m) $  does not depend on the subset $ \{g_1, \cdots g_m\} \subset H $, we see 
as in the proof of Lemma 2.2, that $H$ is cyclic and that $ \dim \sigma = 1$.  For one-dimensional representation, however, the equation (2.8) turns into
$ m^2 = m^2/n^2 $ and  we find that $ n = 1 $ as stated by (i). Finally,  if  $ |H| = m$ then 
taking traces on both sides of  (2.8) and using Lemma 2.1 we get 
\begin{align}
	(m/\dim\sigma) \sum_{h \in H} \chi_{\sigma}(h)\chi_{\sigma}(h^{-1}) = 
  \frac{ m^2}{n^2} 	\dim \sigma  \nonumber
\end{align}
that is equivalent to $ \dim \sigma  = n $.
\begin{corollary}
	Faithful irrep $\theta$ of a group of order $m$ satisfies identity (2.6) if and only if $\dim\theta =n \;(\equiv \dim \rho) $ 
\end{corollary}
\noindent  Combining lemmas  2.8, 2.2 and Corollary 2.1 we get
\begin{proposition}
	If identities (2.6) and (2.4) hold in $ \sigma$ then  $ |H| \leq m $. If  $|H| = m $ then $  \dim \sigma = \dim \rho \; (\equiv n)$ and $
	\textnormal{range} ( \chi_{\sigma})
	\subset  	\textnormal{range} ( \chi_{\rho}) $
\end{proposition}

\noindent  Here is another variant of  dimension identity

\begin{example}
	The following identity holds in $\rho$ 
	\begin{align}
		\mathcal{C}_m (X_m) \mathcal{C}_m (Y_m) \left(
		\sum_{i=1}^{m} \Psi_m( x_i, Y_m )  \Psi_m( x_i^{-1}, Y_m )    -  m^3/n^2 \right)  \nonumber 
	\end{align}	

\end{example}

\subsection{Character Range Identities}
Fix an additional set $ V_m = \{v_1, \cdots, v_m \}\subset Y, \; ( m = |G| )   $ of free variables. Speaking informally, one can say that the identity (2.4) restricts the range of a character. Somewhat similar identity can be used to make sure that a character range contains a given value

\begin{lemma}	 
	For any $ \xi \in 	\textnormal{range} ( \chi_{\rho} ) $, the  identity 
	\begin{align}
	 \mathcal{R}_{\xi}	\equiv \mathcal{R}_{\xi}(\rho) = \mathcal{C}_m (X_m) \mathcal{C}_m (Y_m) 
		\prod_{ i=1}^m \left(  \Psi_m( x_i, Y_m )  - (m/n)\xi  \; \right) v_i
	\end{align}
holds in $\rho $. 
 If  $ |H| =m $ and identity $\mathcal{R}_{\lambda}(\rho)   $ holds in $\sigma$ then 
	 $ \lambda \in 	(\dim \rho \;/\dim\sigma)\textnormal{range} ( \chi_{\sigma} ) $

\end{lemma}
\begin{remark}
	It is quite obvious  that identity $  \mathcal{R}_{\epsilon \xi} $ holds in Galois conjugate irrep $ \epsilon \rho $ for any 
	$ \epsilon \in \mathfrak{G}(G) $.  On the other hand
	the identity $  \mathcal{R}_{\epsilon \xi} $ holds in $ \rho$ as well, since  
	$ \epsilon \xi \in \textnormal{range}(\chi_{\rho}) $ (cf. Section 2.1).    
\end{remark}

\begin{lemma}	 
 If the identity $	\Psi_m(\rho)$ (2.4) and all identities $\mathcal{R}_{\xi} , \; \xi \in 	\textnormal{range} ( \chi_{\rho} ) $ hold in $\sigma$ then 
	\begin{enumerate}
		\item[(ii).1] $ |H| \leq m $
		\item[(ii).2] $ \textnormal{range} ( \chi_{\sigma} ) = 	(\dim \sigma/\dim\rho) \textnormal{range} ( \chi_{\rho})	 $  if $ |H| = m $
	\end{enumerate}
\end{lemma}
\noindent The statement (ii).2 probably requires some explanation. If $|H| = m$ then using lemmas 2.2 and 2.9 we have 
$$  \textnormal{range} ( \chi_{\sigma} ) \subset 	(\dim \sigma/\dim\rho) \textnormal{range} ( \chi_{\rho}) 
\subset	(\dim \sigma/\dim\rho)	(\dim \rho/\dim\sigma) \textnormal{range} ( \chi_{\sigma} )
= \textnormal{range} ( \chi_{\sigma} )$$ 
\begin{question}
	Let $ \chi_1, \chi_2$ be faithful  irreducible characters of a finite group. If $ \chi_2 = r\chi_1 $ for some real number $ r > 0 $   then,  $ \chi_1  = \chi_2$ and $ r = 1 $.
	Is it true that $ \textnormal{range} ( \chi_1) = r \cdot \textnormal{range} ( \chi_2) $ implies $ r = 1 $ ?
	
\end{question}
\noindent Combining lemmas 2.2, 2.8 and 2.10 we get
\begin{proposition}	If irrep  $ \sigma$ satisfies character identity $\Psi_m(\rho)$ (2.4),   dimension identity $	\mathcal{D}_{n=\dim\rho} $ (2.6) and character range identities $ 	\mathcal{R}_{\xi}(\rho) $ (2.9)  for all $ \xi \in \textnormal{range} ( \chi_{\rho}) $ then
 either $ |H| \leq m $ or  $|H| = m $,  $  \dim\sigma = \dim \rho   $ and $
	\textnormal{range} ( \chi_{\sigma})
=	\textnormal{range} ( \chi_{\rho}) $
\end{proposition}
\begin{corollary}
	Let $ |H| = m \; (= |G|) $. Suppose  that the irrep $\sigma$   satisfies identities $ \mathfrak{D}(G)$ (cf. Lemma 1.1) in addition to the list of identities in Proposition 2.2.  Then 
	$$ H \approx G, \;  \dim \sigma= \dim \rho   \textnormal{ and } 
	\textnormal{range} ( \chi_{\sigma})
	=	\textnormal{range} ( \chi_{\rho}) $$
\end{corollary}
\begin{corollary}
	Let $ \phi $ and $ \psi  $ be  faithful irreps of  a finite group $\Gamma$. Identities  $ \Psi_{|\Gamma|}(\phi), \; 	\mathcal{D}_{\dim \phi} $ and   
	$ 	\mathcal{R}_{\xi}(\phi) $ for all $ \xi \in  \textnormal{range} ( \chi_{\phi} ) $ are satisfied by irrep $ \psi $ if and only if 
	$   \dim \psi = \dim \phi $ and $
	\textnormal{range} ( \chi_{\psi})	=	\textnormal{range} ( \chi_{\phi})
 $
\end{corollary}
\noindent In other words,
\begin{enumerate}
\item[] acting group, dimension and 
character values of an exact irrep of a finite group of a given order are determined by the finite set of identities (1.2), (2.4), (2.6) and (2.9) 
\item[] for a fixed acting group, dimension and 
character values of an exact irrep are determined by the set of identities (2.4), (2.6) and (2.9)
\item[] however, the set of identities (2.4), (2.6) and (2.9) do not distinguish between Galois conjugate irreps (cf. Section 2.1)
\end{enumerate}

\begin{remark}
It is not true that faithful irreps with the same character ranges are similar (see e.g. Examples 2.7, 3.1 below) as there are numerous examples of non-isomorphic groups with the same character tables (see for example \cite{Dade}-\cite{Meir}, \cite{Seg}-\cite{Davydov}).  It seems, therefore, that looking for  identities that distinguish between irreps does make sense 
\end{remark}

\begin{example}
	 The symmetric group $ S_5$ has two non-similar four-dimensional irreps that have different character value sets (cf. \cite{S5}).
	 Therefore, by Theorem 1.1, there is an identity that holds in one of these irreps but not in the other, and Corollary 2.6 can be used to write down the separating identity explicitly.   On the other hand,  character ranges of five-dimensional non-similar representations of $S_5$ are equal and therefore,  Proposition 2.2 cannot be used to distinguish between them (cf. Remark 2.7 below) 
\end{example}

\noindent We now turn to identities that  distinguish between character values and corresponding  conjugate class sizes.

\noindent
\subsubsection{Character Level Sets }

Define \emph{range signature }
of a character $ \chi $ as the set of number pairs $$ R(\chi) = \{ (\lambda, \; | \chi^{-1}(\lambda) | )\}, \;   \lambda \in  \textnormal{range}(\chi) \} $$    
For any $ i = 1, \cdots, k  = |  \textnormal{range}(\chi_{\rho}) | $ set
 $ t_i = | \chi^{-1}_{\rho}(\chi_i) | $ and 
\begin{align}
	\!\!\!\!\!\!\!\!\!\!\!\			
	\mathcal{S}_i  =
	\prod_{S \subset X_m,\; |S|= t_i }   \left( \sum_{x \in S} ( \Psi_m( x, Y_m ) -  \frac{m}{n}\chi_i      )
	(   \Psi_m( x^{-1}, Y_m ) -  \frac{m}{n}\bar{\chi}_i      ) \right) v_S 
\end{align}
where $ v_S \in Y $ is a set of free variables in $ Y $ indexed by subsets of $ X_m $ and, as a reminder, $m= |G|, \; n = \dim \rho, \; 
\{ \chi_1, \cdots, \chi_k \} = \textnormal{range}(\chi_{\rho}) $. 
Let also
\begin{align}
	\!\!\!\!\!\!\!\!\!\!\!\	\mathcal{R}_i( \chi_{\rho} ) \equiv	\mathcal{R}_i(\rho) \equiv	\mathcal{R}_i(\rho, \chi_i) = \; \mathcal{C}_m (X_m) \mathcal{C}_m (Y_m) \mathcal{S}_i, \; i = 1, \cdots, k 
\end{align} The intuition behind     (2.11) can be described as follows. First of all,  we have  an obvious
\begin{lemma} Let $ \Gamma $ be a group of order $m$. If $ h \in \Gamma, \; g = \{g_1, \cdots g_m \} \subset \Gamma  $ are evaluated in a faithful $n$-dimensional unitary irrep of $ \Gamma $ (cf. remarks 2.1, 2.2), then 
\begin{equation*} 
( \Psi_m( h, g ) -  \frac{m}{n}\lambda_i I     )
	(   \Psi_m( h^{-1}, g ) -  \frac{m}{n}\bar{\lambda}_i I  ) \; = \; \parallel 
	\Psi_m( h, g ) -  \frac{m}{n}\lambda_i \parallel^2 I  \nonumber
\end{equation*}
\end{lemma} 
\noindent We can say informally, that internal sums in (2.11)  are taken over  square distances of generic (scaled) character values from a fixed character value in the range.     
These sums are then multiplied over all $ t_i$-size subsets of the variable set $ X_m $  
\begin{lemma}\

 \begin{enumerate}
\item[(i)]  Assume that  $ |H| = m $. If an  identity 	$\mathcal{R}_i(\rho) \equiv \mathcal{R}_i(\rho, \chi_i)$ (2.11) holds in $ \sigma$ for some $ 1 \leq i \leq k$ 
then  $ | \chi^{-1}_{\sigma}( \frac{\dim\sigma}{\dim\rho}\chi_i) |  \geq t_i$
\item[(ii)] If $ H \geq m $ and all the identities $\mathcal{R}_i(\rho, \chi_i), \; i =1, \cdots , k $ hold in $\sigma $ then the identity $\Psi_m(\rho) $ (2.4) holds in $ \sigma$   
 \end{enumerate}

\end{lemma} 
\noindent Proof. We begin with the statement (i). Since $|H | = m $, all the   $ \Psi$-terms of the multiplier terms $ \mathcal{S}_i$ of $\mathcal{R}_i(\rho, \chi_i)$ will evaluate to $ (m/\dim \sigma) $ times an appropriate value of $ \chi_{\sigma} $, when computed in $ \sigma$. Assuming that $ \mathcal{R}_i(\chi_{\rho}) $ holds in $ \sigma$ we observe that an internal sum of (2.10) must vanish on at least one of  $t_i$-size subsets of $ H $ and  (i)  follows from Lemma 2.13. It is easy to see that under conditions of statement (ii) all possible values of terms $ \Psi_m(x,Y_m) $ belong to the set   
$  (\dim\sigma / \dim\rho) \textnormal{range}(\chi_{\rho}) $ and, therefore, the statement (ii) follows from the proof of  Lemma 2.2 

\begin{theorem}
	Identities 	$\mathcal{R}_i(\rho) $ (2.11) hold in irrep $ \rho$ for all $i = 1, \cdots, k $. If  all identities $\mathcal{R}_i(\rho)  $ hold in $\sigma$ then
	\begin{enumerate}
		\item [(i)] $ |H| \leq m $ 
			\item[(ii)] if  $ |H| = m $ then $ \dim \sigma= \dim \rho  $ and characters
		 $ \chi_{\rho} $ and $ \chi_{\sigma}$ have the same range signatures
	\end{enumerate}
\end{theorem}    
\noindent Proof. 
Evaluating 	$ \mathcal{R}_i(\rho)$ (2.11) in irrep $ \rho$ on any elements of $ G $, we note that  one of the factors  in  (2.10)  will match the value of $\chi_i$  in the range signature  $ R(\chi_{\rho}) $ and therefore will vanish (in  $ \text{End}(V)$) by Lemma 2.1. Hence all the identities (2.11) hold in $ \rho$. By Lemma 2.12 (ii),  $	\Psi_m(\rho)$ (2.4) holds in $\sigma$ and statement (i)  follows from Lemma 2.2. Assuming, therefore, that  $ |H| =  m $ and using Lemma  2.11 (i) and Lemma 2.12 we find that
$$   \left \vert \chi^{-1}_{\sigma}\left( \frac{\dim\sigma}{\dim\rho}\chi_i \right) \right \vert  = t_i, \; i = 1, \cdots, k  $$ meaning that sets 
$ \chi_{\sigma}^{-1}(\; (\dim \sigma / n ) \chi_i\;), \; i = 1, \; \cdots, k $ form a partition of $ H $.    
Now, to verify the statement (ii), one can  proceed as follows: 
\begin{align}
	|G| = |H| = \sum_{h \in H }  \chi_{\sigma}(h)\chi_{\sigma}(h)^{*} \; = \;	
	\sum_{i=1}^{k} \; \sum_{\chi_{\sigma}(h) = (\dim \sigma /n)\chi_i } \chi_{\sigma}(h)\chi_{\sigma}(h)^{*} \; = \;  \nonumber \\
	\; = \; (\dim\sigma/n)^2	\sum_{i=1}^{k}\sum_{g \in G, \;\chi(g) = \chi_i } \chi_{\rho}(g)\chi_{\rho}(g)^{*} \; = \;  (\dim\sigma/n)^2 |G|  \nonumber
\end{align}
Hence, $ \dim \sigma = n $ and   $     | \chi^{-1}_{\rho}(\chi_i)|  =  | \chi_{\sigma}^{-1}(\; \chi_i\;)| $ for all $ \; i = 1, \cdots k $
\begin{corollary}
	Let $ |H| = m \; (= |G|) $. If in addition to identities $\mathcal{R}_i(\rho), \; i = 1, \cdots, k $ the irrep $ \sigma$ satisfies identities  $ \mathfrak{D}(G)$ (see Lemma 1.1) then 
	$ H \approx G, \;  \dim \sigma= \dim \rho \equiv n $ and characters
	$ \chi_{\rho} $ and $ \chi_{\sigma}$ have the same range signatures
\end{corollary}

\begin{remark} 
Some of the non-similar same-dimension irreducible characters of the symmetric group $ S_5 $ have identical range signatures (cf. \cite{S5}). The same is true for $S_6$ (cf. \cite{S5}, \cite{S6})
 and therefore identities specified by Theorem 2.1 can not  be used to distinguish between irreps of $S_5$ or $ S_6$  
\end{remark}

\begin{conjecture}
Same-dimension irreducible characters   of a symmetric group have identical
 character range signatures if and only if  they have identical ranges 
\end{conjecture}

\subsubsection{Conjugate Class Identities}

Character level sets are unions of conjugate classes. We will construct an identity that takes into account not only character level set sizes (range signatures) but the sizes of conjugate classes themselves. Let 
\begin{align}
		 C_1, \; C_2 , \cdots, C_s = \{1\}
\end{align}
be a list of all the conjugate classes of $ G $ 
ordered by size, so that 
\begin{align}
	t_1 =| C_1 | \; \geq \; t_2 = |C_2| \; \geq \cdots \; \geq \; t_s = |C_s| = 1   
\end{align}

\noindent Take also a similar list of conjugate classes of $ H $
\begin{align}
		t'_1 =| C'_1 | \; \geq \; t'_2 = |C'_2| \; \geq \cdots \; \geq \; t'_{s'} = |C'_{s'}| = 1     
\end{align}
and let  $ \chi_{\rho,i}, \;i = 1, \cdots, s $ and
$ \chi_{\sigma,j}, \; j = 1, \cdots , s'$ be corresponding character values. 
\noindent Take $s$  disjoint sets of free variable  $Y_i \subset Y, \; |Y_i| = t_i, \; i = 1, \cdots, s$ and another
 set of variables $x_1, \cdots, x_s \in Y $. 

Using sufficient supply of   free variables $  v_i, u_i, u_{ij} \cdots \in Y $,  set
\begin{align}
\!\!\!	\mathcal{S}_r(x_r,Y_r) =   \! \prod_{ \; i < j, \; y_i \in Y_r }^{t_r} ( 1 - [x_r, y_i y^{-1}_j] ) u_{ij}, \; r =1, \cdots, s-1; \;	\mathcal{S}_s = 1  
\end{align}

\begin{align}	
	\mathcal{C}_{pq} = \mathcal{C}_{pq}(x_p, x_q; Y_p) \; = \;  \!\!\! \prod_{  y \in Y_p} ( x_q - y ^{-1} x_p y ) v_i, \; p,q = 1, \cdots, s; \; q > p   \nonumber
\end{align}
\begin{align}
\!\!\!\!\!\!	\mathcal{D}_{1} = 1, \;		\mathcal{D}_{q} = \prod_{r=1}^q \mathcal{S}_r(x_r,Y_r) \prod_{1\leq a< b \leq q}^q	\mathcal{C}_{ab}, \; q = 2, \cdots , s  \nonumber
\end{align}
\begin{align}
	\mathcal{E}_{ab} \equiv 	\mathcal{E}_{ab}(x_a, Y_a ) =   (\Psi_{t_a}(x_a,Y_a) - (t_b/n) \chi_{b,\rho} )  	, \; a,b = 1, \cdots, s  \nonumber
\end{align}
Let $ c_1 > c_2 > \cdots > c_l $ be downward ordered  sizes of conjugate classes in the list (2.12). For any $ 1 \leq i \leq l $, let $A_i$ be a set of indices of  conjugate classes
of the same size $c_i$ in the list (2.13) and let 
$ S(A_i)$ be a full group of permutations of the set $ A_i$. Set 
\begin{align}
	\mathfrak{E}_i(\nu) = \sum_{a \in A_i} 	\mathcal{E}_{a \nu(a)} \mathcal{E}_{a \nu(a)}^{*}, \; \nu \in  S(A_i) 
\end{align}
and    
\begin{align}
	\mathcal{E}_{i} = \prod_{\nu \in S(A_i) } 		\mathfrak{E}_i(\nu), \; i = 1,\cdots, l   \nonumber
\end{align}
Finally define the identity 
\begin{align}
	 \mathcal{L}(\rho) \equiv \mathcal{L}(\chi_{\rho}) =  	\mathcal{D}_{s} \sum_{i=1}^l 	\mathcal{E}_{i} 
\end{align}
\begin{remark}
	The  expressions somewhat similar to  (2.15) are well known as "principal centralizer identities" (cf.e.g. \cite{P})  
\end{remark}
\begin{theorem}
	\noindent 				 
	\begin{enumerate}
		\item[(i)]  Identity $ 	\mathfrak{L}(\rho) $ (2.17) holds in $ \rho$
		\item[(ii)]  If  
		 $ G \approx H, \; \dim \rho = \dim \sigma $ and  $ 	\mathfrak{L}(\rho)$ holds in $ \sigma$	then there is a conjugate class size preserving permutation $ \beta$ of the list of conjugate classes (2.13) such that 
		$$ \chi_{\rho,j} = \chi_{\sigma, \beta(j)}, \; j= 1,\cdots, s $$ 
	\end{enumerate}
	
\end{theorem}
\noindent We will sketch a proof of the statement (ii). By construction, $ \mathcal{D}_{s}$ is non-identity of representation $ \sigma$ only if  
$ s' \geq s $ and $ t'_i \geq t_i, \; i = 1, \cdots, s $. For example,  $	\mathcal{S}_1(x_1,Y_1)$ is non-identity of $\sigma$ only if  $ t'_1 = |C'_1|  \geq t_1 $ with the only choice for "non-vanishing" variable assignment being  $	 x_1 \rightarrow h_1 \in C'_1 $   and   
$ \{ y_1, \cdots, y_{t_r} \} \rightarrow $ "a set of pairwise distinct coset representatives of the centralizer of $h_1$ in $H$". 
Further along the formula (2.15), the expression
$$
	\mathcal{D}_{2} = 	\mathcal{S}_1(x_1,Y_1)  	\mathcal{C}_{12}(x_1, x_2; Y_1) 
$$
where 
$$
	\mathcal{C}_{12}(x_1, x_2; Y_1)  =   \prod_{y \in Y_1  }( x_2 - y ^{-1} x_1 y ) v_i
$$
does not vanish in $ \sigma$ only if there is an assignment  $ x_2 \rightarrow h_2 $ such that 
$ h_2 $ is not conjugate to $h_1$. Continuing in this fashion  we note
that variable assignment for which $	\mathcal{D}_s \neq 0$ essentially fixes 
values $ h_r \leftarrow x_r, \; r = 1, \cdots , s $ to be
representatives of the conjugate classes $ 	C'_1,  C'_2, \cdots C'_{s'} $ such that $  s' \geq s  $ and $	|C'_i| \geq 	C_i, \; i = 1, \cdots, s $. In turn, 
the only possible choice of values assigned to $ Y_r$ is a set of pairwise distinct coset representatives of the centralizers of $ h_r, \; r = 1, \cdots , s $. By our assumptions, $H$ is isomorphic to $G$ and therefore  $ 	|C'_i| =  | C_i|, \; i = 1, \cdots, s = s' $.  With the variable assignment  described above ,  the $\Psi_a$-terms in $ \mathcal{E}_i , \; i = 1, \cdots, l $  will evaluate to $ (a/n)\chi_{\sigma}(h_r) I $ for all $ r = 1, \cdots , s$ and $ \mathcal{E}_i$-terms themselves will turn out to be non-negative real scalar matrices that must vanish in irrep $\sigma$ in order to satisfy identity $ \mathcal{L}_s $. From  Lemma 1 and the structure of terms (2.16)-(2.17), it follows, however,  that for this to happen, the values 
of $ \chi_{\sigma}$  should match some permutation of corresponding values of  $\chi_{\rho}$ on same-size conjugate classes, and that is  exactly what is claimed by (ii)
\subsubsection{Character Tables and Symmetric Groups} Denote by $\textnormal{irr}(\Gamma)$ the set of pairwise nonequivalent irreducible characters of a finite group $\Gamma$ and let 
$ \textnormal{cl}(\Gamma) $ be the set of conjugate classes of $ \Gamma $. Following e.g. \cite{Nenciu}, \cite{Wildon_L} we 
will say that two irreducible characters   $ \chi_1, \;\chi_2 \in  \textnormal{irr}(\Gamma)$ (or corresponding irreps themselves)  are \textit{table-equivalent}  
if there is a bijection $\beta$ of the set $\textnormal{cl}(\Gamma) $ such that $ \chi_1(C) = \chi_2 (\beta(C)) $ for any  class $ C \in  \textnormal{cl}(\Gamma) $. We will call  $ \chi_1, \;\chi_2 \in  \textnormal{irr}(\Gamma)$ \textit{strongly table-equivalent} if 
the bijection $ \beta $ preserves the conjugate class sizes, i.e $ |C| = | \beta(C)|, \; C \in \textnormal{cl}(\Gamma) $ (cf. Theorem 2.2 (ii)).
\begin{conjecture}
	Table equivalent irreps  are strongly table equivalent  
\end{conjecture}
\begin{example}
	Similar irreps are strongly table equivalent
\end{example}
	
\noindent In this terminology the statement of Theorem 2.2 (ii) can be read as  

\begin{corollary}
	Let  $ \rho_1, \; \rho_2$ be  faithful irreps of a finite group $\Gamma$ with characters $ \chi_1, \; \chi_2$.	If  identities $ \mathfrak{L}(\rho_1) $ (2.17)
	and $ 	\mathcal{D}_{\dim \rho_1} $ (2.6) hold in $ \rho_2$   
	then  $ \chi_1 $ and $ \chi_2 $
	are strongly table-equivalent
\end{corollary}
\noindent It is obvious that similar characters are (strongly) table-equivalent. In some cases the reverse statement is also true. The following theorem is a reformulation of a result of Wildon    
\begin{theorema}(Wildon, \cite{Wildon_L})
	Let $ \Gamma $ be an altetrnative group $ A_n $ (for any n) or a symmetric group $ S_n, \; n \neq 4$ and let   $ \rho_1, \; \rho_2$ be  faithful irreps of  $\Gamma$. Then 
	$ \chi_{\rho_1} $ is table equivalent to 	$ \chi_{\rho_2} $ if and only if representations $ \rho_1$
	and $ \rho_2 $ are similar
\end{theorema}
\begin{remark} An immediate consequence of this theorem is  that  table equivalent irreps of symmetric or alternative groups are strongly table equivalent (cf. Example 2.7, and Table 1 for the special case of $S_4$).  
\end{remark} 
\noindent In combination with  Corollary 2.10, this result of Wildon can be used to construct explicit identities that separate non-similar irreps of symmetric or alternative groups, as we have
\begin{corollary} With the exception of $S_4$, faithful irreducible representations $ \rho_1, \; \rho_2 $ of a symmetric group $S_n$ or an alternating group $A_n$ are similar if and only if   irrep $ \rho_2$ satisfies  identities $\mathcal{L}(\rho_1)$  (2.17) and  $ 	\mathcal{D}_{\dim\rho_1} $ (2.6) 
\end{corollary}
\begin{corollary}
	Let $ |H| = m \; (= |G|) $ and suppose that $G$ is a symmetric group $S_n, \; n \neq 4$ or an alternative group $A_n$ (for any $n$).  If  identities  $ \mathfrak{D}(G)$ (see Lemma 1.1), $\mathcal{L}(\rho)$  (2.17) and  $ 	\mathcal{D}_{\dim\rho} $(2.6)  hold in $\sigma$ then irreps $ \rho $ and $ \sigma$ are similar
\end{corollary}

\begin{question}
	When  characters with the same signature ranges are table equivalent?
\end{question}

\noindent 	Let's look at the  exception in Corollary 2.11  more closely. The character table of $ S_4 $ is reproduced from \cite{S4} in Table 1 
\begin{example}
	Faithful non similar characters $\rho_4, \; \rho_5$ of the character table (cf. Table 1) of $S_4$   differ on conjugate classes of size $6$. One of these conjugate classes  corresponds to  the full cycle  $(1234)$ and another one corresponds to the partition $(2, 1, 1)$ The characters $\rho_4$ and $\rho_5$  are strongly table equivalent
	and none of the identities discussed so far can be used to distinguish between them. However, using methodology  described above, it is not hard to produce an identity that separates $ \rho_4$ from $ \rho_5$. Note, that if an element    $ x \in S_4$ is not a cycle of length $4$  then $  x^6 = 1$. Therefore, just by looking at the appropriate line of Table 1 one can write down an identity satisfied by the representation with character $ \rho_4$:
	\begin{align}
		\mathcal{C}_{24}(y_1, \cdots, y_{24}) (x^6 - 1 ) y_{25} ( \Psi_{24}(x,Y_{24}) - 24/3)  \nonumber
	\end{align}	
Corresponding identity for $ \rho_5$ is 
	\begin{align}
		\mathcal{C}_{24}(y_1, \cdots, y_{24})  (x^6 - 1 ) y_{25}( \Psi_{24}(x,Y_{24}) + 24/3)  \nonumber
	\end{align}	
	and either of these two identities distinguishes between $ \rho_4$ and $\rho_5$.
\end{example}
\noindent Following this example it is natural to define \textit{disjunctive trace identities} (cf. \cite{RT}-\cite{PT2} and \cite{PK}). Let $ u_1,  \cdots u_p; \;  v_1, \cdots v_q  \in F(Y)  $. Using formal trace symbol
$ Tr(f), \; f \in F(Y)  $  consider disjunctions over the free group $ F(Y) $ of the form 
\begin{align}
	u_1 = 0 \vee \cdots \vee u_p = 0 \vee P_1 = 0
	\vee \cdots \vee P_s = 0  
\end{align} 
where $ P_j, \; i = j, \cdots, s$ are polynomials in formal variables
$  Tr(v_i), \; i = 1, \cdots, v_q $. Any such polynomial, $P$ can be turned into a "valid" element $ \hat{P}$ of $kF(Y)$ by substitution 
$$ Tr(v_i) \leftarrow (n/m)\Psi_m(v_1, Y_m)   $$
and it is obvious that if the disjunctive formula (2.18) identically holds in the irrep $ \rho$ then $ \rho $ satisfies the identity
\begin{align}
		\mathcal{C}_{m}(Y_m) u_1 x_1 \cdots  u_p x_p \hat{P}_1  x_{p+1} \cdots x_{p+s-1} \hat{P}_s \tag{2.18*}
\end{align}
where $ Y_m = \{y_, \cdots, y_m \} \subset Y, \; x_i \in Y\setminus Y_m, \; i = 1, \cdots, p+ s - 1 $ and it is assumed that none of these variables appear in the reduced form of free group elements  $ u_1,  \cdots u_p; \;  v_1, \cdots v_q  $. 
 
 It is probably  a straightforward exercise to formalize the notion of  "disjunctive trace identities" along the lines suggested in \cite{PK}, \cite{RT}-\cite{PT2} and it should be clear from the preceding discussion that if faithful finite dimensional irreps $\rho$ and $\sigma$ have the same disjunctive trace identities then at the very least $ G \approx H $ and $  \textnormal{range}(\rho) = \textnormal{range}(\sigma) $. 
 
 \begin{example}
 	Let  $ D = \{ d_1, \cdots, d_o \} $ be a set of all possible element orders in (order statistics of) the group G. Let $ \{ d_{ij},  \; i = 1, \cdots, o; \; j = 1, \cdots, l_i \} $  be a set of   fixed point set dimensions of cyclic subgroups of $ G $ of order $ d_i $ in irrep   $ \rho $.  For any $ i =1, \cdots, o $ set $ D_i = D \setminus \{i\}$ and  set
 	$$ w_i =   \bigvee_{d \in D_i } 
 	(x^d = 1  ) \bigvee \bigvee_{ 1 \leq j \leq l_i } \left( \frac{1}{d_i}\sum_{j=0}^{d_i - 1}Tr(x^j) = d_{ij} \right)  , \; x \in Y $$
 	By definition  all disjunctions $ w_i $  hold in   $\rho$ and the recipe (2.18) can be used to write "fixed point identities" for an irrep. 
 	In particular (cf. 2.18*), if $ G $ acts without fixed points then $ \rho $ satisfies identities
\begin{align}
		\mathcal{C}_{m}(Y_m) \prod_{d \in D_i } 
		(x^d - 1  )y_d \sum_{j=0}^{d_i - 1}\Psi_m(x^j, Y_m), \; y_d \in Y, \; d \in D, \; i = 1, \cdots, o \nonumber
\end{align} 	 
 \end{example}
 \begin{conjecture} Faithful irreps of finite groups that have the same disjunctive trace identities
 	are similar.
\end{conjecture}       
\begin{table} 
	$$
	\begin{array}{c|rrrrr}
		\rm class&\rm1&\rm2A&\rm2B&\rm3&\rm4\cr
		\rm size&1&3&6&8&6\cr
		\hline
		\rho_{1}&1&1&1&1&1\cr
		\rho_{2}&1&1&-1&1&-1\cr
		\rho_{3}&2&2&0&-1&0\cr
		\rho_{4}&3&-1&-1&0&1\cr
		\rho_{5}&3&-1&1&0&-1\cr
	\end{array}
	$$
	\caption{Character table of $S_4$ (\cite{S4})}
\end{table}

\section{Trace  Identities}
Aby $ n \times n $ matrix $A$ satisfies  \textit{Cayley-Hamilton
identity} (cf. e.g. \cite{TF}) 
\begin{align}
	p_A(x) = 	\sum_{0}^n (-1)^{i} \text{tr}(\wedge^i A )x^{n-i} \nonumber
\end{align}
where $ \text{tr}(\wedge^i A )$ is a trace of the $i$-th exterior power of the linear operator $A$. In turn, coefficients 
of Cayley-Hamilton polynomial $ p_A(x) $ are polynomials in
 $\text{tr}(A), \; \text{tr}(A^2), \; \cdots $ (see \cite{TF}-\cite{TF1} and refrences therein).
Actually, (cf. e.g. \cite{PT2}) any matrix $ A\in M_n(k) $ satisfies a one-variable \textit{ trace identity  }
\begin{align}
T_n(x) = x^n + \sum_{i=1}^n (-1)^i \sigma_i(x) x^{n-i}
\end{align}
where 
\begin{align}
	\sigma_i(x) = \frac{1}{i!} \det
	\begin{bmatrix}
		\!\text{tr}(x) & i-1 & 0 & \cdots & 0 \\
		\text{tr}(x^2) & 	\text{tr}(x) & i-2 & \cdots & 0 \\ 
		\cdot & 	\cdot  & \cdot  & \cdots &  \cdot \\
			\cdot & 	\cdot  & \cdot  & \cdots &  \cdot \\
		\;\;\;\text{tr}(x^{i-1}) & \;	\text{tr}(x^{i-2}) & \cdot & \cdots & 1 \\
		\text{tr}(x^{i}) & 	\; \text{tr}(x^{i-1}) & \cdot
		 &  \cdots & \text{tr}(x)
	\end{bmatrix}   \nonumber
\end{align} 
are polynomials in $ \text{tr}(x), \; \text{tr}(x^2), \cdots $. 
 For example, if $n=2$ then  
\begin{align}
	T_2(x) = x^2 - \text{tr}(x)x + \textnormal{det}_2(x)  \textnormal{ where }\textnormal{det}_2(x) = 1/2 ( \textnormal{tr}(x)^2 -  \textnormal{tr}(x^2) )  
\end{align} 
Formally replacing traces with $ \Psi$-terms   in  $ T_n$ and using Lemma 2.1 we get a generic identity that holds in any $n$-dimensional exact irrep of any finite group of  order $m$
\begin{example}
Faithful $n$-dimensional irrep of a group $\Gamma$ of order $ \leq m $ satisfies the Hamilton-Cayley identity
\begin{align}
	\mathcal{C}_m (y_1, \cdots, y_m) T_n(x) |_{	\text{tr}(x^k) \leftarrow (n/m)\Psi_m( x^k, \; y_1, \; \cdots, \; y_m )} 
\end{align}  
\end{example}
\begin{example}
	Using the same substitution as in (3.3) one can turn  the   matrix polynomial $ \sigma_i(x) $ 
	into an element in $ kF(Y) $
	\begin{align}
		\hat{\sigma}_i(x,Y_m) = \sigma_i(x) |_{	\text{tr}(x^k) \leftarrow (n/m)\Psi_m( x^k, \; y_1, \; \cdots, \; y_m )} \nonumber 
	\end{align}
	thus obtaining a set of identities of $ \rho$ 
	\begin{align}
		\sigma_i(\rho)  =	\mathcal{C}_m (y_1, \cdots, y_m)\prod_{d \in \Delta_i} (		\hat{\sigma}_i(x,Y_m) - d) v_d, \;\; i = 1, \; \cdots, \; n  
	\end{align}
	where $ \Delta_i =  \{  \sigma_i(\rho(g )), \; g \in G \} $ is the set of all possible values of the trace polynomial  $ \sigma_i $ on $ \rho(G) $ and  $\{ v_d, \; d \in \Delta_i \} \subset Y \setminus Y_m \setminus \{x\} $ is a set of free variables. In particular, $\sigma_1(\rho) $ is the character value identity (2.4) and $  D_n(\rho) \equiv \sigma_n(\rho) $ is the determinant value identity that among other things, can be used to characterize irreducible representations of groups of order $m$ into special unitary group $SU(n)$. More precisely,  the identity  
	$$   \mathcal{C}_m (Y_{m}) (		\hat{\sigma}_n(x,Y_m) - 1 ) $$   
	holds in $\rho$ if and only if $ \rho(G) \subset 
	\textnormal{SU}(n)$
\end{example}

\subsection{Adams Operations and Gassmann Equivalence}
In a way quite similar to the Example 3.2 one can use Adams operations (cf. e.g. \cite{Serre2}, \cite{Adams}, \cite{Meir}) to "customize"  generic Cayley-Hamilton identity (3.3). We define  the "Adams map" $$ A \equiv A_{\rho}: \rho(G) \rightarrow \mathbb{C}^{n=\dim\rho} $$ associated  to the faithful irrep  $ \rho $   
as 
\begin{align}
 A(\rho(g)) = ( \chi_{\rho}(g), \chi_{\rho}(g^2), \cdots, \chi_{\rho}(g^n)  ) \nonumber
\end{align}

\noindent Let $ r $ be the cardinality of the image of $ A $. Level sets of the map $A$  form a partition
\begin{align}
 \mathfrak{A} \equiv \mathfrak{A}(\rho) = G_1 \cup \cdots \cup  G_r   \tag{3.5.0}
\end{align}
 of $ \rho(G) \subset \textnormal{GL}(V) $ and  we have

 \begin{lemma} There are distinct conjugacy classes $ O_1, \cdots, O_r $ of $ \textnormal{GL}_n(\mathbb{C}) $ such that $ G_i = O_i \cap \rho(G), \; i = 1, \cdots , r $  
 \end{lemma}
\noindent Proof. Two unitary matrices are conjugate in   $ \textnormal{GL}_n(\mathbb{C}) $ if and only if
they have the same spectrum counting  multiplicity. On the other hand, for any $  g,h \in G$ one has  $ A(\rho(g)) = A(\rho(h))  $ if and only if characteristic polynomials 
of $ \rho(g) $ and $ \rho(h)$ are the same (cf. (3.1)).
\newline
\newline
We can reformulate Lemma 3.1 as follows. 
For $ g \in G $ let $ S( \rho(g))  $ denote the  spectrum  of $ \rho(g) $ (as a size $n$ multi-set of complex numbers). Thus we have another map $ S \equiv S_{\rho} : \rho(G) \rightarrow \mathbb{C}^n $. It follows (e.g. from Lemma 3.1) that  level sets of the maps $A$ and $S$ coincide and therefore both maps define the same partition $ \mathfrak{A}_{\rho} $ of $ \rho(G)$. 
Speaking informally we will call the set $ S(\rho(G)) \subset \mathbb{C}^n$ the spectrum of  irrep $ \rho$ and we will say   that faithful representations $ \rho_1 $ and $ \rho_2$ \textit{have equal spectrum} if $\dim \rho_1 = \dim \rho_2$  and $ S_{\rho_1}(\rho_1(G)) = S_{\rho_2}(\rho_2(G)) $.  
Define  \textit{spectral signature} of $ \rho$ as  a
set of pairs 
$ \{ (|G_i|, S(G_i) \}, \; i = 1, \cdots r  $ (cf. 2.3.1). 
 It is clear that each $ G_i \in  \mathfrak{A}_{\rho} , \; i=1,\cdots ,r $ is a disjoint union of $\rho$-images of some conjugate classes  of $ G $ and therefore $ r \leq s $ (see the introduction).
\newline\newline
\noindent The image of the map $ A_{\rho} $ is  a set of $r $ points in $ \mathbb{C}^n $ and, therefore, it can be identified (row per point) with an $ r \times n $ matrix as follows   
\begin{align}
 \textnormal{"image of } A" = [a_{i,k}], \; i = 1, \cdots, r; \; k = 1, \cdots , n 
\end{align}
By analogy with (2.4) and (2.9)-(2.11) set 
\begin{align}
\!\!\!\!\!\!\!\!\!		\mathcal{A}_i(\rho) \equiv		\mathcal{A}_i(\rho)(x,Y_m) = \sum_{k=1}^{n}(\Psi_m(x^k,Y_m) - \frac{m}{n} a_{i,k})(\Psi_m(x^{-k},Y_m) - \frac{m}{n}\bar{a}_{i,k})
\end{align}
and 
\begin{align}
\!\!\!\!\!	 \mathcal{A}(\rho) \equiv \mathcal{A}(\rho)(x,Y_m)
   = \;
	\mathcal{C}_m(Y_m)  
\mathcal{A}_{1}(\rho) v_1 \cdots v_{r-1} \mathcal{A}_{r}(\rho) \\
		\!\!\!\!\!\!\!	\mathcal{S}_i(\rho)	 \equiv	\mathcal{S}_i(\rho)(X_m,Y_m) = \mathcal{C}_m (X_m) \mathcal{C}_m (Y_m) 
		\prod_{ j=1}^m  \mathcal{A}_{i}(\rho)( x_j, Y_m )    w_j	 \\
		\!\!\!\!\!\!\!\!\!\!\!\			
	\mathcal{G}_i(\rho,X_m) \equiv		\mathcal{G}_i(\rho)  = \mathcal{C}_m (X_m) \mathcal{C}_m (Y_m) 
		\prod_{S \subset X_m,\; |S|= |G_i| }   \left( \sum_{x \in S}  \mathcal{A}_{i}(\rho)( x, Y_m ) \right) v_S
\end{align}
where 
 $ X_m =
\{x_1, \cdots x_m\} \subset  Y \setminus Y_m $, $ v_1, \cdots, v_{r-1}, w_1, \cdots w_m, v_S \in Y \setminus Y_m \setminus X_m $ are distinct free variables and $ i = 1, \cdots,  r $. 

To describe  properties of irreps that are characterized by identities (3.7-9) we need some terminology that invariably shows up in a context of \textit{manifold isospectrality} (cf. e.g.    \cite{Sutherland}, \cite{Mackey}-\cite{Sunada}, \cite{Pointwise} and Section 3.1.2 below) 
\begin{definition}
	Faithful representations  
	$ \rho_i : \Gamma_i \hookrightarrow \textnormal{GL}_n(k), \; i = 1,2  $ of finite groups $\Gamma_1, \; \Gamma_2 $ are called \textbf{Gassmann
	equivalent} (or \textbf{almost conjugate})  if the following equivalent conditions hold:
	\begin{enumerate}
		\item[(a)] $\rho_1(\Gamma_1) $ and $ \rho_2(\Gamma_2) $ have the same spectral signature
		\item[(b)] 
		$ | \rho_1(\Gamma_1) \cap O | =   | \rho_2(\Gamma_2) \cap O |  $
		for any orbit $ O $ of  adjoint $ 	\textnormal{GL}_n(k)$-action
		
		\item[(c)]  there is a spectrum preserving bijection $ \gamma: \Gamma_1  \leftrightarrow \Gamma_2$, i.e. $\rho_1(h) $ and $  \rho_2(\gamma h) $ are conjugate 
		in  $ \textnormal{GL}_n(k) $  for all $ h \in \Gamma_1$ 
	\end{enumerate}
\end{definition}
\begin{remark}
	This is not the most general definition For more details see the references cited above.
	
	It should be clear that Definition 3.1 could be applied verbatim to a pair of subsets in $\textnormal{GL}_n(k)$. In fact, one can say that representations $ \rho_1 : \Gamma_1 \hookrightarrow \textnormal{GL}(V) $ and $ \rho_2 : \Gamma_2 \hookrightarrow \textnormal{GL}(V) $  are Gassmann equivalent if (and only if) subsets $ \rho_1(\Gamma_1), \; \rho_2(\Gamma_2) \subset \textnormal{GL}(V) $ are Gassmann equivalent
\end{remark}
\noindent The following lemma lists some basic properties of Gassmann equivalency and provides an illustration for the definition 
\begin{lemma} 
	\; 
	\begin{enumerate}
		\item[(i)]  Faithful similar representations are Gassman equivalent
		\item[(ii)]  Faithful Galois conjugate representations are Gassman equivalent
		\item[(iii)] If  the bijection $\gamma$ in a statement (c)  of Definition 3.1 is an isomorphism then representations $ \rho_1 $ and $ \rho_2 $ are similar
		\item[(iv)] (cf. e.g \cite{Pointwise}). If representations of 
		finite abelian groups $A $ and $ B $ are Gassman equivalent `then $ A \approx B $ 
			\item[(v)]  Faithful representations of cyclic groups are  almost conjugate if  and only if they are similar (cf. Appendix 2) 
	\end{enumerate}.
\end{lemma}
\noindent 	Proof. If faithful representations  
$ \rho_i : \Gamma_i \hookrightarrow \textnormal{GL}_n(k), \; i = 1,2  $ of finite groups $\Gamma_1, \; \Gamma_2 $ are similar then there is 
an isomorphism $ \alpha : \Gamma_1 \rightarrow \Gamma_2 $ and an invertible matrix $ X \in  \textnormal{GL}_n(k)$ such that $ \rho_2(\alpha g) = X^{-1} \rho_1(g) X$ for all $ g \in \Gamma_1$, hence the isomorphism $ \alpha$  is a bijection required by definition 3.1 (c). 

 Turning to the statement (ii), replace $\rho_2$ by an equivalent representation if necessary, to get $ \rho_2(g) = \epsilon \rho_1(g), \; g \in \Gamma $ for some $ \epsilon \in \mathfrak{G}(\Gamma)  $.
 By Lemma 2.6 there is an integer $t$ co-prime
 to $ |\Gamma| $  such that  $ \rho_1(g^t) $ and $  \epsilon \rho_1(g) $
 have the same spectrum for all $  g \in \Gamma $.
 Hence the condition (c) of Definition 3.1 can be satisfied by a bijection $ g \rightarrow g^t$

 If  there is an isomorphism $ \alpha : \Gamma_1 \rightarrow \Gamma_2 $ such that 
$ \rho_2(h)$ is conjugate to $ \rho_1 \alpha (h) $ for all $ h \in \Gamma_1$  then  $ \chi_{ \rho_1 \alpha } = \chi_{\rho_2} $ and (iii) follows from definition 1.1.

To verify the statement (iv), split $ A$ and $ B $ into direct sum of cyclic subgroups and count the number of same order elements (cf. e.g. \cite{Pointwise})

Finally, if $ \rho_1, \rho_2$ are almost conjugate (hence faithful) representations of cyclic groups  $ \Gamma_1, \;  \Gamma_2 $ then by (iv) $ \Gamma_1 $ is isomorphic to $ \Gamma_2$ and we can assume   that $ \rho_1, \; \rho_2 $ are almost conjugate  representations of the same cyclic group $  \Gamma_1 \approx \Gamma_2 \approx \Gamma \; = \;\; <h>$. By Gassman equivalence  there is an integer
$ t $ such that $ \rho_1(h) $ has the same spectrum as 
$ \rho_2(h^t) $. Clearly such $t$  must be coprime to $|\Gamma|$ and it follows then, that automorphism  $ h \rightarrow h^t$ entails a   similarity  between $ \rho_1 $ and $ \rho_2 $

\begin{proposition} Identities (3.7-9) hold  in $\rho $ and   if $|H| < |G| $ all of them hold in  $ \sigma $.  If
	$  \dim \sigma = \dim \rho $ and $ |H| = |G|$ then
	\begin{enumerate}

		\item[(i)] Identity $\mathcal{A}(\rho) $ (3.7) hold in $ \sigma$ if and only if the spectrum of  $ \sigma $ is contained in  the spectrum of $ \rho$, i.e., 
	 for any $ h \in H $ there is $ g \in G $ such that $ \sigma(h) $ and $ \sigma(g) $ have the same spectrum  
	
	\item[(ii)] Identity $\mathcal{A}(\rho)$  (3.7) and all identities $\mathcal{S}_{i}, \; i = 1, \cdots, r $ (3.8) hold in $ \sigma$ if and only if   $ \rho $ and $\sigma$ have the same spectrum 

		\item[(iii)]   All identities $\mathcal{G}_{i}, \; i = 1, \cdots, r $ (3.9) hold in $ \sigma$ if and only if
	 $\rho$ and $ \sigma $ are Gassman equivalent 
	\end{enumerate}
 
\end{proposition}
\begin{corollary}
Suppose that $ |H| = m = |G| $.
	If identities $\mathcal{D}_n $ (2.6), $ \mathfrak{D}(G) $
	(Lemma 1.1 (2)) and  $\mathcal{G}_{i}, \; i = 1, \cdots, r $ (3.9) hold in $\sigma$ 
then $ H \approx G $ and  $\sigma$ is Gassmann equivalent to  $ \rho$
\end{corollary}
\begin{remark}
Corollary 3.1 is a significant refinement of the results of the previous section. Essentially we can say now that
	acting group, dimension, character values and 
	spectrum  of an exact irrep of a finite group of a given order are determined by a (some selection from) well defined  sets of identities (1.2),  (2.4), (2.6), (2.9), (3.7-9). Moreover,  identities (3.9) determine an irrep up to Gassmann equivalency
\end{remark}

\noindent  We omit the a proof of the Proposition 3.1, since the following explanation should be sufficient.
 Identities (3.7-3.9) and corresponding identities (2.4), (2.10), (2.11) are structurally similar. The only difference is that terms $ \Psi_{m,i} $ (2.3) in (2.4), (2.10), (2.11) are replaced with terms 
$ \mathcal{A}_i $ (3.6) in (3.7-3.9). The meaning of this substitution is not hard to understand. The term (2.3) dictates the character value while the term (3.6) (in similar context) determines the spectrum (values of coefficients of characteristic polynomial). Note also, that entries
$ a_{i,1}$ of the matrix (3.5) are still character values, and therefore, the spectral signature is just a refinement of the character range signature. 

It is now quite obvious, that the identity $\mathfrak{L}(\rho)$ (2.17) can be also modified to control Gassmann equivalency. Indeed, the 
terms $ (\Psi_{t_a}(x_a,Y_a) - (t_b/n) \chi_{b,\rho} )  $ in (2.16) can be replaced with appropriate $ 	\mathcal{A}_i(\rho) $ terms of (3.7) leading to  the identity of the irrep $ \rho$ that we will denote by    
$ \mathfrak{L}'(\rho) $. 

\begin{definition}
 We will say that Gassmann equivalent representations  are
\textit{strongly Gassman equivalent} if the map $ \gamma$ in Definition 3.1 (c) takes a conjugate class into a conjugate class
\end{definition}
\noindent It is easy to see that faithful strongly Gassman equivalent irreps are strongly table equivalent. Thus, as a consequence (from the proof of) Theorem 2.2 we have
\begin{corollary}
		\noindent 				 
	\begin{enumerate}
		\item[(i)]  Identity $ 	\mathfrak{L}'(\rho) $  holds in $ \rho$
		\item[(ii)]  If  
		$ G \approx H, \; \dim \rho = \dim \sigma $ and  $ 	\mathfrak{L}'(\rho)$ holds in $ \sigma$	then $\sigma$ is strongly Gassmann equivalent to $\rho$

	\end{enumerate}
\end{corollary}
 \begin{remark}
	A representation of an abelian group can be Gassman equivalent to a representation of non-abelian group  (see e.g. \cite{Sutherland}). Hence, (at  least non-irreducible) representations can be Gassmann equivalent but not strongly Gassman equivalent.  It is probably an open question, however, whether Gassmann equivalent irreps of the same group are  strongly Gassmann equivalent. Note also  that, as before, 
Galois conjugate irreps can not be distinguished by identities (3.7-9) as these identities are invariant under action of Gallois group (cf. Lemma 2.5)

\end{remark}
\begin{example} 
	Adams operations can be used to  extend the character table of a finite group $\Gamma$ into a $3$-tensor  $$a_{ijk} = \Psi^k(\chi)(c), \; \chi \in \textnormal{irr}(\Gamma), \; c \in 
	\textnormal{cl}(\Gamma),  \; 1 \leq k < \mathfrak{e} $$ where  $\mathfrak{e}$ is an exponent of $\Gamma$.     
	There are examples  of finite $p$-groups that are not determined by this tensor not up to isomorphism and  not even up to similarity (cf. \cite{Dade} and a discussion in \cite{Meir}).
	It is easy to see that the  example in  \cite{Dade}  demonstrates  existence of Gassmann equivalent  faithful irreps of non-isomorphic finite  $p$-groups. Therefore, two faithful irreps of non-isomorphic groups can satisfy same identities (3.7-9)  
\end{example}
\noindent  Next we discuss some  examples related to Proposition 3.2. In particular, we will see that Gassmann equivalent irreps of  \textit{a given finite group} are not necessarily similar or Galois conjugate.     
\subsubsection{Spectrum and group element order statistics (cf. \cite{S}-\cite{S2})}

\noindent In contrast with Example 3.3, the identities (3.7-9) do imply group isomorphism in case of irreps of simple groups.  This fact is a direct consequence of a remarkable group theory result  (see \cite{S1} and a discussion in \cite{S}, \cite{S2}) that we will now briefly describe. 
Following \cite{S} denote by $ \pi_e(G) $ the integer set of all element orders  in a group $G$.    

\begin{theorema} (\cite{S1}) If $ G $ is a simple finite group then $ H \approx G $ if and only if  $ |H| = |G|$
and $ \pi_e(G) = \pi_e(H) $ 	
\end{theorema}
\noindent If  irreps $ \rho $ and $\sigma$ have the same spectrum   then, of course,  $ \pi_e(G) = \pi_e(H) $ and  this theorem can be straightforwardly combined with  Proposition 3.2 (ii), yielding
\begin{corollary} (cf.\cite{S1}).	
Suppose that $ G $ is simple and that $ |G| = |H| $. If identities $ \mathcal{A}(\rho) $ (3.7),  $\mathcal{S}_{i}(\rho), \; i = 1, \cdots, r $ (3.8) and $	\mathcal{D}_{\dim\rho} $ (2.6) hold in $ \sigma$ then 
	$	 G \approx H  $,  $  S(\rho(G)) = S(\sigma(H)) $ and $ \textnormal{range}(\rho) = \textnormal{range}(\sigma)  $
\end{corollary} 
\begin{corollary} (cf. Corollary 2.12).
	Let $ |H| = m \; (= |G|) $ and suppose that $G$ is a simple  alternating group.  If  identities   $\mathcal{L'}(\rho)$  (2.17) and  $ 	\mathcal{D}_n $ (2.6)  hold in $\sigma$ then irreps $ \rho $ and $ \sigma$ are similar
\end{corollary}

\subsubsection{Gassmann  equivalence and similarity. Orbifolds and Spherical Space Forms (\cite{Swa}, \cite{Gor}, \cite{Ikeda}-\cite{Sutton})}
Bearing in mind Corollary 3.1, one may ask if  two Gassmann quivalent  irreps of a given finite group are similar.
It turns out that the answer to this question is  "no", and this fact was known for some time  in a context of  \textit{isospectral spherical space forms} (cf. \cite{Ikeda}, \cite{Wolf}, \cite{Wolf2}). We will review a relevant example in some detail bellow. Before doing that, however, we would like to point out that it was established in \cite{Swa} that the notion of representation similarity plays an important role in geometry of more general orbit spaces (and orbifolds). In fact, it is shown in \cite{Swa} (see also \cite{Gor}) that same-dimension exact orthogonal representations of  finite groups are similar if and only if their orbit spaces are isometric. We will briefly discuss some details, 

Let $ V $ be a real Euclidean vector space on which the orthogonal group $ O(V)  $ naturally acts. Let $ S(V)  $ denote a unit sphere in $ V $. A finite subgroup $ \Gamma \subset O(V) $ acts on $S(V)$ and the orbit space of this action
$ X = S(V)/\Gamma $ is one of the basic    examples
of \textit{good} orbifolds (cf. e.g. \cite{Davis}, \cite{Gor}).
\begin{remark}
	General theory of (Riemannian)
	orbifolds (cf. e.g. \cite{Davis}, \cite{Gor}) is well beyond the scope of this  paper.
	We will just point out that the orbifold  $ X$ inherits its metric from $S(V)$. The following   quote from \cite{Swa} could be helpful: "In general $X$ will have singularities but, except for $S(\mathbb{R}^1)/\{ I \}$, it will still be
	a length space. That is, the distance between two points in X is the infimum
	of the lengths of all paths between the two points"
\end{remark}
\noindent  
\begin{lemma}(\cite{Swa}).  Let $ \rho_i : \Gamma_i \rightarrow \textbf{O}(V) $ be real orthogonal exact representations of finite groups $ \Gamma_i, \; i = 1,2 $. The orbifolds 
	$S(V)/\rho_i(\Gamma_i) $ are  isometric if  and only if representations $ \rho_i, \; i = 1,2 $ are similar
\end{lemma}
\noindent Combining this result with Theorem 1.1 we unexpectedly get 
 \begin{corollary} 
 	 Exact orthogonal irreps $ \rho_i : \Gamma_i \rightarrow \textbf{O}(V)$ of finite groups   have  same identities
 	 if and only if  
 	  orbifolds 
 	$S(V)/\rho_i(\Gamma_i), \; i = 1,2  $ are  isometric
 \end{corollary} 
\begin{remark}
	It is kind of obvious that there is a connection between identities of representations and polynomial invariants. In turn, polynomial invariants are very much related to orbit spaces. What seems to be surprising is that representation identities could be related to (Rimannian) metrics on  orbit spaces   
\end{remark} 
\noindent  Unitary representations
$ \rho_i : \Gamma_i \rightarrow \textbf{U}(V), \; i = 1,2; \; \dim V = n $ can be viewed as $2n$-dimensional real representations $ \prescript{}{\mathbb{R}}\rho_i : \Gamma \rightarrow \textbf{O}(\prescript{}{\mathbb{R}}V) $ where the left index $\prescript{}{\mathbb{R}}*$ means the restriction of the field of scalars. The groups $\Gamma_i$ act on a $(2n-1)$-dimensional unit sphere 
$ S(V) = S(\prescript{}{\mathbb{R}}V)  $ since   $ \textbf{U}(V) \subset \textbf{O}(\prescript{}{\mathbb{R}}V) \approx  \textbf{O}(2n-1) $. If
representations $ \rho_i$  are similar then  representations $  \prescript{}{\mathbb{R}}\rho_i  $
 are also  similar and 
 orbifolds $ S(V)/\rho_i(\Gamma_i) $ are isometric  by (an easy part of) Lemma 3.3. It is easy to see, however, that orbit spaces of complex conjugate representations  are isometric. In any case, Lemma 3.3 directly yields  the following 
 \begin{corollary}
 	 Let $ \rho_i : \Gamma_i \rightarrow \textbf{U}(V) $ be exact unitary representations of finite groups $ \Gamma_i, \; i = 1,2 $. 
 	\begin{enumerate} 
 		\item[(i)]
 		If $ \rho_1 $ and $ \rho_2 $ are similar then orbifolds  	$S(V)/\rho_i(\Gamma_i), \; i = 1,2 $ are  isometric
 		\item[(ii)] If orbifolds  	$S(V)/\rho_i(\Gamma_i), \; i = 1,2 $ are  isometric then $ \Gamma_1 \approx \Gamma_2$
 		\item[(iii)] Assuming  that  representations $\rho_i, \; i = 1,2 $ are irreducible, suppose that  orbifolds  	$S(V)/\rho_i(\Gamma_i) $ are  isometric. 
 		Then  $ \rho_2 $ is similar either to $\rho_1$ or to $ \bar{\rho}_1 $ 
 	\end{enumerate} 
 \end{corollary} 
\noindent  Proof. If representations 	$S(V)/\rho_i(\Gamma_i), \; i = 1,2 $ are  isometric then $  \prescript{}{\mathbb{R}}\rho_i  $ are similar as real representations by Lemma 3.3 and therefore  $ \Gamma_1 \approx \Gamma_2$.  Hence, to prove (iii) we can assume that $ \Gamma_1 = \Gamma_2 = \Gamma $. Again by Lemma 3.3,  irreps $  \prescript{}{\mathbb{R}}\rho_i, \; i = 1,2  $ are similar as real (orthogonal) representations, so there is an automorphism $\alpha$ of $ \Gamma $ such that 
\begin{equation}
\textnormal{ complexification of }   \prescript{}{\mathbb{R}}\chi_{\rho_2}  =  \chi_{\rho_2} +   \bar{\chi}_{\rho_2}   = \chi_{\rho_1}\alpha +   \bar{\chi}_{\rho_1}\alpha 
\end{equation}      
\noindent\noindent All four characters in equation (3.10) are irreducible, and therefore either $  \chi_{\rho_2} = \chi_{\rho_1}\alpha $   or $  \chi_{\rho_2} = \bar{\chi}_{\rho_1}\alpha  $ as  required by (iii).
\newline\newline

\noindent Non-similar Gassmann equivalent  representations are one of the main sources of examples of non-isometric but isospectral Riemannian manifolds  (cf. e.g. \cite{Sunada}). Following \cite{Ikeda} - \cite{Wolf} we will  take a look at a specific example that was introduced by   Ikeda (see \cite{Ikeda}, \cite {Ikeda_1}) in a context of \textit{spherical space forms} (cf. \cite{Wolf}, \cite{Wolf2}). 

Let $ M $ be a compact Riemannian manifold and let $ E_{\lambda}(M), \; \lambda \in \mathbb{R}$ denote 
the $\lambda$-eigenspace of the Laplacian (Laplace-Beltrami operator) $\Delta$ on $ C^{\infty}(M) $. It is well known that  the spectrum of Laplacian is real, discrete and all the multiplicities $ \dim E_{\lambda}(M)$ are finite.    
Riemannian manifolds $M$ that have the same spectrum  of  Laplacian (counting multiplicities) are called \textit{isospectral}. 

A finite group of isometries  $ \Gamma \subset   \textnormal{Aut}(M) $  is called \textit{ fixed point free (fpf)}  if the length of any orbit $ \Gamma m , \; m \in M $ is equal to $ |\Gamma| $. If that is the case,  the orbit space $M/\Gamma$ is also a compact Riemannian manifold in induced Riemannian metric. 
In particular, let $ M = S(V) $ be  $(2l-1)$-dimensional  (unit) sphere  in a complex $l$-dimensional vector space $ V $ so that $  \textnormal{Aut}(S(V)) \approx \textbf{O}(2l) $. If a finite group $ \Gamma$ has faithful fixed point free (fpf) orthogonal representation  $ \pi : \Gamma \rightarrow  V $  then the orbit space $ S(V)/\pi(\Gamma $) is called a \textit{spherical space form} (with fundamental group  $ \Gamma $)(\cite{Wolf}).
The book \cite{Wolf} (see also \cite{Allcock}, \cite{Gilkey}) contains a full classification of  spherical space forms ("Vincent Programme", cf \cite{Vincent}, \cite{Wolf}). We will restrict our attention here to  the so called fpf groups of type $I$. These are metacyclic groups  that we will define following  \cite{Wolf2}.

Let  the group $\Gamma_d \equiv \Gamma_d(m,n,r)$ be defined by generators $ A, B $ that satisfy relations  
\begin{gather}  
 A^m = 1 = B^n, \; BAB^{-1} = A^r  \nonumber \\ m,n \geq 1, \; (m, (r-1)n) = 1, \; r^d \equiv 1 \!\! \!\!\! \mod m  
 \end{gather}
with additional assumptions on integer parameters $ m, \; n, \; r $: 
\begin{gather}
	 d \textnormal{ is the order of the residue class of r in }  \mathbb{Z}_m^{\ast} \nonumber \tag{3.11.1} \\  
	 n = n'd \nonumber \tag{3.11.2} \\
	 n' \textnormal{ is divisible by any prime divisor of } d
	 \nonumber \tag{3.11.3}  
\end{gather}

\noindent It is easy to see  that $ \Gamma_d $ is a semidirect product of its cyclic gsubroups $<A>$ and $ <B> $ and that the cyclic subgroup $ <B^d>  $ is the center of $ \Gamma_d$.   

Let's summarize classification of fpf irreps of the group $ \Gamma_d$ as it appears in \cite{Wolf},  \cite{Wolf2}.  Take integers $k,l, \; (m,k) = 1 =  (l, n) = 1 $ and define one-dimensional characters $$ \sigma'_k: \; <A> \; \rightarrow \mathbb{C}^{\ast}, \; \sigma''_l : \; <B^d> \;  \rightarrow \mathbb{C}^{\ast}$$ by setting $ \sigma'_k(A) = \exp(2\pi i k/m), \;  \sigma''_l  (B^d) = \exp(2\pi i l/n')$. Let 
\begin{align}
	\pi_{k,l} = \textnormal{Ind}^{\Gamma}_{\Sigma} \; \sigma'_k \otimes \sigma''_l, \; \textnormal{where } \Sigma = A \; \times <B^d> \; \approx \; <A,B^d> \;  \subset \Gamma
\end{align}
One checks  that 
\begin{gather}   
\pi_{k,l}(A) = \textnormal{diag}(\; \exp(2\pi i k/m), \; \exp(2\pi i kr/m), \cdots, \; \exp(2\pi i kr^{d-1}/m) \;  )
 \nonumber \\ 
\pi_{k,l}(B) =
\begin{bmatrix}
	0 & 1 & 0 & \cdots & 0 \\
	0 &  0 & 1 & \cdots & 0 \\ 
	\cdot & 	\cdot  & \cdot  & \cdots &  \cdot \\
	\cdot & 	\cdot  & \cdot  & \cdots &  \cdot \\
	 0 & 0 & \cdot & \cdots & 1 \\
	\exp(2\pi i l/n') & 	0 & \cdot
	&  \cdots & 0
\end{bmatrix}     
\end{gather}
\noindent In other words,   $ \pi_{k,l}(B)$ is a $ d \times d $ (Toeplitz) matrix with 
the bottom left corner entry $ \exp(2\pi i l/n')$, ones immediately above the main diagonal and zeroes everywhere else. 
Complete classification of fixed point free irreps of the group 
$ \Gamma_d $ is given by the following theorem of Wolf.  
\begin{thA}(Wolf, \cite{Wolf2}, \cite{Wolf}) 
\begin{enumerate}
	\item[(i)]	
All representations $ \pi_{k,l} $ are irreducible and fixed point free. Any fpf irrep of $ \Gamma_d$ is equivalent to $\pi_{k,l} $ for  some integer $k, \; l \; ((m,k) = 1 = (l,n) = 1 )$
\item[(ii)]  $\pi_{k,l} \approx \pi_{k',l'} $ if and only if
$ l \equiv l' \!\!\! \mod \! n' $ and $ k = k' r^c \!\!\! \mod \! m $ with some integer $c$ such that $ 0 \leq c < d $.
\item[(iii)]  Irreps $\pi_{k,l} $ and $\pi_{ak,bl} $ are Gassmann equivalent for all integers $a,b$ such that $ (a,m) = 1  = (b,n) $
\item[(iv)]  Irreps $\pi_{k,l} $ and $\pi_{ak,bl} $ in (iii) are similar if and only if $b \equiv 1  \!\!\! \mod  \! d$
\end{enumerate}
\end{thA}
\begin{remark}
	The statement (iii) of the Theorem A  can be understood in the  following way (cf. \cite{Gilkey}).     
	Note  that  
	$$ \textnormal{Gal}(\mathbb{Q}({\omega_{mn}}) / \mathbb{Q}) =  \textnormal{Gal}(\mathbb{Q}({\omega_{m}}) / \mathbb{Q}) \times  \textnormal{Gal}(\mathbb{Q}({\omega_{n}}) / \mathbb{Q})$$ 
	and therefore
	characters 
	$ \sigma'_{ak} \otimes \sigma''_{bl} $ and 
	$ \sigma'_{k} \otimes \sigma''_{l} $ are
	Galois conjugate. It follows   that  characters 
	$	\pi_{ak,bl} = \textnormal{Ind}^{\Gamma}_{\Sigma} \; \sigma'_{ak} \otimes \sigma''_{bl}  $ and
	$	\pi_{k,l}  = \textnormal{Ind}^{\Gamma}_{\Sigma} \; \sigma'_{k} \otimes \sigma''_{l}  $ are also  Galois conjugate  (e.g. by Lemma 2.7) and hence Gassmann equivalent by Lemma 3.2 (ii).
\end{remark}
\noindent  
 According to Proposition 3.2,  identities (3.7-9)) do not distinguish  between Gassmann equivalent irreps  $\pi_{k,l} $ 
and $\pi_{ak,bl} $ of Theorem A (iii). However,  unless $b \equiv 1  \!\!\! \mod  \! d$, by Theorem 1.1 there is an identity 
that holds in  $\pi_{k,l} $ but does not hold in  $\pi_{ak,bl} $. Our goal is to write down such an identity   explicitly.  To simplify matters we will assume that $ n' $ is divisible by $ d $ (cf.  3.11.3).

Following  Remark 2.2 we will identify  generators $ A, B $ with their matrix images under irrep $ \pi_{k,l} $ (3.13). Note that $ \det(A)  = 1$ and therefore,  possible values for determinants of matrices in  $\pi_{k,l}(\Gamma_d) $ are 
$$ \exp( 2\pi i l t ), \; t = 0, \cdots, n' - 1 $$ 

\begin{proposition} 
	If $ n'$ is divisible by $ d $ (cf. (3.11.3)) then 
	\begin{enumerate}
		\item[(i)] 
	the identity 
\begin{align}
\!\! \!\!	\!\! \!\!\!\! \!\!	\!\! \!\!	\mathcal{C}_{mn} (Y_{mn}) \prod_{t=2}^{n'-1} ( \hat{\sigma}_d(y,Y_{mn}) - \exp( 2\pi i l t )) ( yx^{n}y^{-1} - x^{rn}), \; y \in Y \setminus Y_{mn}
\end{align}		
\noindent holds in irrep  $\pi_{k,l} $ (see Example 3.2 for the definition of $\hat{\sigma}_d$) 
	\item[(ii)]
	the identity (3.14) holds in  $\pi_{ak,bl} $ if and only if  $b \equiv 1  \!\!\! \mod  \! d$
		\end{enumerate}
\end{proposition}
\noindent Proof. We start with the statement  (i).  Let the variable $ x $ take a value $ g \in \pi_{k,l}(\Gamma_d) $.  Assuming that the variable $ y $ takes a value $ B^u A^v $ we have $ \det(B^u A^v) = \exp( 2\pi i l u ) $ and if the number $ \exp( 2\pi i l u )$ does not occur in one of "determinant" terms in (3.14) then $ u \equiv 1 \! \!\!\mod n' $ and hence $ u \equiv 1 \! \!\!\mod d $. In this case, using the fact that  $ g^n \in \; <A> $ we get  $ B^u(g^n)B^{-u} =  g^{nr^u} = g^{nr} $, because $ r^d = 1 \mod m $ (cf. (3.11), (3.11.1)). The conclusion is that the last term of (3.14)  vanishes under any assignment of variables $x, \; y$.

On the other hand, suppose  that the identity (3.14) holds in irrep 
 $\pi_{ak,bl} $ and from now on let the generators $ A, B $ belong to 
  $\pi_{ak,bl}(\Gamma_d) $. All possible values of the determinant
  are $ \exp( 2\pi i bl t ), \; t = 0, \cdots, n' - 1 $ and since 
  $(b,n') = 1= (l, n' ) $ there is $t $ such that $ bt \equiv 1  \! \!\!\mod n' $.  Assign the value $ B^t$ to the variable $ y $ 
for  this particular choice of $ t $ and let $A$ be the value of $x^n$. For identity (3.14) to hold, one must have $  A^{r^t} = B^{t}A B^{-t} = A^r $ (cf. (3.11)) or equivalently
$ r^{t-1} \equiv 1  \! \!\! \mod m  $. Again, recalling the condition (3.11.1) we  get  $ t \equiv 1  \! \!\! \mod d $ 
(cf. \cite{Wolf2}). Since $ bt \equiv 1  \! \!\!\mod n' $ and
$ n'$ in turn is divisible by $d$, this proves (ii).

\begin{remark}
	The identity (3.14) is constructed by the 
	recipe   (2.18*). Most probably, the same  technique can be used to handle the general case of $ n' $ not necessarily divisible by $ d $. That will require  expanding the "disjunction" (3.14) by terms  corresponding to  residue classes in  $ \mathbb{Z}_{n'} / \mathbb{Z}_d$ 
\end{remark}

\noindent We conclude this section with a well known example that emphasizes a connection between some properties of spherical space forms of type $I$ (\cite{Wolf2}, \cite{Ikeda}, \cite{Ikeda_1}, \cite{Gilkey})
 \cite{Ikeda}) and identities of corresponding irreducible representations.

\noindent \begin{example} Gassmann equivalent  but non-similar irreps  of finite groups appear in the following context:
\begin{enumerate}
	\item[1.] (\cite{Sunada}). If $ \Gamma_1, \Gamma_2 \subset \textnormal{Aut}(M) $ are Gassmann equivalent fixed point free subgroups of the group of isometries of a compact Riemannian manifold $ M$ then   manifolds $ M /\Gamma_i, \; i = 1,2 $ are isospectral
	\item[2.] (\cite{Ikeda}-\cite{Ikeda_1}, \cite{Wolf}, 	\cite{Wolf2}). If $\pi_1, \pi_2 $ are faithful real orthogonal  
 fpf representations of a finite group  $ \Gamma_d $ (3.11)  then spherical space forms $ S^{2d-1}/\pi_i(\Gamma), \; i = 1,2 $ are isospectral. They are isometric if and only if $ \pi_1$ and $ \pi_2$ are similar (cf. Lemma 3.3)
	\item[3.] (\cite{Ikeda}-\cite{Ikeda_1}, 
	\cite{Wolf}, \cite{Wolf2}) There is an infinite series (indexed by natural number $d$) of metacyclic  groups  $\Gamma_d$ (cf. (3.11) and Theorem A) of  type $I$ such that  :
		\begin{enumerate} 
			\item [3.1] 	$ \Gamma_d$ has at least  two  non-similar Gassmann equivalent  real orthogonal irreps $ \rho_1, \rho_2 : \Gamma_d \hookrightarrow \textbf{O}(2d)$ with corresponding  non-isometric but isospectral space forms   
		\item[3.2] Complexification of $\rho_i, \; i = 1,2 $ in 3.1 splits into direct sum of complex conjugate unitary irreps  $ \prescript{}{\mathbb{C}}\rho_{i} = \sigma_i \oplus \bar{\sigma_i} $. Unitary irreps $\sigma_i, \;  i =1,2$ are Gassmann equivalent (cf. also \cite{isospectral-CR}) and Galois conjugate  but they are not similar. 
			By Proposition 3.1 both irreps $ \sigma_1$ and $  \sigma_2$ satisfy identities $\mathcal{A}(\sigma_1)$  (3.7) and  $\mathcal{G}_{i}(\sigma_1 )$ (3.9). On the other hand, 
	the identity (3.14) separates irreps $\sigma_1$ and
		 $\sigma_2$  at least in some generic cases 
	\end{enumerate}
	
\end{enumerate}
\end{example}
\subsubsection{Gassmann Equivalence and Atyah-Tall Theorem}
A forthcoming question is: when Gassmann equivalent irreps are Galois-conjugate?
Trying to answer this question one notes that if representations $ \rho_i : \Gamma \rightarrow \textnormal{GL}(V), \; i = 1,2  $ of a finite group $\Gamma $ are Galois conjugate then  restricted  representations  $ \rho_i | \Sigma, \; i = 1,2 $ are also Galois conjugate for any subgroup $ \Sigma \subset \Gamma $ .
\begin{definition}
Let's call representations $ \rho_i : \Gamma \hookrightarrow \textnormal{GL}(V) $ \textit{uniformly Gassmann equivalent} if restrictions $ \rho_i | \Sigma, \; i = 1,2 $ on any subgroup $ \Sigma $ of $ \Gamma$ are Gassmann equivalent.
\end{definition}

\noindent It is probably worth mentioning that irreducible representations of finite  nilpotent groups are uniformly Gassmann equivalent if and only if they are Galois conjugate.  This fact is a direct consequence of Atyah-Tall Theorem (cf. \cite{Atyah}, \cite{JO}). We will provide some details here.       

\begin{lemma}(\cite{JSegal})	
	If  representations  
	$ \rho_i : \Gamma_i \hookrightarrow \textnormal{GL}(V), \; i = 1,2  $ of finite groups $\Gamma_1, \; \Gamma_2 $
	are Gassmann equivalent then their fixed point subspaces have the same dimension: $ \dim V^{\rho_1(\Gamma_1)} = \dim V^{  \rho_2(\Gamma_2)} $
\end{lemma}
\noindent The proof (cf. e.g. \cite{JSegal}) is a straightforward application of the Molien theorem (cf. e.g. \cite{Stanley}). 
\begin{theorema}(Molien) Let $ \phi : \Gamma \rightarrow  
	\textnormal{GL}(V) $ be a finite dimensional representation of a finite group $ \Gamma$. Then the Poincaré series of the ring  $ k[V]^{\Gamma}$ of polynomial $\Gamma$-invariants on $ V $ is given by 
	\begin{align}
		P(k[V]^{\phi(\Gamma)}, \lambda ) = \frac{1}{|\Gamma|} \sum_{g \in \Gamma} \frac{1}{\det(I - \phi(g)^{-1}\lambda)} 
	\end{align} 
and is a rational function of $\lambda$
\end{theorema}
\noindent Lemma 3.4 is now obvious - the  average of inverted characteristic polynomials is preserved by Gassmann equivalency and so are (by (3.15))  dimensions of  spaces of homogeneous polynomial invariants (cf. e.g. \cite{JSegal}) in general and  dimensions of fixed point sets in particular.
\newline
\newline
\noindent Applying Lemma 3.2 (ii) to every subgroup of the acting group, one gets the following
\begin{lemma}
	Galois conjugate representations are uniformly Gassmann equivalent 
\end{lemma}
\noindent To state the Atyah-Tall theorem we need the definition of $J$-equivalent representations (\cite{Atyah}, \cite{JO}).
\begin{definition}
	Let
$ \rho_i : \Gamma \rightarrow \textnormal{GL}(V_i), \; i = 1,2  $ be unitary representations of a group $\Gamma $. 
The  map  $$ f : S(V_1) \rightarrow S(V_2) $$
 is called ($(\rho_1, \rho_2)$-)equivariant  if
$  f(\rho_1(g)v) = \rho_2(g)f(v) $ for all $ v \in S(V_1) $ 
\end{definition}
\noindent For a continuous map  $ f : S(V_1) \rightarrow S(V_2) $
between unit spheres of the same dimension we denote by $ \deg f $
the topological degree of $f$.
\begin{definition} (\cite{Atyah})
Two	unitary representations $ \rho_i : \Gamma \rightarrow \textnormal{GL}(V_i),  \; i = 1,2  $ of a
finite group $\Gamma$  are said to be J-equivalent if $ \dim V_1 =\dim V_2 $ and there are continuous       	equivariant     
maps $ f_1 : S(V_1) \rightarrow S(V_2) $ and
$ f_2 : S(V_2) \rightarrow S(V_1) $
such that $(\deg f_1 , |\Gamma| ) = 1 = (\deg f_2, |\Gamma| )$
\end{definition}
\noindent A version of Atya-Tall theorem can be stated as follows 
\begin{theorema}(Atyah-Tall, \cite{Atyah}, \cite{JO}).	
	Let $ \rho_i : \Gamma \rightarrow \textnormal{GL}(V), \; i = 1,2  $ be faithful unitary irreps of a finite nilpotent group $\Gamma $. The following conditions are equivalent
\begin{enumerate}
	\item[(i)] Irreps $  \rho_1 $ and
	$\rho_2$ are uniformly Gasmann equivalent
	\item[(ii)] $ \dim V^{\rho_1(\Sigma)} = \dim V^{\rho_2(\Sigma)} $ for any subgroup $ \Sigma $
	of $ \Gamma  $   
	\item[(iii)] $ \dim V^{\rho_1(\Sigma)} = \dim V^{\rho_2(\Sigma)} $ for any cyclic subgroup $ \Sigma $
	of $ \Gamma  $
	\item[(iv)] Irreps $  \rho_1 $ and
	$\rho_2$ are J-equivalent
	\item[(v)] Irreps $  \rho_1 $ and
	$\rho_2$ are Galois conjugate
\end{enumerate} 
\end{theorema} 
\noindent To sum up,  implication $(i) \implies (ii) $ follows from Lemma 3.3.
The equivalence of statements (ii) - (v) is the subject of Atyah-Tall theorem (see \cite{Atyah}, \cite{JO} and references therein). The implication $(v) \implies (i) $ is the statement of Lemma 3.5

\begin{remark}
Atya-Tall theorem is actually a statement about arbitrary (not necessarily irreducible) representations   (see \cite{Atyah}, \cite{JO}). In its original more general setting, the condition (v) requires that each of the representations $\rho_i, \; i = 1,2 $ splits into direct sum of irreducible components $ \overset{t}{\underset{j=1}\oplus} \rho_{i,j}$  in such a way that $ \rho_{1,j} $ is Galois conjugate to  $ \rho_{2,j}, \; j = 1, \cdots,t $. With this modification the implication $(v) \implies (i) $ is no longer valid when $ t >1  $. For a simple counter-example, consider  sums of exact one-dimensional representations of a cyclic $p$-group.  
\end{remark}	
	
\begin{example}
	Under conditions of Atyah-Tall theorem  both  irreps $ \rho_1 $ and  $\rho_2$ satisfy  identities $\mathcal{G}_{i}(\rho_1) $ (3.9) (cf. Proposition 3.1) 
\end{example}
	
\noindent It is easy to see that Gassmann equivalence is inherited by restrictions to Sylov subgroups of  nilpotent groups. Indeed, 
let $ G = P_1 \times \cdots \times P_t $ be a finite nilpotent group decomposed into a direct product of its Sylow $p$-subgroups (cf. e.g. \cite{Serre2}). Any faithful finite dimensional irrep  
$ \rho :  G \rightarrow \textnormal{Aut}(V) $ uniquely  decomposes into a tensor product of faithful irreps of its subgroups $ P_i, \; i =1, \cdots, t $ (cf. e.g \cite{Serr}), so we have
\begin{equation}
	G = \prod_{i=1}^t P_i, \;	\rho = \otimes^{t}_{i=1} \rho_i, \;  \rho_i : P_i \rightarrow \textnormal{Aut}(V_i) \nonumber
\end{equation}
Let $ H = \prod_{i=1}^s Q_i $ be a Sylow decomposition of another nilpotent group and let
\begin{equation}
	H = \prod_{i=1}^s Q_i, \;	\sigma = \otimes^{s}_{i=1} \sigma_i, \;  \sigma_i : Q_i \rightarrow \textnormal{Aut}(W_i)  \nonumber
\end{equation}
be a tensor decompostion  of a faithful irrep  	$ \sigma :  H \rightarrow \textnormal{Aut}(W) $ 
\begin{lemma}
	Irreps $ \rho $ and $ \sigma$ are Gassmann equivalent if and only if 
	$ s = t $ and (up-to reordering)
	$ \rho_i $ is Gassmann equivalent to $ \sigma_i$ for all $ i = 1, \cdots, s $
\end{lemma}
\noindent The lemma  follows from definitions, since: (a) Gassmann equivalence is preserved by tensor products, (b) images of a $p$-Syllow subgroup in Gassmann equivalent representations are Gassmann equivalent and (c) Galois automorphism takes a $p$-group into a $p$-group. 
\newline
\newline
\noindent  
Denote by $\mathfrak{Z}_i(\Gamma) $ the $i$-th term of upper central series of a group $\Gamma$, so that $  \mathfrak{Z}_0(\Gamma) = {1}, \; \mathfrak{Z}_1(\Gamma)  $ is a center of $\Gamma, \; \mathfrak{Z}_2(\Gamma)/\mathfrak{Z}_1(\Gamma)$ is a center of $ \Gamma/\mathfrak{Z}_1(\Gamma)$ and so on.

 It is well known and it is easy to check that a finite $p$-group $P$ admits a faithful irrep if and only if 
$ \mathfrak{Z}_1(P)  $ is cyclic. 
\begin{question}
	Let $P$ be a finite non-abelian $p$-group with a cyclic center 
 and let $ \phi, \; \psi $  be  faithful Gassmann equivalent  irreps of $ P $ of the same dimension. Is it true that restricted representations 
	$ \phi |  \mathfrak{Z}_i(\Gamma) $ and $ \psi |  \mathfrak{Z}_i(\Gamma) $ are Gassmann equivalent for all $ i = 1, \; 2, \; \cdots $ ? 
\end{question} 
\noindent It is not difficult to write identities that impose a positive answer (see Section 3.1.5). 
What seems to be more interesting, however,  is the following 

\begin{proposition}
	Let $P$ be a finite non-abelian $p$-group with a cyclic center $ C  = \mathfrak{Z}_1(P)  $
	and let $ C_2  = \mathfrak{Z}_2(P) $. Let $ \phi, \; \psi $  be  faithful irreps of $ P $ that have the same dimension.
	Then 
	restricted representations $ \phi | C_2 $ and $\psi|C_2$ are Galois conjugate
\end{proposition}

\noindent We postpone the proof in order to review some additional tools. 
The following lemma is self explanatory.
\begin{lemma}(\cite{Phase}, Lemmas 22, 26).
	In notation of  Proposition 3.3,  the following holds for any $ a \in C_2 \setminus C $
	\begin{enumerate}
		\item[(i)] The subgroup $A = \; <a,C>$ spanned by $ a $ and $ C $ is a  normal abelian subgroup of $ P $
		\item[(ii)] The map  $ \mathfrak{c} : P \rightarrow C $ defined by $ x \rightarrow [x,a] $ is a homomorphism
		\item[(iii)] A centralizer of $a $ in $ P$ is a normal subgroup $ P_a = \ker \mathfrak{c}  $.  Therefore, $ C \ni  [x,a] \neq 1 $ 
		for any $ x \in P \setminus P_a $. 
		\item[(iv)] The factor-group $P/P_a $ is cyclic. Let $ x \in P \setminus P_a $ be any preimage in $ P $ of a generator of $ P/P_a$. Then  $ P/P_a \approx \; <z> \; \approx  A/C $, where  $ z = [x,a] \in C  $. 

	\item[(v)] Let 
	$\pi : P \hookrightarrow \textnormal{Aut}(V) $  be a  faithful irrep of  $P$ in a  vector space  $V$. Then there are exactly $q$  isotypical components $ W_1, \cdots, W_q $ of the restricted representation $ \pi | A $.  Representations $  \pi_{a,i} : P_a \rightarrow \textnormal{Aut}(W_i) $ are irreducible and 
	$ \pi = \textnormal{Ind}_{P_a}^P \pi_i $ for all   $ i = 1, 2, \cdots,  q $.    
	\end{enumerate}
\end{lemma}
\begin{lemma}
	Under conditions of Lemma 3.7,  set  $ q = |P/P_0| = |A/C| = | <z> | $ and let $ \pi(z) = \lambda I, \; \lambda \in \mathbb{C}^* $. The following holds:
\begin{enumerate}
\item[(i)] $ n = \dim \pi $ is divisible by $ q $
 \item[(ii)] 	
$ \sigma_i(\pi(a^t)) = 0  = \sigma_i(\pi(x^t)) $ (cf. (3.1)) for  all integers $i, \; t $ not divisible by $q$ 
\item[(iii)] In particular,	$ \chi_{\pi}(a^t) = 0 = 
\chi_{\pi}(x^t) $  
for  all integer $t $ not divisible by $q$ 
\item[(iv)] and therefore characteristic polynomials $ p_a(T), \; p_x(T) $ of $ \pi(a) $ and $ \pi(x) $ are polynomials of $ T^q $ 
 \end{enumerate}	
\end{lemma}
\noindent Let's identify the group  $P$  with its  $\pi$-image in $ \textnormal{GL}(V)$. 

Proof of statement (i).  One has $ 1 = \det (z) = \lambda^n $
because $ z = [a,x] $ is a commutator. Hence, $n$ is divisible by the order $q $ of $z$ in $C$. 

Proof of statement (ii).
 By Lemma 3.7 (iv), $ z^t \neq 1 $ for all $ t = 1,\;\cdots , \; q-1$. For integers $t$ in this range we have $ xa^tx^{-1} = a^t z^t, \; z^t \neq 1 $. 
  Remembering that $ z = \lambda I$  is a scalar matrix, evaluate polynomial $\sigma_i$ on both sides of this equation. 
  The homogeneous polynomial $ \sigma_i $  has degree $i$ and  is invariant under adjoint action of the group, so we get $ \sigma_i(a^t) = \lambda^{it} \sigma_i(a^t) $. Hence $ \sigma_i(a^t) = 0 = \sigma_i(x^t) $ for all $ i,t $ not divisible by $ q $

\begin{lemma}
	If	in addition to  conditions of Prposition 3.3 the irreps $\phi $ and $\psi$ coincide on the center, i.e.
	$ \chi_{\phi} | C = \chi_{\psi} | C $ then   $ \phi(a) $ and  $ \psi(a) $ have the same spectrum for any $ a \in C_2$. In particular, $ \chi_{\phi} | C_2 = \chi_{\psi} | C_2 $ and representations
	$ \phi |C_2 $ and $ \psi | C_2 $  are equivalent
\end{lemma}
\noindent  Indeed, by Lemma 3.8 we must have
$ \chi_{\phi}(a^t)  =  \chi_{\psi}(a^t) $ for all integer $ t $. Therefore, characteristic polynomials, and  hence spectra, of $ \phi(a) $ and $ \psi(a) $ do coincide.   

\paragraph{Proof of Proposition 3.3} Since $C$ is cyclic   there is a Galois automorphism $\epsilon$ such that $ \epsilon \phi | C =
   \psi | C $ and we can reference Lemma 3.9 to finish the proof.   
\newline
\newline
\noindent Bearing in mind Lemma 3.8 we also have
\begin{corollary} (Cf. \cite{Drinfeld} Appendic C).
	The invariants  $ \sigma_i, \; i = 1, \cdots $ (cf. (3.1)) and in particular a character of an excat irreducible representation of a finite two-step nilpotent group $\Gamma$
	are identically zero on $ \Gamma \setminus   \mathfrak{Z}(\Gamma) $ 
\end{corollary}
 
\begin{corollary} (Cf. e.g. \cite{Stone}, \cite{Howe}, \cite{Heisenberg}, \cite{Bump}).
	Faithful irrep $\theta$ of a finite two-step nilpotent group $\Gamma$ is determined by its restriction $ \theta |  \mathfrak{Z}(\Gamma) $  
\end{corollary}   

\begin{corollary} 
	Same dimension faithful irreps of a finite two-step nilpotent  group are Galois conjugate and therefore J-equivalent and uniformly Gassmann equivalent  
\end{corollary}

\begin{example} Disjunctive trace identity $ [x,y] = 1 \;\vee \; Tr(x) = 0 $ holds in any faithful finite dimensional irrep of a two-step nilpotent group. This is just a reformulation of Corollary 3.8. Moreover, it is easy to see that disjunctive trace identities 
	$$ [ [x_1,x_2], \cdots ] , x_t] = 1  \;\vee \; Tr(x_t) = 0 $$ 
	$$ [ [x_1,x_2], \cdots ] , x_t] = 1  \;\vee \; Tr([ [x_1,x_2], \cdots ] , x_{t-1}]) = 0 $$
	hold in any faithful finite dimensional irrep of a $t$-step nilpotent group ($t > 2 $)
\end{example} 
\begin{example}
	Take a prime number  $p>2$ and let $H_p$ be a group of upper unitriangular matrices in
	$ \textnormal{GL}_3(\mathbb{Z}_p) $
	  This is a minimal Heisenberg $p$-group of order $p^3 $  (cf. e.g. \cite{semidirect}).  Exact  irreps of this group are induced from normal abelian subgroups of order $p^2$ (e.g. by Lemma 3.7 (v)) and are  Gallois conjugate to each other by Corollary 3.9. Any automorphism of $\mathfrak{Z}(H_p)$ can be extended to an automorphism of $H_p$ and therefore all exact irreps of $ H_p$ are  similar as well. Thereby these irreps have the same identities. Let $ \theta_0$ be one such  $p$-dimensional (exact) irrep of  $ H_p$. Clearly $ \theta_0$  satisfies       the following  list of identities
	  	\begin{enumerate}
	  		\item[(a)] exponent identity $ x^p - 1$ 
	  		\item[(b)] order identity  $	\mathcal{C}_{p^3+1} $ (cf. (1.3))
	  		\item[(c)] standard identity $s_{2p} $ (cf. Example 1.2)  
	  	\end{enumerate}

\noindent Let's check that all identities of $\theta_0$  follow from identities (a)-(c). Let $\theta$ be a finite dimensional exact irrep of a finite group  $ \Gamma $ that satisfies identities (a)-(c).  According to standard methodology (cf. \cite{P}-\cite{PK}) we need to check  that 
\begin{enumerate}
\item[(*)] there is  a direct summand of a restriction of $ \theta_0$ unto some subgroup of $ H_p $ that modulo its kernel is similar to $ \theta$  
\end{enumerate}
\noindent Indeed, it follows from (b) that $ |\Gamma| \leq p^3$
and it follows from (a) that $ \Gamma $ is a $p$-group of exponent $p$. Therefore, $ \Gamma$ is either  abelian or is isomorphic to $H_p$. In the latter case $ \theta $ is similar to $\theta_0$ as was explained above. In the former case the irrep $\theta$ is a one-dimensional representation of a cyclic group of order $p$ that clearly satisfies the condition (*).     
\end{example}
    
\noindent Our next example shows that Gassmann equivalent irreps of nilpotent groups of class greater than two  are not necessarily Galois conjugate 
\subsubsection{Gassmann equivalence does not imply Galois conjugacy for irreps of finnite $p$-groups}
Fix a prime number $p, \; p > 2 $ and let $ A = \mathbb{Z}^{p}_p  $ be   a direct sum of $ p$  copies the cyclic $p$-group 
$ \mathbb{Z}_p$.
Denote by $ \pi : S_p \rightarrow  \textnormal{GL}_p(\mathbb{Z}_p) \approx \textnormal{Aut}(A)  $ coordinate index permutation action of the symmetric group $S_p$, that is   
\begin{equation}
\pi(\tau )a = ( a_{\tau(1)}, \cdots , a_{\tau(p)}) , \; a = (a_1, \cdots , a_p ) \in A, \; \tau \in S_p 	
\end{equation} 
In particular, we fix an action of the cyclic group $ C_p$ on $ A$ by mapping a generator of $C_p$ into a cycle $ \sigma = (1,2,\cdots, p) $. Let $ \Gamma = A C_p   $ be a semidirect product (cf. e.g \cite{Serre2}) that arises from this action. Recall that multiplication in $  \Gamma$ is defined by the rule 
$$ (a_1\sigma_1)(a_2\sigma_2) = (a_1 + \pi(\sigma_1)(a_2)) \sigma_1 \sigma_2, \; a_1,a_2 \in A, \; \sigma_1, \sigma_2 \in C_p     $$
\begin{lemma}
	If  automorphism $ g \in \textnormal{Aut}(A) $ commutes with $ \pi(\sigma) $ then it extends to an automorphism of $ \Gamma$
\end{lemma}
\noindent Proof. Define a map $ g: \Gamma \rightarrow \Gamma $ by the rule $ g(a\tau) = g(a) \tau $. Then 
\begin{gather}
 g( (a_1\sigma_1)(a_2\sigma_2) ) = g( ( a_1 + \pi(\sigma_1)(a_2) ) \sigma_1 \sigma_2 ) = 
g( ( a_1 + \pi(\sigma_1)(a_2) ) \sigma_1 \sigma_2 =  \nonumber \\
= (ga_1) \sigma_1\sigma_2 + (g\pi(\sigma_1)(a_2))\sigma_1\sigma_2 
= g(a_1)\sigma_1 \sigma_2 + (\pi(\sigma_1)g(a_2)) \sigma_1\sigma_2 = \nonumber  \\
= (g(a_1)\sigma_1) (g(a_2)\sigma_2) = g(a_1\sigma_1) g(a_2\sigma_2)
\nonumber
\end{gather}

\noindent A generic linear character $ w' : A \rightarrow \mathbb{C}^* $ can be defined as follows (cf. e.g \cite{Fourier}). Take a linear form  $  w : A \rightarrow \mathbb{Z}_p $
represented by a vector $(w_1, \cdots, w_p) \in \mathbb{Z}^p$, via the standard scalar product $$ w(a) \; = \; <w,a> \;= \; 
\sum_{i=1}^p w_i a_i, \; a \in A$$
and set 

 \begin{equation}
  \omega'  (a ) = \exp \left(\frac{2\pi i}{p} <w, a> \right)
 \end{equation} 
\noindent  To be precise, the coordinates $ a_i, w_i, \; i = 1, \cdots p $ on the right hand side of (3.17) are understood as integers that correspond to their residuals in $ \mathbb{Z}_p$. This short-cut should not cause any confusion.   
\begin{remark} (see also examples in appendix 2). 
	 Non-trivial characters  (3.17) are all similar to each other.  
Indeed, for  $ \omega_1 \neq 0 \neq \omega_2 \in \mathbb{Z}^p $ there is a matrix $ g \in \textnormal{GL}_p(\mathbb{Z}_p) \approx \textnormal{Aut}(A) $ such that  $ \omega_2 = g^T \omega_1 $ and we have 
$$ w'_2 = w'_1 g \Longleftrightarrow  w'_1(ga) = w'_2(a)
\Longleftrightarrow \; < g^Tw_1, a > \; = \; <w_2,a>
, \; a \in A $$
\end{remark} 

\noindent  Let  $ \rho_w = \textnormal{Ind}^{\Gamma}_{A} \omega'    : \Gamma \rightarrow \textnormal{Aut}(V) $ be  a representation induced by  $w'$. 
It is easy to see that in the basis  $ e_1 , \; e_2 = \sigma e_1, \; \cdots , \; e_p = \sigma^{p-1} e_1  \in V  $ the  matrix of $ \rho_w(a), \; a \in A $  is
\begin{gather}
		\textnormal{diag}  \!	\left( \; 
	\exp \left(  \gamma_p <w,a> \right ) , \; 
	\exp  \left( \gamma_p < \sigma_{1}^T w,a \!> \right), \cdots, 
	\exp  \left( \gamma_p < \sigma_{p-1}^Tw,a > \right)\!
	\; \right)  \nonumber = \\
	\!\!\!\!\! =  \;\; 	\textnormal{diag} \!	\left( \:
	\exp \left( \gamma_p <w,a> \right),\; \cdots,
	\exp  \left( \gamma_p <w,\sigma_{p-1}a \!> \;\right)
	\right)    
\end{gather} 
\noindent where $ \gamma_p = 2\pi i / p $ and 
$ \sigma_i = \sigma^i, \; i = 1,  \cdots p-1 $.

Denote by $ \tilde{w} $ the circulant matrix that has rows  rows $ w, \; w\sigma_1, \; \cdots, \; w \sigma_{p-1} $. Using a shorthand
\begin{align}
 \exp(v) = (\exp(v_1), \cdots, \exp(v_p) )) , \; v =(v_1, \cdots,v_p) \in \mathbb{C}^p  \tag{3.18'}
\end{align}
we can rewrire (3.18) as
\begin{equation}
 \rho_w(a) = \exp( \gamma_p \tilde{w}a) 
\end{equation}

\begin{lemma}  $\rho_w $  is irreducible and exact if and only if $ \; \sum_{i=1}^p w_i \neq 0 $   	
\end{lemma}
\noindent Proof. If $ \sum_{i=1}^p w_i = 0 $ then $\rho_w $ is trivial on $ \mathfrak{Z}_1(\Gamma ) $ and the "only if" part of the Lemma follows from the definition (3.17).
   
 Suppose that$ \; \sum_{i=1}^p w_i \neq 0 $.    Let $ a^{\gamma} \equiv \gamma^{-1} a \gamma, \; \gamma \in \Gamma   $ denote the (adjoint)  action of $ \Gamma $ on $A$ and let  $ H \equiv H_{w'} = \{  h \; | \;  w'(a^h) = w'(a), \; h \in \Gamma , \; a \in A\} $  be the group of inertia of the character $w'$. It is easy to see (cf. e.g. \cite{Serr} or \cite{semidirect})  that representation $ \textnormal{Ind}^{\Gamma}_{H} w' $ is irreducible. Note, however,  that conditions of the our lemma imply that  $ H = A $. Indeed, if $ w_1 = w_2 = \cdots = w_p $ then   $ \; \sum_{i=1}^p w_i = 0 $, so $ w_i \neq  w_j $ for some pair of indices of $ 1 \leq i < j \leq p $. Let $ \delta_i \in A \approx \mathbb{Z}^{p}_p $ denote a vector with $i$-th coordinate equal to $1$  and all other coordinates equal to zero. It is clear that $$ w_i = \; <w,\delta_i> \; \neq  \;
   <w,\delta_{i}^{\sigma^{j-i}}> \;= \; <w, \delta_j > \;= w_j $$ and therefore if the character $ w' $ is invariant under the adjoint action of  $\sigma^k \in \Gamma $ then $ k = 0 $. Hence $ H = A $,  $ \rho_w $ is irreducible and $ \dim \rho_w = p $.
Suppose that  $  a \in \ker \rho $ for some $ a \in \Gamma $.
Since any normal subgroup  intersects the center of a nilpotent group non-trivially (cf. e.g. \cite{Serre2}) we can assume without loss of generality that
$ a \in \mathfrak{Z}_1( \Gamma ) \cap \ker \rho $ in which case (cf. (3.16))
 all coordinates of $a$ must be equal,  i.e $ a_i = a_1, \; i = 2, \cdots,  p $. It follows that $ <w,a> \; = \; a_1 \sum_{i=1}^p w_i = 0 $ as it is assumed  that $ a \in \ker \rho$. Conditions of the Lemma now imply  that $ a_1 = 0 $ and we must conclude that  $ \ker \rho  = \{1\} $. 
\newline
\newline
\noindent Call a linear form  \textit{admissible} if it satisfies the condition of Lemma 3.11. 
\newline\newline
\noindent The matrix of $ \pi(\sigma) $ (3.15)  in $ \textnormal{GL}_p(\mathbb{Z}_p) \approx \textnormal{Aut}(A) $  looks as follows
\begin{gather}
\pi(\sigma) =
\begin{bmatrix}
	0 & 1 & 0 & \cdots & 0 \\
	0 &  0 & 1 & \cdots & 0 \\ 
	\cdot & 	\cdot  & \cdot  & \cdots &  \cdot \\
	\cdot & 	\cdot  & \cdot  & \cdots &  \cdot \\
	 0 & 0 & \cdot & \cdots & 1 \\
	1 & 	0 & \cdot
	&  \cdots & 0
\end{bmatrix}     
\end{gather}

\noindent Extending the homomorphis $ \pi $ (3.16) by linearity,  set $$ M(\sigma) = \pi(\mathbb{Z}_p C_p ) \subset 
	\textnormal{M}_p(\mathbb{Z}_p ) \approx 	\textnormal{End}(\mathbb{Z}_p^{p}) \approx \textnormal{End}(A) $$
This is a well known algebra of \textit{circulant} matrices
(see e.g. \cite{circulant} and references therein).	
	 Few useful facts  about the algebra $M(\sigma)$ are listed below. 
\begin{lemma}\
 \begin{enumerate}
 	\item[(i)] $ M(\sigma) $  is isomorphic to the group algebra 
 	$ \mathbb{Z}_p C_p $
 
 	\item[(ii)] $M(\sigma)$ is the set of all circulant matrices
 	\item[(iii)]   $ M(\sigma) $ is a centralizer of $\pi(\sigma) $ in the full matrix algebra $\textnormal{M}_p(\mathbb{Z}_p ) $
 	\item[(iv)] There is a unit  $ h \in M(\sigma)^{*} \setminus \pi(C_p) $  of order $ p $
 	\item[(v)] A form $w$ is admissible if and only if its circulant $ \tilde{w}$ is invertible 
 	\item[(vi)] As an element in $M(\sigma) $ the circulant matrix of a form $w =  (w_1, \cdots, w_p) \in \mathbb{Z}^p $ can be expressed as   
 	
 \begin{align}
 	\tilde{w} = \sum_{i=1}^{p} w_i \pi(\sigma^{i-1})  \tag{3.20'}
 \end{align} 
 \end{enumerate}  
\end{lemma}
\noindent Proof.  Let $ \mathfrak{i}(M) $ be the set of index pairs of non-zero entries of a matrix $M $.  
Due to the shape of the matrix  (3.20), we have:  (a)  $ \mathfrak{i}(\pi(\sigma^{i} ) \cap  \mathfrak{i} (\pi(\sigma^{j} )) = \varnothing $ if $ i \neq j$, and (b)
 $ M \in M(\sigma) $ if and only if the  value of an entry $ M_{i,j}  $ depends only on  the residual $ (i - j) \mod p $.
  Now (i)  trivially follows from (a) while (ii) and (iii)  follow from (a) and (b) since $ \sigma^{-1}  M \sigma = M $ if and only if $ M $ satisfies (b).
 
 Further, it is well known that the group of units of the algebra  	$ \mathbb{Z}_p C_p $ is an elementary abelian $p$-group  of rank $ p -1 $ (cf. \cite{Sandlinng}, \cite{Bovdi}). Hence, (iv) follows from (i) as we assumed from the start that $ p > 2$.
 
  The formula (vi) directly follows from (3.20).
  
   Let  $ \Delta $  be the kernel (augmentation ideal) of the homomorphism $    \mathbb{Z}_p C_p  \rightarrow  \mathbb{Z}_p $ that sums up the  of elements of the group algebra. 
Recall that the ideal $ \Delta$ is maximal and nilpotent 
(cf. e.g.  \cite{Sandlinng}, \cite{Bovdi}). Hence, its complement  	$ \mathbb{Z}_p C_p \setminus \Delta  $ is a multiplicative set of invertible in $	 \mathbb{Z}_p C_p $ elements and we see that (v) follows from (3.20')
 \begin{lemma}
 Let $ g \in M(\sigma)    $ and let $w$ be an admissible form   in $  \mathbb{Z}^p $. Then  
 \begin{enumerate}
 	\item[(i)] if $ g $  is invertible in $M(\sigma) $ then the form $gw $    is also admissible
 	\item[(ii)] $ \widetilde{gu}= g\tilde{u} $ and  $ g \tilde{u} = \tilde{u} g $ for any form $ u \in  \mathbb{Z}^p $ 
 \end{enumerate} 
 \end{lemma}
 \noindent Proof.  Statement (ii) follows from commutativity of $ M(\sigma)$,  Lemma 3.12 (iii) and (vi). In turn, (i) follows from  (ii) and Lemma 3.12 (v)

\begin{proposition}
	If $ h \in M(\sigma)$ satisfies conditions of Lemma 3.12 (iv) then exact irreps
	$ \rho_w $ and $ \rho_{hw} $ are similar and therefore Gassmann equivalent (Lemma 3.2(i)). However, these irreps are not Galois conjugate   
\end{proposition}
\noindent Proof. By lemma 3.13 we have (cf. 3.19)
\begin{equation}
	\rho_w(ha) = \exp( \gamma_p \tilde{w}ha) = \exp( \gamma_p \widetilde{hw}a) =  \rho_{hw}(a), \; a \in A   \nonumber
\end{equation}
\noindent meaning that automorphism $h$ induces similarity between $ \rho_{w} | A $ and  $ \rho_{hw} | A $. 
By Lemma 3.10, $ h \in  \textnormal{Aut}(A) $ trivially  extends to an automorphism of $ \Gamma $ and therefore the irreps $ \rho_w $ and $ \rho_{hw} $ are similar.	
	Suppose that $ \rho_w $ and $ \rho_{hw} $ are Galois conjugate. Then  $ \rho_{w} | A $ and  $ \rho_{hw} | A $ are also Galois conjugate and there is a Galois automorphism $ \epsilon \in \mathfrak{G}(\Gamma)  $ such that representations  $ \epsilon (\rho_{w} | A ) $ and  $ \rho_{hw} | A $ are equivalent.
	In other words, there are permutation $ s \in S_p$ and an integer   $t \in \mathbb{Z}_p, \; t \neq 0$   such that
 \begin{equation}
 \exp(  \gamma_p s\tilde{w}a) = \exp( t \gamma_p \widetilde{hw}a), \; a \in A 
 \end{equation}  Applying Lemma 3.13 once again ($\widetilde{hw}= h\tilde{w}$), rewrite  (3.21) as
\begin{align}
	s = t \tilde{w}h \tilde{w}^{-1} \Longrightarrow s = th \tag{3.21'}
\end{align} 
Hence $ s \in M(\sigma)$ and, raising both sides (3.21') to a power $p$, we see that $ t = 1$. Therefore $s = h$ is a full cycle that commutes with $\sigma$, hence $ s \in \pi(C_p)$.  This, however, is in contradiction with  $ h \in  M(\sigma)^{*} \setminus \pi(C_p) $ (cf. Lemma 3.12 (iv)) 

\begin{corollary}
If  linear forms $ w_1, \; w_2 $ are admissible  then irreps 	$ \rho_{w_1} $ and $ \rho_{w_2} $ are similar. 
\end{corollary}
\noindent This fact follows from the proof of Proposition  3.4 and from the following 

\begin{corollary}
The group of units $M(s)^{*} $ acts transitively on the set of admissible linear forms  
\end{corollary}
\noindent 
Proof. If $w$ is an admissible linear form, then $ w = w_0 \tilde{w} $  where $ \tilde{w} \in M(s)^{*}  $ and $ w_0 = (1, \; 0, \; \cdots, \; 0 )$
\newline\newline
We thus have the following (counter-)example
\begin{example}
	If $ p > 2 $ then all faithful irreps of the semidirect product $ \Gamma = \mathbb{Z}_p^{p} C_p$ are similar and therefore Gassmann equivalent. However, there are faithful irreps of $ \Gamma $ that are not Galois conjugate and therefore not uniformely Gassmann equivalent. The group $ \Gamma$ is a \textit{wreath product} of two cyclic groups of order $p $ and according to \cite{Liebeck}, the nilpotence class of $ \Gamma $ is exactly  $p$. Hence, corollary 3.9 is not true for nilpotent groups of class 3.
\end{example}

\begin{conjecture}
	Let $P$ be a finite non-abelian $p$-group with a cyclic center and let $ \phi, \; \psi $  be  faithful Gassmann equivalent irreps of $ P $.
Then there is a Galois automorphism $ \epsilon$ and an automorphism $ \alpha $ of $P$ such that
irreps $ \epsilon\phi \alpha $ and $  \psi $ are equivalent
\end{conjecture} 

\subsubsection{Identities that impose Gassmann equivalence on upper central series  } 
 We need a slight modification of the construction of  identities $	\mathcal{G}_i(\rho)$ (3.9). Let $ S $ be a subset of $ \rho(G) $. Intersecting $S$ one-by-one with 
sets of the partition  $ \mathfrak{A}(\rho) $ (3.5.0) and removing empty sets if necessary, we get a partition of $S$  
\begin{align}
 \mathfrak{A}(\rho, S) = G'_{1}  \cup \cdots \cup  G'_{r'}, \;\; 0 < r' \leq r
\end{align}    
Denote by $ \mathcal{G}_i(\rho,X_s), \; i = 1, \cdots, r'; \; s =|S|  $ the expression (3.9) obtained from the partition $ \mathfrak{A}(\rho, S)$ (3.22) instead of the partition (3.5.0), where as before  $ X_s = \{x_1, \cdots x_s \} \subset Y $ is a set of free variables. Set  $ S = \rho(G) \setminus \rho(\mathfrak{Z}_t(G)) $ and    let $ U_{x_i} = \{ u_{i,1}, \cdots, u_{i,t} \}, \; i = 1, \cdots, s $ be disjoint sets of  variables in $Y$.   
Assuming that $Y_t=\{y_1,y_2, \cdots, y_t\} \subset Y $, define a shorthand $ c_t(x,Y_t) = c_t(x,y_1,y_2, \cdots, y_t) $ for a length $t+1$ commutator $ [[[x,y],y_2],\cdots , y_t] $. Let    
\begin{align}
   \mathcal{G}_{i}(S) =   \mathcal{G}_i(\rho,X_s) \prod_{x \in X_s} 
 (c_t(x,U_{x}) - 1) v_x, \;\; i = 1, \cdots, s  	
\end{align}
\noindent where as usual $ v_x \in Y $ are additional variables.
\begin{proposition}
	All identities (3.23) hold in $\rho$.
	If  irrep $\sigma$ of $G$ is Gassmann equivalent to $\rho$ then all identities (3.23)  hold in  $\sigma$  if and only if  restricted representations 
		$ \rho |  \mathfrak{Z}_t(\Gamma) $ and $ \sigma |  \mathfrak{Z}_t(\Gamma) $ are Gassmann equivalent
\end{proposition}
\noindent Proof. Let's show first that (3.23) holds in $\rho$. Due to the guard term $\mathcal{C}_s (X_s)$ (that is a multiple of  $ \mathcal{G}_i(\rho,X_s) $, cf. (3.9)) we can assume that all the variables in $ X_s $ take pairwise distinct values in $ \rho(G)$. If one of these values happens to belong to  $ \rho(\mathfrak{Z}_t(\Gamma) ) $ then corresponding commutator term in (3.23) vanishes. If, however, $ \rho(X_s) \subset \rho(G) \setminus \rho(\mathfrak{Z}_t(G)) $ then   vanishes the term $ \mathcal{G}_i(\rho,X_s) $.
 
The proof of the second part of the Proposistion 3.5 is completely analogous. If   $ \mathcal{G}_i(\rho,X_s) $  holds on any pairwise distinct set of elements in $ \sigma(G) \setminus \sigma(\mathfrak{Z}_t(G)) $ then the sets $ \sigma(G) \setminus \sigma(\mathfrak{Z}_t(G)) $ and  $ \rho(G) \setminus \rho(\mathfrak{Z}_t(G)) $ are Gassmann equivalent (cf. Remark 3.1). Since we have assumed that $ \rho(G) $ is Gassmann equivalent to $\sigma(G) $ the same must be true for $  \sigma(\mathfrak{Z}_t(G)) $ and  $\rho(\mathfrak{Z}_t(G)) $  

\begin{remark}
	If $ \mathfrak{Z}_t(G) = \{1\} $
	then $  \mathcal{G}_i(\rho,X_s) \equiv  \mathcal{G}_i(\rho,X_m) $.  In this case (3.23) does not differ from (3.9)  
\end{remark}

\section{Non-character identities}
\subsection{Minimal Polynomial Identity}
 Let $ \phi : \Gamma \rightarrow \textnormal{GL}(V) $ be a finite dimensional representation of a finite group $\Gamma$. Denote by $ \textnormal{eig}(\phi(g))  $ a set of distinct eigenvalues of $ \phi(g), \; g \in G $ and let $$ Eig(\phi) = 
 \bigcup_{g \in G } eig(\phi(g)) $$
 
 \noindent Then representation $\phi$ satisfies one-variable polynomial identity 
\begin{align}
  	P(\phi)  = \prod_{\lambda \in Eig(\phi)} (x - \lambda) 
\end{align}
\noindent The identity (4.1) can be refined in the following way. Let $ eig(\phi) $ be an ordered by inclusion poset of eigenvalue sets $  \textnormal{eig}(\phi(g)),  g \in \Gamma $ and  let $ \{\Lambda_1, \cdots \Lambda_l  \} $ be a collection of maximal elements in  $ \textnormal{eig}(\phi)$. Set 
\begin{align}
	 P_i(\phi)(x) \equiv P_i(x) = \prod_{\lambda \in \Lambda_i}(x-\lambda), \;  i = 1, \cdots, l   \nonumber 
\end{align} 
and let 
\begin{align}
	\text{M}(\phi)  = P_1(x) v_1 P_2(x) v_2 \cdots v_{l-1} P_l(x)  	
\end{align}
\noindent where  $ x, v_1, \cdots, v_l \in Y \subset F(Y)  $ are (pairwise distinct) free variables. 	
We have by design
\begin{proposition}\
\begin{enumerate}
\item[(i)]
	The representation $\phi$ satisfies identities $ P(\phi) $ (4.1) 
	and $ M(\phi) $ (4.2) 

\item[(ii)]
The representation $ \phi' : \Gamma' \rightarrow \textnormal{GL}(V') $ satisfies identity $ P(\phi) $ (4.1) if and only if
$ Eig(\phi') \subset Eig(\phi) $
\item[(iii)] Suppose that  representation $ \phi' : \Gamma' \rightarrow \textnormal{GL}(V') $ is exact and irreducible. Then the identity $ M(\phi) $ (4.2)  holds in $ \phi' $ if and only if 
 for any $ h \in \Gamma' $ there is $ g \in \Gamma $ such that
	 $ eig(\phi'(h)) \subset eig( \phi(g)) $ 
\end{enumerate}
\end{proposition}
\begin{definition} (\cite{P}) We will say that an exact irrep $ \sigma   $ of a group $ H $ is an (irreducible) \textit{factor} (factor-representation) of a representation $ \rho $  of a group $ G $ if the following conditions hold:
\begin{enumerate}
	\item[(a)] $ H $ is a section of $ G $, i.e. $ H \approx H'/K $ where $ H' $ is a subgroup of $ G $ and $ K $ is a normal subgroup in $ H '$  
	\item[(b)] there are  $H'$-invariant subspaces $ V'' \subset  V' \subset V $ 
	such that the action of $K$ on $ V'/V''$ is trivial and corresponding action of $ H'$ on $ V'/V'' $ is similar to $\sigma$ 
	
\end{enumerate}	  
\end{definition}
\noindent By definition any representation identity holds in any of its factors.  
\begin{question}
	Suppose that $ H $ is a section of $ G $ and let $\dim \sigma \leq \dim \rho $.  
\begin{enumerate}
	\item[(a)] Is it true that the inclusion $ Eig(\sigma) \subset Eig(\rho) $ implies that $ \sigma $ is a factor of $ \rho$ ?
	\item[(b)] Same question as (a) assuming in addition that  for any $ h \in H $ there is $ g \in G $ such that
	$ eig(\sigma(h)) \subset eig( \rho(g)) $
\end{enumerate}	

\end{question}
\begin{lemma}
	Suppose that irrep $ \sigma$ is one-dimensional and that $  Eig(\sigma) \subset Eig(\rho) $. Then $ \sigma$ is  a factor of $ \rho$
\end{lemma}
\noindent Proof. If $ \dim \sigma =1 $ then $ H = \; <h> $ is cyclic. Let $ \sigma(h) = \lambda \in \mathbb{C} $ and take $ g \in G $ such that one of the eigenvalues of $ \rho(g) $ is $\lambda$. If $ v \in V $ is a correponding eigenvector then representation $ \rho \; | <g> $ has an irreducible component $\rho'$  defined by  $ \rho'(g)v = \lambda v $.
\newline
\newline
\noindent Take a finite subgroup  $ \Gamma $ of $ \textnormal{SU}(2) $ and suppose that natural representation 
$ \theta: \Gamma \hookrightarrow \textnormal{SU}(2) $ is irreducible.
Most probably, Lemma 4.1 could be used to completely describe identities of such representations.  We will restrict ourselves to a  case of a binary tetrahedral group (a double cover of alternating group $A_4$, cf. \cite{binarytetrahedral}). This group is a   semidirect product of a quaternion group $ Q $ and a cyclic group $C_3 $ of order $3$.

\begin{example}
	All identities of the natural two dimensional irrep $ \theta : \Gamma = QC_3  \hookrightarrow \textnormal{SU}(2)  $ of binary tetrahedral group are consequences of the following list
	\begin{itemize}
		\item[(1)] disjunctive identities $\mathfrak{D}(\Gamma)$ (cf. Lemma 1.1) and Remark 1.3
		\item[(2)] standard identity $s_4$ (cf. Example 1.2)
		\item[(3)]  character identity $ \Psi(\theta) $ (2.4)
		\item[(4)] eigenvalue identity $ P(\theta) $ (4.1) 
	\end{itemize} 
Indeed, let $ \theta' : \Gamma' \rightarrow \textnormal{GL}(V) $ be a faithful  irrep that satisfies identities (1)-(4). As in Example 3.7 we need to check (cf. \cite{P}) that  $\theta'$ is a factor of $ \theta$.
It follows from (1) that $\Gamma' $ is a section of $ \Gamma $ and it follows from (2) that $\dim V \leq 2 $. A proper section of $ \Gamma $ is either (a) abelian,  or (b) group of quaternions Q, or (c) alternating group $ A_4$. 

In case (a) $\theta'$ is a factor of $\theta $ by (4), Proposition 4.1 and Lemma 4.1. Case (c) can be discarded because $ A_4$ does not have exact two-dimensional irreps. 
In case  (b) we note that the group of quaternions has only one irrep of dimension two.

Assuming therefore, that $ \Gamma' \approx \Gamma $ we see that (3) implies equivalence of $ \theta' $ and $\theta$ (cf. Lemma 2.2), because character values of the natural irrep of $\Gamma$ differ from character values of other $\Gamma$-irreps  (cf. e.g. \cite{binarytetrahedral}).          
\end{example}

\begin{example}
	It is easy to see (cf. e.g. \cite{GW}) that  set of eigenvalues of the tetrahedral representation $ \tau $ of the alternating group $A_4 $ is $$  \{ \; \{ 1,-1 \}, \;
	 \{ 1, e^{2\pi i/3}, e^{4\pi i/3} \}  \; \} $$  and therefore the minimal polynomial identity $ M(\tau) = (x^{2} - 1) v (x^3-1) $ is a consequence of the disjunctive formula $ x^2 = 1 \; \vee \;  x^3 =1 $ satisfied by $ A_4 $. On the other hand, minimal polynomial identities  $ M(\rho_4) $ and $ M(\rho_5)$ (cf. Example 2.9)  contain terms $ (x \pm 1) (x+i)(x-i)$ and  neither of these identities can be derived from disjunctive identities of $ S_4$ (note that the cycle of length $4$ has eigenvalues  $ \pm 1,\;  \pm i$ in  $\rho_4, \rho_5 $ respectively)  
\end{example}

\subsection{Central Partitions}
As before (see the Introduction), consider a partition $P$ of the group $G$ into conjugate classes
\begin{align}
	P = (C_1, \; \cdots, C_s ), \;\; |G| = m = \sum_{i=1}^s t_i, \; t_i = |C_i|, \; \chi_{\rho}(C_i) = \chi_i     
\end{align} 
For any set of free variables $ X'  \subset Y $ define   
\begin{align}
	  kF(Y) \ni Av(X') = (1/|X'|)\sum_{x \in X' } x \tag{4.3'}
\end{align} 
Let $ \mathcal{P} = \mathcal{P}(X_m) $ be the set of all  partitions of $X_m$ of the same type as  (4.3), i. e. 
\begin{align}
	( X_1, \cdots X_s ) \in \mathcal{P} \; \iff \; X_i \subset X_m, \;   |X_i| =t_i, \;  i = 1, \cdots, s 
\end{align}  
\noindent Set
\begin{align}
\mathfrak{P}(\rho, P) =	\mathcal{C}_m (X_m) \prod_{P \in \mathcal{P}} 
	\left( \sum_{\;\; i=1, \; X_i \in P }^s
	(Av(X_i) - (1/n)  \chi_i)   (Av(X)^{*} - (1/n) \bar{\chi}_i ) \right) v_P   
\end{align}
where $ v_P \in Y, \; P \in \mathcal{P}  $ is a set of free variables indexed by partitions in $ \mathcal{P} $  
\newline\newline\noindent
By analogy with Theorem 2.1 (section 2.3.1) we have 
\begin{proposition}\
\begin{itemize}
	\item[(i)] Irrep $ \rho$ satisfies the identity $\mathfrak{P}(\rho, P)  $ (4.5) 
	\item[(ii)] If  (4.5) holds in $ \sigma $  then 
	\begin{itemize}
		\item[(ii.1)] $ |H| \leq |G| \equiv m  $  
		\item[(ii.2)] if $ |H| = |G| $ then there is a partition $ P_H = (H_1, \; \cdots,\; H_s  )  $  of the group $ H $  (of type (4.4))  such that 
		$$ \sum_{h \in H_i } \sigma(h) = \left(\frac{t_i}{n} \chi_i \right)  I_W, \;  t_i = |H_i|, \;  i = 1, \; \cdots , \; s  $$ 
	\end{itemize}
\end{itemize}	
 
\end{proposition}
\noindent Proof. The statement (i) is obvious - an interior sum of (4.5) will vanish for a partition of $ \rho(G) $ that corresponds to (4.3). Further, if $ |H| < |G| $ then the guard term $ 	\mathcal{C}_m (X_m) $ vanishes identically on $ \rho(H) $. Hence, assume that $ H \geq m $. If the identity $\mathcal{P}(\rho)$ (4.5) holds in irrep $ \sigma $  then there is a partition $ P' = (H'_1, \; \cdots, H'_s   ), \;  |H'_i| = t_i, \; i = 1, \; \cdots, \; s $  of any $m$-element subset $ H' $ of $ H $ that satisfies the condition 
\begin{align}
	\sum_{i = 1, \; H'_i \in P' }^s \left( Av( \sigma(H'_i) ) -  (\chi_i/n ) I_W \right)  \left( Av( \sigma(H'_i) ) -  (\chi_i/n ) I_W \right)^{*} = 0
\end{align}  
\noindent The right        
hand side of (4.6) is a sum of positive semidefinite matrices and therefore 
\begin{align}
 \sum_{h' \in H'_i } \sigma(h') = \left(\frac{t_i}{n} \chi_i \right)  I_W, \;  i = 1, \; \cdots , \; s  \tag{4.6.1}
\end{align} 	 
If $ |H| > m $ then  there is  another $m$-element subset $ H'' $ of $ H $ that is obtained by replacing some $ h' \in H' $ by $ h'' \in  H \setminus H' $. Arguing as in the proof of Lemma 2.2 we find that $ \sigma(h') - \sigma(h'') $ is a scalar matrix and moreover, $H''$-version of (4.6.1)  implies that  $ \dim \sigma =1 $,  $H$ is cyclic and  $ \textnormal{range}(\chi_{\sigma} )  \subset 
\mathbb{Q}(\sqrt[m]1) $. This  contradiction with $ |H| > m $ assumption  finishes the proof of statements (ii).1 and (ii).2 
\begin{definition}
	Let $ \theta : \Gamma \rightarrow \textnormal{GL}(V) $ be an exact irrep of a finite group $\Gamma$. We will say that a partition $ P_{\Gamma} = (\Gamma_1, \; \cdots, \; \Gamma_l) $ of the group $\Gamma$ is \textit{central} $\theta$-partition if $ \sum_{g \in \Gamma_i}
	\theta(g) = \lambda_i I_V $ is a scalar matrix for all $ i = 1, \; \cdots, l $.
	We will call the (unordered) list of pairs $ (|\Gamma_i|, \; \lambda_i) $ the signature 
	of  central $\theta$-partition $ P_{\Gamma} $. Let's also say that central partition is nontrivial if it contains subsets that are not unions of full conjugate classes.  
\end{definition}
\noindent It is clear from the above that the recipe (4.5) yields an identity  
$\mathfrak{P}(\theta, P_{\Gamma})$  of $ \theta $ and as the proof of Proposition 4.2 suggests we have 
\begin{corollary}\
\begin{itemize}
	\item[(i)] Irrep $ \rho$ satisfies the identity $\mathfrak{P}(\rho, P_G)  $ for any central $\rho$-partition of $ G $ 
	\item[(ii)] If  the identity $\mathfrak{P}(\rho, P_G)  $ holds in $ \sigma $  then 
	\begin{itemize}
		\item[(ii.1)] $ |H| \leq |G| \equiv m  $  
		\item[(ii.2)] if $ |H| = |G| $ then there is a central $\sigma$-partition   
		of $ H $ that has the same signature as $ P_G $ 
	\end{itemize}
\end{itemize}
\end{corollary}
\noindent It is clear that any subset $X$ of $\Gamma$ (cf. definition 4.2) such that
$ \sum_{x\in X} \theta(x) $ is a scalar matrix gives rise to a central $\theta$-partition $ (X, \; \Gamma\setminus X ) $. The following example shows that nontrivial central partitions do exist.   
\begin{example}
	Let $ \Gamma = QC_3$  (a binary tetrahedral group) be a   semidirect product of a quaternion group $ Q $ and a cyclic group $C_3 = \; <h> $ of order $3$.
	This group has a natural embedding $ \theta: \Gamma \hookrightarrow \textnormal{SU}(2) $.
	It is obvious that $ \theta(h) + \theta(h^{-1}) $ is a scalar matrix and therefore we have a central $\theta$-partition $ (\{h, h^{-1}\}, \; \Gamma \setminus  \{h, h^{-1}\}) $. This central partition is nontrivial. Indeed, let's check  that  $ \{h, h^{-1}\} $ is not a conjugate class in $  \Gamma $. To evaluate  $  ghg^{-1}, \; g \in \Gamma $ we can assume that $ g \in Q$ and therefore  $ ghg^{-1} = (ghg^{-1}h^{-1})h  $ where the expression in brackets belongs to $ Q $. It is thus obvious that $ ghg^{-1} \neq h^{-1} $ for any $ g \in \Gamma $. It is no less obvious that there is $ g \in Q $ such that   $ghg^{-1}h^{-1} \neq  1 $ and hence $ ghg^{-1} \neq h $

\end{example}	
\begin{question}
	Is there a simple reason for existence of nontrivial central partitions? 
	Are these partitions "rare"?
\end{question}

\subsection{Relation probability}
Commuting elements $x,y$ of a group $\Gamma$ give rise to a relation $ [x,y] - 1 = 0 $ in any of its representations. 
It seems reasonable, therefore, to extend the notion of \textit{commuting probability} (cf. e.g.   \cite{Erdos}, \cite{Howe1}, \cite{commprob} and  references therein) in order to include more general relations.
   
If $ \theta : \Gamma \rightarrow \textnormal{GL}(V) $ is a representation of a group $ \Gamma $ then any $u= u(x_1, \cdots, x_p ) \in kF(Y)$ can be viewed as a function 
$$  u(\theta) : \Gamma^p \rightarrow \textnormal{End}(V), \; u(\theta)(g_1, \cdots, g_p ) =  u( \theta(g_1), \cdots, \theta(g_p)), \; g_1, \cdots, g_p \in \Gamma $$    
Let the group $\Gamma$ be finite. Define the probability  of the relation $ u $ over representation $\theta$  as  
\begin{align}  \textnormal{Pr}(u, \theta) = \frac{| \{ g \in \Gamma^p \; | \; u(\theta)(g) =0  \}|} { |\Gamma|^p} 
\end{align}	
\noindent 	     
Obviously $ \textnormal{Pr}(u, \theta)  = 1 $ if and only if $ u $ is an identity of $ \theta$. 
\newline\newline
Set $ \mathfrak{R}(u,\theta) = \{ g \in \Gamma^p \; | \; u(\theta)(g) =0  \} $.  If the representation $ \theta $ is unitary, then
\begin{enumerate}
	\item[(a)]   $ \mathfrak{R}(u,\theta) =  \mathfrak{R}(uu^{*},\theta) $ 
		\item[(b)]   $ \mathfrak{R}(u,\theta) \cap  \mathfrak{R}(v,\theta) =  \mathfrak{R}(uu^{*} + vv^{*},\theta) $ 
\end{enumerate}
\noindent Hence, one can define conditional probability (cf. \cite{commprob})
 $$\textnormal{Pr}(u|v, \theta) = \frac { | \mathfrak{R}(uu^{*} + vv^{*},\theta)|} {| \mathfrak{R}(vv^{*},\theta)}   
$$
\noindent where $ v $ is assumed to be non-identity of $\theta $.

\begin{proposition}
	Let  $ \theta : \Gamma \rightarrow \textnormal{U}(V) $ be an exact unitary irrep of a finite group $\Gamma$. 
	Then for any $ u  \in kF(Y) $ and any $ 0 <  r \leq 1 $ there is  $ u_r \in kF(Y) $  that is an identity of  $ \theta $ if and only if    
	$ \textnormal{Pr}(u, \theta) \geq  r $  
\end{proposition}
\noindent Proof.  Suppose that $u= u (x_, \cdots,x_p ) \in kF, \; x_1, \cdots, x_p \in Y   $ depends on $p$ free variables. Set $ m = | \Gamma| $, take $ p $ pairwise disjoint subsets $ X_1, \cdots , X_p \subset Y $ of size $ m $ each and set $ X = X_1 \times \cdots \times X_p $. Let 
\begin{align}
	u'_t = \prod_{i=1}^p \mathcal{C}_m(X_i) \prod_{S\in X, \; |S| = t } \; \sum_{ y \in S } u(y)u(y)^{*}  
\end{align}
\noindent where  $y$ runs over all tuples $ (y_1, \cdots, y_t ) \in S $. It is easy to see that Proposition 4.3 is equivalent to the following 
\begin{lemma}
	The identity $ u'_t $ (4.7) holds in $ \theta$ if and only if  	$ \textnormal{Pr}(u, \theta) \geq  t/|\Gamma|^p $ 
\end{lemma}   
\noindent This last statement, however,  is easily verified by noticing that (4.8) holds in $ \theta$ if and only if the number of tuples $ g =(g_1, \cdots, g_p) \in \Gamma^p $ such that $u(\theta)(g)= 0 $ is no less than $ t =|S| $.   

\begin{remark}
Essentially the definition (4.7) interprets $ u(\theta) $ as a random variable on $\Gamma^p$. Note that we have few times encountered expectations of random variables of this kind (cf.  remark 2.5, (4.3') and corollaries 2.4, 2.5 and 6.4). An application of expectation estimators for random variables related to commuting probability can be found in \cite{comm}
\end{remark}

\section{Theorem 1.1 and Central Polynomials}
We will present a new proof of Thorem 1.1 that works for nonmodular (ordinary) exact irreps over algebraically closed fields.  Sticking to original convention (cf. section 1.1)  we will be dealing with  arbitrarily fixed  faithful irreps 
$ \rho $ and $ \sigma$ over  $ k = \mathbb{C}$. We will see in section 5.3  that general case of ordinary representations is no different. 
The key to the proof we have in mind is  the  central polynomial of Razmyslov (5.2). Some steps in the proof seem to be interesting on their own and will be used to extend the result to a case of linear algebraic groups in characteristic $0$ (see Section 6). We begin, however, with a construction of a \textit{central Laurent polynomial} that can be associated with  an irrep of a finite group.

 Let 
 $
 \theta : \Gamma \rightarrow \textnormal{GL}(V)  
$ 
 be an exact irrep of a finite group $ \Gamma $. Following  well established terminology (cf. e.g \cite{Artin}) let us say that  $ c = c(y_1, \cdots, y_t )   \in kF $ is a central polynomial (central Laurent polynomial) of the irrep $ \theta$ if  $c$ is not an identity of $\theta$ but  $ c( \theta(g_1), \cdots, \theta(g_t) ) $ is a scalar matrix for any $ g_1, \cdots, g_t \in \Gamma $. We have already encountered (cf. Example 4.3) the Lurant polynomial $ y + y ^ {-1} $  that is central on any subgroup of $ \textnormal{SL}_2 (k) $.  
 Over algebraically closed field the linear span of $\theta(\Gamma)$ coincides with the full matrix algebra
 $ \textnormal{End}(V) \approx M_{\dim V}(k) $ and therefore, 
   any  multilinear central polynomial (cf. e.g. \cite{Razmyslov}) of  the matrix algebra $ M_{\dim V}(k) $ is also a central polynomial of $ \theta $. It is not hard, however, to come up with a  central Laurent polynomial of an  irrep of a finite group of a given order.
 \begin{lemma}
 	  Set  (cf.(1.3))
 \begin{align} 
 	  c_m \equiv	\Psi_m(\mathcal{C}_m(Y_m), Y_m )  
 \equiv  \sum_{k=1}^m  y_k \left( u_0\prod_{1 \leq i<j \leq m}^m  (y_i - y_j) u_{ij} \right) y_{k}^{-1} 
 \end{align} 
 where $ u_0, u_{ij} \in Y $ are pairwise distinct free variables in $ Y \setminus Y_m $. Then $ c_m$ is a central Laurent polynomial of $ \rho $  
 \end{lemma}  
 \noindent Proof. Let $  y_i \rightarrow g_i \in \rho(G), \; i = 1,  \cdots m $ 
 be an arbitrary variable assignment. If $ g_i =g_j $ for some $ i \neq j $ then $c_m $ vanishes in $ \rho $. Otherwise, it follows from  Lemma 2.1 that the value of $c_m$ is a scalar matrix 
 \begin{align}
 c_m(g,h) = (m/n)\textnormal{tr} \left( h_0\prod_{1 \leq i<j \leq m}^m  (g_i - g_j) h_{ij} \right)I_V  \tag{5.1'}
 \end{align} 
 where $  h_0,    h_{ij} \in  \rho(G)  $ is an arbitrary assignment of values for variables $ u_0, u_{ij} $. To finish the proof we need to show that $ h_0, h_{ij}$ can be chosen in such a way that expression under trace is nonzero. But this is exactly the claim of Lemma 1.1 (3).

\begin{example}
Let  $ \textnormal{Val}(p, \theta) = \{ \lambda_1 I_V, \cdots , \lambda_s I_V \} $ denote a set of all values $$ p(\theta(g_1), \cdots, \theta(g_t)) , \; g_1, \cdots , g_t \in \Gamma $$ 
of a central polynomial $p$ of the irrep $\theta$.  
 
		Let $ \gamma =  \{  \gamma_i \in k, \; i = 1, \cdots , r  \}   $ be a set of scalars. 
			  Set
		\begin{align} 
			q' = (q - \gamma_1)u_1 \cdots u_{t-1}(q - \gamma_r) \nonumber
		\end{align}
		  where $ q = q(y, \cdots, q_t)  \in kF$   
 and $ u_i \in Y , \; i = 1, \cdots , t-1  $ are additional free variables. 
Then  
	$ q' $  is an identity of the irrep $\theta$  if and only if $q$ is  central polynomial of $\theta$ such that $ \textnormal{Val}(q, \theta) \subset  \gamma I_V$. 
\end{example}
\begin{remark}
	The guard term (1.3) can be replaced with the central polynomial (5.1) in all the identities discussed so far.      
\end{remark}
\begin{example}
	If  a   representation $
	\theta : \Gamma \rightarrow \textnormal{GL}(V)  
	$  is not exact then $ c_{_{|\Gamma|}} $ is an identity of $\theta$. In particular, if $ \Gamma $ does not have exact irreps then  $ c_{_{|\Gamma|}} $ is an identity of the regular representation of $ \Gamma$ 
\end{example}
\subsection{Razmyslov’s central polynomial (\cite{Razmyslov})}
 This 
is a multilinear polynomial 
\begin{align}
 R_n  \equiv R_n(x_1, \cdots x_{n^2} ; \; y_1, \cdots, y_{n^2 + 1}  )  
\end{align}  in $ t = 2n^2 + 1 $ free variables. 
It  has the following remarkable properties :
\begin{enumerate}
	\item[1)] $R_n$ is skew symmetric in $x_1, \cdots, x_{n^2}$
	\item[2)] there are $ a_1, \cdots a_{t} \in M_n(k) $ such that $ R_n(a_1, \cdots a_{t} ) \neq 0 $
	\item[3)]   $ R_n(a_1, \cdots a_{t} )   $ is a scalar matrix for any  $ a_1, \cdots a_{t} \in M_n(k) $
	
\end{enumerate}	
As a direct consequence of 1) we note that the matrix algebra $ M_n(k) $ satisfies  polynomial identity in $(2n^2 + 2)$ variables 
\begin{align} 
	R'_{n+1} = R_{n}( x_1, \cdots x_{n^2}; \; \cdots )x_0 + 	\sum_{i=1}^{n^2}  (-1)^i R_{n}(  x_1, \cdots x_{n^2}; \; \cdots )|_{x_i \leftarrow x_{0}} x_i 
\end{align}
where $ x_0 $ is an additional free variable 
\subsection{Theorem 1.1 over the field of complex numbers}

 $N$-variable identities of $ \rho $ form an ideal $ U_N = U_{\rho, N } \subset kF_N $ where  $  F_N = F(Y_N) $ is a free group generated by $N$ free variables. Set $ A_N = A_{\rho,N} =  kF_N/U_N $ and
let $ C_N $ be the center of $ A_N$.   
Assuming that $ N \geq  t + s   $  set $$ C_N \ni R_0 = R_n(x_1, \cdots, x_t ) \!\! \!\!\! \mod \! U_N , \; C_N \ni c_0 = c_m(x'_1, \cdots, x'_s) \!\! \!\!\! \mod \! U_N$$ where $ t =2n^2 + 2, \; s = m(m+1)/2 $ and all variables $ x_i, \; x'_j \in Y, \; i = 1,\cdots, t, \; j = 1, \cdots, s $ are pairwise distinct. 

Let $ \rho' : G' \rightarrow V' $ be a representation that satisfies all the identities of $ \rho$.  Any value assignment 
 $ F_N \ni y \rightarrow g \in G', \; y \in Y_N  $ uniquely extends to a group  
 homomorphism $ \phi : F_N \rightarrow G' \approx \rho'(G') \subset \textnormal{End}(V)  $. 
 Let's denote the (finite) set of all such homomorphisms by  $  \textnormal{H}({\rho'})$. Any $ \phi \in \textnormal{H}({\rho'})$  uniquely extends to a  homomorphism $ \phi : kF_N \rightarrow \textnormal{End}(V')  $ which in turn factors through into
 (denoted by the same letter) homomorphism $ \phi : A_N  \rightarrow  \textnormal{End}(V') $.
 For  $ \phi \in \textnormal{H}({\rho})$  set $  \mathfrak{p}_{\phi} =  \ker \phi \cap C_N = \ker (\phi | C_N ) $.

Any identity of $ \rho $ obviously holds in all factors of $\rho$ (cf. Definition 4.1) and we get the following 
 
\begin{lemma}\
	\begin{enumerate}
	\item[(i)] For any  $ \phi \in \textnormal{H}(\rho) $ there are  irreducible factors $ \rho_i : G_i \rightarrow \textnormal{End}(V_i), \; i = 1, \cdots, l $ of $ \rho $ and epimorphisms
	$ \textnormal{H}(\rho_i) \ni \phi_i : F_N \rightarrow G_i \approx \rho_i(G_i)  $ such that 	
\begin{align}
\ker (\phi | C_N ) = \bigcap_{i=1}^l   \ker (\phi_i | C_N ) \nonumber
\end{align} 
All ideals $ \ker (\phi_i | C_N ), \; i = 1, \cdots,l $ are maximal 	
 
	\item[(ii)]  If   $ \phi \in \textnormal{H}(\rho) $ is an epimorphism (or if $ \phi(R_0) \neq 0 $) then $ \mathfrak{p}_{\phi}  $ is a maximal ideal in $ C_N$
	\end{enumerate}
\end{lemma}
 
\noindent Proof.  Under conditions of statement (ii) 
$\phi | C_N $ is actually an epimorphism onto the ground field.   
To prove (i), suppose that $ \phi(F_N) $ is a (proper) subgroup $ \rho(G') \subset \rho(G)$ of $ G $. 
The representation  $ \rho|G' : \rho(G')  \hookrightarrow \textnormal{GL}(V)  $ splits into a direct sum of irreducible components
$ \rho_i : G' \rightarrow \textnormal{GL}(V_i), \; V_i \subset V , \; i = 1, \cdots, l $.  Therefore $ \phi(C_N )  = \oplus_{i=1}^l \pi_i \phi(C_N)$  where $ \pi_i $ is a linear  projector onto $ V_i$. Clearly, all $ \phi_i = \pi_i \phi \in \textnormal{H}(\rho_i) $ are surjective and the result now follows from (ii)

\begin{lemma} \

	\begin{enumerate}
		\item[(i)] $ C_N \ni R_0 c_0 \; (\neq 0) $ is a central (Laurent) polynomial of $ \rho $
		\item[(ii)] If   $ \phi \in \textnormal{H} (\rho)   $ is not surjective then $ \phi(c_0) = 0 $
		\item[(iii)]  If   $ \phi \in \textnormal{H} (\rho)   $  is a homomorphism such that $ \phi(R_0) \neq 0 $ then
		$ \ker \phi = \mathfrak{p}_{\phi} A_N  $ 
	\end{enumerate}
\end{lemma}    
\noindent Proof. Statements (i) and (ii) follow from  definitions.To prove (iii), take $ a \in \ker \phi$ and  
note that by (5.3) the following equation holds in $ A_N $   
\begin{align} 
	 R_{n}( x_1, \cdots x_{n^2}; \cdots )a \; = \; - 	\sum_{i=1}^{n^2}  (-1)^i R_{n}(  x_1, \cdots x_{n^2}; \cdots  )|_{x_i \leftarrow a} x_i 
\end{align}
\noindent where only the first $ n^2$ variables of $ R_0 $
are shown. It is easy to see that  $  R_{n}(  x_1, \cdots x_{n^2}; \cdots  )|_{x_i \leftarrow a}  $ belong to  $ \mathfrak{p}_{\phi} $ for all $i = 1, \cdots, n^2$. Further,
$ \mathfrak{p}_{\phi}$ is a maximal ideal of $ C_N$  and therefore $  R_0 = R_{n}( x_1, \cdots x_{n^2}; \cdots )$ is invertible  in $ C_N $  modulo $  \mathfrak{p}_{\phi} $ i.e. there is $ r \in C_N $ such that $ r R_0 = 1 + b, \; b \in  \mathfrak{p}_{\phi} $ and we have 
$  r R_0 a  = (1 + b ) a \in   \mathfrak{p}_{\phi} A_N \implies a \in    \mathfrak{p}_{\phi} A_N $.
\newline\newline\noindent
We will need below a standard fact of commutative algebra (cf. e.g. \cite{Miln-CA}). Let $ A $ be a commutative ring. Suppose that the  spectrum   of $ A $  is a finite union of irreducible components 
\begin{align}
	\textnormal{spec}(A) = V_1 \cup \cdots \cup V_l 
\end{align}
\noindent then the irreducible decomposition (5.5) (when it exists) is said to be irredundant if there are no inclusions between components $ V_i $. Clearly, for corresponding prime ideals $ P_i \subset A $
one has $ P_1 \cap \cdots \cap P_l = \{0\} $ with no inclusions between $ P_i $.   A well known simple fact that we need is   
\begin{lemma}(cf. e.g. \cite{Miln-CA}). An irredundant finite decomposition of a spectrum of a commutative ring is unique up to a permutation of components.
\end{lemma}

\noindent Consider now the spectrum $ X  = \textnormal{spec}(C_N) $ of the ring $ C_N$. 

\begin{lemma} (cf. \cite{Miln-CA})
	\begin{enumerate}
		\item [(i)] 	There are homomorphisms  $ \phi_i : F_N \rightarrow G,  \; i = 1, \cdots, l  $ such that $ \mathfrak{p}_{\phi_1} \cap \cdots \cap \mathfrak{p}_{\phi_l} = \{0\}  $
		\item [(ii)] 	$\textnormal{spec}(C_N) $ is finite, i.e. there is an irredundant set $	\mathfrak{m}_1 , \cdots , \mathfrak{m}_q $ of  maximal ideals of $ C_N$ such that
		\begin{align}
			\mathfrak{m}_1 \cap \cdots \cap \mathfrak{m}_q = \{0\}
		\end{align} 
		\item [(iii)]  The list (5.6) contains an ideal 
		 $ \mathfrak{m} = \mathfrak{p}_{\phi} = \ker \phi \bigcap C_N = \ker (\phi | C_N ) $ such that $ \phi(R_0 c_0) \neq 0 $

	\end{enumerate}	
\end{lemma}
\noindent Proof. First of all, the set  $ \textnormal{H}(\rho) $ is finite.  If  $ u \in  \cap_{\phi \in \textnormal{H}(\rho)} \mathfrak{p}_{\phi} $  then $u$  is an identity of $ \rho$, or in other words   $ u = 0 \!\!\mod U_N$. That proves (i). The existence of unique (up to reordering)  decomposition of $ C_N$ with corresponding list of maximal ideal (5.6) follows from (i) and  Lemma 5.4.  By (i) (and Lemma 5.3 (i)), the list (5.6) contains an ideal  $  \mathfrak{m}$ such that  $ R_0 c_0 \notin\mathfrak{ m} $ and  by Lemma 5.3 (i) and Lemma 5.5 (i) there is $ \phi 
\in \textnormal{H}(\rho) $ such that   $ \mathfrak{p}_{\phi} \subset \mathfrak{m} $. By Lemma 5.4 (ii) $  \mathfrak{p}_{\phi} $ is maximal  and therefore the ideal $ \mathfrak{p}_{\phi} = \mathfrak{m} $ satisfies the claim (iii).   
\newline\newline

\subsubsection{Conclusion of the proof} It is easy to see (cf. e.g. sections 1-3) that if exact irreps $ \rho $ and $\sigma$
have the same identities then $ \dim \rho = \dim \sigma = n $
and $ |G| = |H| = m $. It follows  then from Lemma 5.5 (iii) and Lemma 5.3 that  we have a pair of epimorphisms 
$ \phi : A_N  \rightarrow  G  $ and $ \psi : A_N   \rightarrow H $
such that 
\begin{enumerate}
	\item [(a)] $ \phi(F_N) = G $ and   $ \psi(F_N) = H $ 
	\item [(b)] $ \ker \phi \cap C_N   =\ker \psi \cap C_N $ and 
	$ \phi(R_0) \neq 0 \neq \psi(R_0)$
	\item [(c)] and therefore, by Lemma 5.3 (iii) $ \ker \phi = \ker \psi \; (\subset A_N) $   
\end{enumerate}
Hence there is an automorphism $ \alpha $ of the  matrix algebra  $ \textnormal{M}_n(k) $ such that 
$ \psi = \alpha \phi $. Since all automorphisms of $ \textnormal{M}_n(k) $  are inner,  $ \rho(G) $ is conjugate to $ \sigma(H)$ and that completes the proof.        

\subsection{Theorem 1.1 over algebraically closed fields}
 Let's try to figure out if/what results of section 5.2 remain valid in characteristic $ p$. As before, we will be dealing with  fixed exact irreps $ \rho, \; \sigma $ (cf. Introduction)  but this time  over an algebraically closed field $k$ of arbitrary characteristic $  \textnormal{char}(k) = p$ or zero. We will retain  notation and setup of Section 5.2.

 Lemmas 5.2 and 5.1 (note that
 $  ( p , m ) = 1 \implies (m/n, p) = 1 $) remain valid when $ (p,m) = 1 $  and it is easy to see, therefore, that all statements of the section 5.2 are valid in ordinary case. 
 
   Although in modular case, Lemma 5.2 (i) is no longer valid, we can proceed as follows. Let $ \phi' \in \textnormal{H}(\rho') $ be an epimorphism where $ \rho' : G' \rightarrow \textnormal{GL}(V') $  is an irreducible factor of $ \rho $. Clearly, $C_N$-ideal $ \ker \phi' \cap C_N$ is maximal. Let $ \mathfrak{m}_0$ denote the intersection of all such ideals. 
\begin{lemma}
	The ideal $ \mathfrak{m}_0$ is a nil-radical of $C_N$. Actually, $ u \in \mathfrak{m}_0 $ if and only if $ u^n $ is an identity of $ \rho$.  
\end{lemma}    
\noindent Proof. Take   $ u \in \mathfrak{m}_0 $ and let $ \phi \in  \textnormal{H}(\rho) $. Let $  \phi(F_N)=\rho(G')  $ for some 
subgroup $ G' \subset G $ and let  $ \{0\} \subset V_1 \subset \cdots \subset V_l = V $ be  invariant composition series of $ \rho|G' $. Obviously $ l \leq n $ and  $ \phi(u) $ acts trivially in all  irreducible factors $ W_i = V_{i}/V_{i-1} $. Hence $ \phi(u)^n = 0 $ and since $ \phi \in  \textnormal{H}(\rho) $ is arbitrary, $ u^n$ must be an identity of $ \rho$, i. e.  $ u^n = 0 $.    The reverse statement is obvious. 
\newline\newline\noindent 
Now  we have  a straightforward generalization of  Lemma 5.5 that does not exclude the case of modular representations    
\begin{lemma5}(cf. \cite{Miln-CA})
	\begin{enumerate}
		\item [(i)] 	There are homomorphisms  $ \phi_i : F_N \rightarrow G,  \; i = 1, \cdots, l  $ such that $ \mathfrak{p}_{\phi_1} \cap \cdots \cap \mathfrak{p}_{\phi_l} = \mathfrak{m}_0  $
		\item [(ii)] 	$\textnormal{spec}(C_N) $ is finite, i.e. there is an irredundant set $	\mathfrak{m}_1 , \cdots , \mathfrak{m}_q $ of  maximal ideals of $ C_N$ such that
		\begin{align}
			\mathfrak{m}_1 \cap \cdots \cap \mathfrak{m}_q = \mathfrak{m}_0 \tag{5.6'}
		\end{align} 
		\item [(iii)]   The list of maximal ideals  (5.6') contains an ideal 
		$ \mathfrak{m} = \mathfrak{p}_{\phi} = \ker \phi \bigcap C_N = \ker (\phi | C_N ) $ such that $ \phi(R_0 c_0) \neq 0 $
			
	\end{enumerate}	

\end{lemma5}
\noindent The remaining difficulty is that Lemma 5.1  does not hold for modular representations. 
For example,  if the field characteristic $ p $ divides $m=|G|$ then  the right hand side  of (5.1)  evaluates in $\rho$ to a traceless scalar matrix $ \lambda I_V$  that must be zero if  $(n=\dim \rho, \; p ) = 1$. However, conditions for $c_m$ (5.1) to be a central Laurent polynomial are more or less known within the modular representation theory. To specify  these  conditions we need to review some relevant terminology  (cf. e.g. \cite{Serr}, \cite{Webb})

\subsubsection{Reduction mod-$p$ (cf. e.g. \cite{Serr}, \cite{Webb})}
\noindent  It is easy to see that results discussed so far  for complex representations remain valid over  possibly smaller splitting field(s) (cf. e.g. Lemma 2.4 (iii)). To apply $\mod p$ reduction to  a complex representation $\theta$  we can switch to a splitting field $ K' \subset \mathbb{C} $ of $ \theta$ (that is to an appropriate finite algebraic extension of $\mathbb{Q}$). The field $ K $ can be equipped with a discrete valuation ring  $ \mathfrak{A} \subset K' $  with maximal ideal $\mathfrak{o} \subset A $ such that $ k = \mathfrak{A}/\mathfrak{o} $ (cf. \cite{Miln-ANT}). A standard Cauchy sequence procedure can be then applied to complete the field with respect to the metric induced by the valuation (cf. \cite{Miln-ANT}). Essentially, to quote from \cite{Serr},
"We denote by K a field complete with respect to a discrete valuation ... with valuation ring $\mathfrak{A} $, maximal ideal $\mathfrak{o}$ and residue field
$ k = \mathfrak{A}/\mathfrak{o}$. We assume that $K$ has characteristic zero and that $k$ has
characteristic $p$". It should be clear that we can always assume that $ K $ (and/or $k$) is a splitting field for any representation over $ K $ (and/or $k$) that comes around.    

\begin{definition} (cf. e.g. \cite{Serr}, \; \cite{Webb}).
	 $k$-irrep $\theta $ of a finite group $ \Gamma$  
	is called liftable if there is a $K$-irrep $ \theta' $ of $ \Gamma $ in a vector space $V'$ and an $\mathfrak{A}\Gamma$-invariant submodule $ L \subset V' $ such that $ k\Gamma $-module $ L/\mathfrak{o}L $ is equivalent to $ \theta $. Under these conditions $ \theta'$ is called a \textbf{lift} of $ \theta$ and $ \theta $ is called  a \textbf{reduction} of $ \theta'$. 
\end{definition}

\begin{lemma}
	For the  irrep $ \rho : G \rightarrow \textnormal{End}(V) $ the following conditions are equivalent
	\begin{enumerate}
		\item[(i)] The Laurent polynomial (5.1)
		\begin{align} 
			c_m \equiv	\Psi_m(\mathcal{C}_m(Y_m), Y_m )  
			\equiv  \sum_{k=1}^m  y_k \left( u_0\prod_{1 \leq i<j \leq m}^m  (y_i - y_j) u_{ij} \right) y_{k}^{-1}   \nonumber 
		\end{align}  
		is a central  polynomial of $ \rho$ 
		\item[(ii)] $ V $is a projective $kG$-module
		\item[(iii)] $ \chi_{ \rho } $ is a $p$-defect zero character, i. e. $ (m = | G| )/(n = \dim V ) $ is relatively prime to $ p$    
	\item[(iv)]	The irrep $ \rho $ is liftable
	\end{enumerate} 
\end{lemma}
\noindent For the equivalence of statements (ii)-(iv) we refer the reader to \cite{Serr}, \cite{Webb}.  We will outline two simple proofs of the equivalence of statements (i) and (ii)-(iv). The first observation is that equivalence of (i) and (ii) follows from Lemma 1.1 ((1)-(3)) and a well known 
\begin{lemma}(cf. e.g. \cite{Serr}).
	A $ kG$-module $V$ is projective if and only if 
	\begin{align}
		A_{\rho}(u) \; \equiv \;  \sum_{g \in G }  \rho(g)u \rho(g ^{-1})    
	\end{align} 
	\noindent is an identity map for some $ u \in \textnormal{End}(V) $ 
\end{lemma}
\noindent  Another observation is that (i) follows from (iv) as we have the following  
\begin{lemma}
	If the (exact) irrep $ \rho $ (over $k$) is  of defect zero then $ c_m $ (5.1) is a central polynomial of $ \rho$  
\end{lemma}
\noindent Proof. As in (the proof of) Lemma 5.1 (cf. also Lemma 1.1) there is an assignment of variables $ y_i \rightarrow g_i, h_0\rightarrow u_0, h_{ij} \rightarrow h_{ij}; \; g_i,h_0, h_{ij} \in \rho(G) $ such that  the expression 
\begin{align}
	u_0\prod_{1 \leq i<j \leq m}^m  (y_i - y_j) u_{ij}
\end{align} 
 evaluates (in $\rho$) to a matrix $ a \in  \textnormal{End(V)}$ such that $ tr(a) \neq 0 $.  
Let $ \tilde{\rho} $ be a $K$-lift of $\rho$. Then we have $ \tilde{a} \in L $ such that $ tr(a) = tr(\tilde{a}) \!\!\!\mod \!\mathfrak{o} $. By Lemma 5.1 (cf. (5.1')), on corresponding lifted variables the Laurent polynomial $ c_m $ (5.1) evaluates 
to $ (m/n)tr(\tilde{a}) $ in $ \tilde{\rho} $ and since $  m/n \neq 0 \!\!\! \mod \! \mathfrak{o}$ the reduction  $ (m/n)tr(\tilde{a}) \!\!\! \mod \!\mathfrak{o} \; = \; (m/n)tr(a)   $ is not  zero either.
\newline\newline
\noindent Hence, as a modification of Theorem 1.1  we have 

\begin{thover} Let the ground field $k$ be algebraically closed.
	Suppose that exact irreps $ \rho: G \rightarrow \textnormal{End}(V) $ and $ \sigma: H \rightarrow \textnormal{End}(W)  $ have the same identities.  If $V$ is a projective $kG$-module then $ \rho $ is similar to $\sigma$
\end{thover}

\begin{remark}
	It should be clear that in ordinary case any $G$-representation is a projective $kG$-module 
\end{remark}

\noindent A group is called $p$-solvable if it has a composition series with each factor either a $p$-group or a group of order prime to $p$. By the Fong-Swan Theorem (cf. e.g. \cite{Serr}) any irrep of a $p$-solvable group is liftable and we have 
\begin{corollary}
Let the ground field be algebraically closed.
Suppose that exact irreps $ \rho: G \rightarrow \textnormal{End}(V) $ and $ \sigma: H \rightarrow \textnormal{End}(W)  $ have the same identities.  If the group $ G $ is $p$-solvable then  $ \rho$ and  $ \sigma$  are similar 
\end{corollary}  

\begin{remark} (cf. \cite{Granville}).
 With a few exceptions, simple finite groups  have a  $p$-defect zero character for any prime $p$  
\end{remark}

\section{Linear Algebraic Groups}
All representations in this section (unless explicitly stated otherwise) are over algebraically closed field of characteristic zero. We will consider  linear  \textit{algebraic groups}, i.e. we assume that
 $ \rho(G) \subset \textnormal{GL}(V)$ and $ \sigma(H) \subset \textnormal{GL}(W) $ are closed  in Zariski topology.
\newline
\newline
\noindent  
As was mentioned above, faithful irreps of finite groups  with the same identities are similar. The same is true for a somewhat opposite case of  finite dimensional irreps of connected algebraic groups (see \cite{irrop}). In fact, faithful finite dimensional irrep of an arbitrary linear algebraic group is determined by  identities. More precisely, we have the following  
\begin{theorem}
Let $ \rho : G \rightarrow \textnormal{GL}(V) $ and $ \sigma : H \rightarrow \textnormal{GL}(W) $ be faithful finite dimensional irreps with the same identities. Suppose that $ \rho(G) \subset \textnormal{GL(V)} $ and $ \sigma(H) \subset \textnormal{GL(W)} $ are linear algebraic groups. Then representations $ \rho$ and $ \sigma$ are similar.
\end{theorem}
\noindent For a subgroup $ \Gamma \subset \textnormal{GL}(V) $ denote by $ \overline{\Gamma} $ its closure in Zariski toplogy (Zariski closure).
It is well known (and easy to see) that  representations $ \Gamma \hookrightarrow \textnormal{GL}(V) $ and  $ \overline{\Gamma} \hookrightarrow \textnormal{GL}(V) $ have the same identities. Hence we have the following 
\begin{corollary}
If $ \rho : G \rightarrow \textnormal{GL}(V) $ and $ \sigma : H \rightarrow \textnormal{GL}(W) $ are finite dimensional irreps with the same identities then  $ \dim V = \dim W $ and  groups $ \overline{\rho(G)}, \; \overline{\sigma(H)} $ are conjugate to each other in $ \textnormal{GL}(V) $ 
\end{corollary}

\begin{corollary}
	If a natural group representation  $ \Gamma \hookrightarrow  \textnormal{GL}_n(k)  $ is irreducible then   Zariski closure $ \overline{\Gamma} $ is a maximal subgroup of $ \textnormal{GL}_n(k) $ that contains $\Gamma$ and still satisfies all  identities of the irrep 
	$ \Gamma \hookrightarrow \textnormal{GL}(k) $. 
\end{corollary}

\noindent The proof of the Thorem 6.1 will be given in  few steps.
\paragraph{Step 1.} We will identify  $ \rho(G) $ with the set of zeroes in $ GL_n(k) \subset M_n(k) $ of a polynomial ideal $ J = J_G  \subset A =  k[T_{11}, \cdots, T_{nn}, \textnormal{det}(T_{ij})^{-1} ] $.
The ring 
\begin{align}
\mathcal{O} \equiv \mathcal{O}(G) =  k[t^{1}_{11},  \cdots , t^{1}_{nn}] 
\end{align}
\noindent of polynomial functions on $ \rho(G) $ (coordinate ring)     is generated by 
images $ t_{ij}$ of generic
matrix coordinate functions  $ T_{ij} $ in $ A/J, \;  i,j = 1, \cdots, n $. Before going forward we will review  some basic information on connected components of algebraic groups (cf.  \cite{Miln-ALG}).  An algebraic group $G$ contains a connected normal subgroup  $G^{0}$ such that factor group $ G/G^{0}$ is a finite algebraic group. 
Let $ \mathcal{O}_0 = \mathcal{O}(G^{0}) $ be a coordinate ring of $G^0$.  The connected components of $G$ are the fibers of the map
\begin{align}
\pi_0 : G \rightarrow   G/G^{0} \approx \pi_0(G)
\end{align}
and  $ \mathcal{O} = \mathcal{O}(G) $ is a direct sum of integral domains that are  isomorphic copies of  $ \mathcal{O}_0$. Each fiber of (6.2) is a
connected coset  $ \tilde{g}G_0 $ where  $ \tilde{g} $ is some  coset representative of  $ g \in G / G^{0} $. It follows from the existence of the   fibration (6.2) that   $ \textnormal{spec}( \mathcal{O}) $ is a dijoint union of open subsets, i.e. there are orthogonal idempotents  $$ e_g \in \mathcal{O}, \;    e_g e_h = \delta_{gh} e_g, \; g, h  \in  G/G^{0} $$  such that
\begin{equation}
\mathcal{O} = \bigoplus_{g \in G/G_0} \mathcal{O}_g \;\; \textnormal{where } \mathcal{O}_g \equiv \mathcal{O}e_g \approx  \mathcal{O}e_1 \approx \mathcal{O}_0
\end{equation}  
   In other words,  there is a maximal  (\'etale) subalgebra $ \pi_0(\mathcal{O}) \subset \mathcal{O}  $ generated by idempotents $ e_g, \; g \in  G/G^{0}  $. The subalgebra   $ \pi_0(\mathcal{O}) $  is a direct sum of $ |\pi_0(G)| \equiv |G/G_0|
$ copies of the ground field $ k $
 and  
$ \mathcal{O} = (\mathcal{O}_0 \equiv \mathcal{O}e_1) \otimes_k \pi_0(\mathcal{O} ) $.
 A fiber (6.1) over $ \pi_0(g), \; g \in G/G^{0} $
is a connected coset of  $G^{0}$ in $ G $ and  as  algebraic subset  of $  \textnormal{GL}(V) $ coincides with $ \textnormal{spec}( \mathcal{O}e_g)   $.
If $ \rho'(G') \subset \textnormal{GL}(V') $ 	is another linear algebraic group then   
\begin{align}
\mathcal{O}(G \times G') =   \mathcal{O}(G) \otimes  \mathcal{O}(G'), \; \;
	 \mathcal{O}( (G \otimes (G')^{0} ) =  \mathcal{O}(G^0 \times G^{'0}) =   \mathcal{O}(G^0) \otimes  \mathcal{O}(G'^0)   
\end{align}
and 
\begin{align}
	\!\!\!\!\!\!\!\!\! \pi_0(\mathcal{O}(G\times G') ) =   \pi_0(\mathcal{O}(G) ) \otimes \pi_0(  \mathcal{O}(G') ), \;\; \pi_0(K \otimes_k \mathcal{O} ) = K \otimes_k \pi_0(\mathcal{O})
\end{align}
for any field extension $ K \supset k $.

\paragraph{Step 2.} Consider   $ N$-fold direct product 
$$ G_N = \rho(G)  \times  \cdots \times \rho(G)  \supset \rho(G^{0}) \times \cdots \times \rho (G^{0}) = G^{0}_{N} $$ 
as a linear
algebraic  subgroup in $ N$-fold direct sum $ \oplus^N M_n(k) \supset \textnormal{GL}(V) \times \cdots \times \textnormal{GL}(V) $.
The coordinate ring of $ G_N $ is  $N$-fold tensor product of the coordinate ring $ \mathcal{O} $ of $ G $
\begin{align}
 \mathcal{O}_{N} \; = \;  \bigotimes^N\mathcal{O} \; \equiv  \; k[t^{1}_{11},  \cdots , t^{1}_{nn}; \; \cdots ; \; t^{N}_{11},  \cdots , t^{N}_{nn}]  	
\end{align}
where coordinate functions $ t^r_{ij}, \; r = 1, \cdots, N; \; i,j = 1, \cdots, n $ have the same meaning as in (6.1).
Take the matrix algebra 
$M_n( \mathcal{O}_N ) = \mathcal{O}_N \tens{k} M_n(k)  $ with coefficients in  $ \mathcal{O}_N $ and set 
\begin{align}
 \mathcal{M}_r = \sum_{i,j=1}^n t^r_{ij} e_{ij}, \; r = 1, \cdots, N  
\end{align}
\noindent where $ e_{ij} = e_i \otimes e^{*}_j$ are standard matrix units.
As in the Section 5.2, let    $ kF_N = kF(Y_N) $ be the group algebra of a free group with a set of free variables $ Y_N = \{y_1, \cdots , y_N \}$.
Define a homomorphism $ \phi : F_N \rightarrow M_n(\mathcal{O}_N) $ by setting
$ \phi( y_i ) = \mathcal{M}_i, \; i = 1, \cdots, N  $ and extend this homomorphism to 
$ kF_N$ by linearity. Let $U_N$ be an ideal of $N$-variable identities of the irrep $ \rho$ (and thereby of the irrep $\sigma$).  Set $ A_N = kF_N/U_N $ and with a slight abuse of notation assume that images of free variables $ y_i, \; i = 1, \cdots, N $ in $ A_N$ and   the image of $ F_N $ in $ A_N $ are denoted by the same letters. The following lemma is quite obvious (cf. \cite{irrop}).
\begin{lemma}   
The ideal 	$ U_N $  of $N$-variable identities of $ \rho$  coincides with the kernel of the homomorphism
\begin{align}
\phi : kF_N \rightarrow \textnormal{M}_n(\mathcal{O}_N), \;\;
 \phi( y_i ) = \mathcal{M}_i, \; i = 1, \cdots, N  \nonumber
\end{align}
and therefore $ \phi $ induces an injection (denoted by the same letter)
$$  \phi : A_N \equiv kF/U_N  \hookrightarrow \textnormal{M}_n(\mathcal{O}_N) $$ 
\end{lemma}

\noindent As in (6.2) we have a fibration
\begin{align}
	G_N \rightarrow \pi_0(G_N) \approx G_N/G_{N}^{0}
\end{align}
\noindent that for the reasons outlined  above (cf. (6.3)-(6.5)) induces a direct sum  splitting of the 
coordinate ring

\begin{align}
	\mathcal{O}_{N} = \underset{g \in G_N/G_{N}^{0}}\bigoplus  e_g	\mathcal{O}_{N} \;\textnormal{  where  }   e_g = e_{(g_1, \cdots, g_N) }, \;  
	g = (g_1, \cdots, g_N), \; g_i \in G/G^{0}, \; i =  1, \cdots , N  \nonumber
\end{align}
where  $g_i $ are fixed representatives of connected  cosets  $  g_iG^0 \subset G ,  \; i = 1, \cdots, N $ and 
$ e_g $ are corresponding orthogonal idempotents in  $\mathcal{O}_{N} $ (cf. step 1). In other words, by  convention (6.3) we write
\begin{align}
	e_g \mathcal{O}_N \; \equiv \; \mathcal{O}_{g} \; = \;  \mathcal{O}_{g_1} \otimes \cdots  \otimes \mathcal{O}_{g_N}
	\approx \overset{N}\bigotimes \mathcal{O}_0	
	, \; g \in G_N/G_{N}^{0} \nonumber
%
%
\end{align}
\noindent where 
  coordinate rings $ \mathcal{O}_g $  are pairwise isomorphic integral domains that correspond to  fibers of the map (6.8). 
\paragraph{Step 3.} Let's fix once and for all a universal (algebraically closed with infinite transcendental degree) extension $K$ of the ground field $ k $.  A field  of fractions of the integral domain $\mathcal{O}_g  $ embeds into $ K $
for any $ g \in G_N/G_{N}^{0} $
and there is an induced homomorphism of matrix rings
\begin{align}
\!\!\!\!\!\!\!\!\!\!\!\!\!\!\!\! 
\pi_{g} :  M_n( \mathcal{O}_N) \rightarrow  
	K \otimes_{\mathcal{O}_g}
	 M_n(\mathcal{O}_g ) 
	 \; \approx \; 
	 M_n(K) 
	 , \;\;  g \in G_N/G^{0}_N   
\end{align}
where, of course, the tensor product depends on embedding of $\mathcal{O}_g$ into $K$.  Composing   homomorphisms (6.9) with the homomorphism  of Lemma 6.1 we get $k$-algebra homomorphisms 
\begin{align}
\pi_g \phi  \equiv\phi_g: A_N \;\; \rightarrow   M_n(K ) , \;
\phi_g(y_r) =  \mathcal{M}_{r,g} \equiv  e_g \mathcal{M}_r, \; r = 1, \cdots, N, \;
 	g \in G_N/G^{0}_N 
\end{align}
\noindent where $ \mathcal{M}_{r,g} = e_g \mathcal{M}_r $ is the matrix with coefficients $  t^r_{ij} | gG_N^0 $ (cf. (6.7)-(6.9)). It is obvious (cf. Lemma 6.1 ) that 
\begin{align}
	\{0\} = \ker \phi \; = \bigcap_{	g \in G_N/G^{0}_N} \ker \phi_g 
\end{align}
\begin{lemma}
	If a set of distinct cosets  in  $	h = (h_1, \cdots, h_N ) \in G_N/G^{0}_N $ is contained in a set  of distinct cosets 
	in $	g = ( g_1, \cdots, g_N )  \in G_N/G^{0}_N $ then $ \ker \phi_g \subset \ker \phi_h $. 
\end{lemma}  
\noindent Proof. Using the definition of te matrices $M_{r,g} $, note  that a polynomial that vanishes on a set  $S$, certainly vanishes on any subset of $S$.
\newline\newline\noindent
From now on, fix an integer $ N $ such that
$ N \geq \max \{ 2n^2 + 2, |G/G_0| \} $  
\begin{lemma}\

	\begin{enumerate}
		\item[(i)] there is $ g  \in G_N/G^0_N $ such that  
		$ R_n(\mathcal{M}_{i_1,g}, \cdots , \mathcal{M}_{i_t,g}) \neq 0 $
			for some set of indices  $ \{ i_1, \cdots, i_{t=2n^2 +1} \} \subset \{ 1, \cdots, N  \} $; the condition 	$ R_n(\mathcal{M}_{i_1,g}, \cdots , \mathcal{M}_{i_t,g}) \neq 0 $ is equivalent to the condition $  K\phi_g(F_N)  =  M_n(K) $ 
		\item[(ii)] any $ 	g = (g_1, \cdots, g_N)  \in G_N/G^0_N $
		such that the set $ \{ g_1G^0, \cdots g_NG^0 \} $ contains all the factor sets of $ G/G^{0}$  satisfies the condition (i) 
%
		
		\item[(iii)] 
		 if $g \in G_N/G^0_N$  satisfies the condition (ii) then   the closure of the group 
		$ \phi_g(F_N) $ in Zariski topology of $\textnormal{GL}_n(K)$  coincides with the group of points  $ G(K) \equiv \textnormal{spec}( K \otimes_k \mathcal{O}) $ 
		
		\item[(iv)] for any  $ g  \in G_N/G^0_N $ the closure of the group 
		$ \phi_g(F_N) $ in Zariski topology of $\textnormal{GL}_n(K)$
		is contained in the group of points  $ G(K) \equiv \textnormal{spec}( K \otimes_k \mathcal{O}) $ 
		
	\end{enumerate}  
\end{lemma}
\noindent To verify the statement (i) recall the  properties of Razmyslov central polynomial (5.2) and use (6.11).  Statement (ii) follows from (i) and Lemma 6.2. Statements (iii) and (iv) are well known properties of algebraic sets over fields of characteristic zero (cf. e.g. \cite{Miln-ALG}, \cite{OV}) 

\paragraph{Step 4.}   Let $ C_N $ be a center of $ A_N  $. Set $  \mathfrak{p}_{g} =  \ker \phi_{g} \cap C_N = \ker (\phi_{g} | C_N ) $. If the ideal $    \mathfrak{p}_{g}$ is prime then denote by $\widetilde{...} $ the localization in $ A_N $ at   $ \mathfrak{p}_{g} $, i.e. for any  $ X \subset A_N $ set $ \tilde{X} = S^{-1}X \subset S^{-1}A_N = \widetilde{A}_N$ where $ S = C \setminus  \mathfrak{p}_{g} 
 $
\begin{lemma} For any $ g  \in G_N/G^0_N $
	\begin{enumerate}
		\item[(i)]  the ideal $\mathfrak{p}_{g} $  is an intersection of prime ideals of $ C_N $
			\item[(ii)] if the ideal  $  \mathfrak{p}_{g} =  \ker \phi_{g} \cap C_N = \ker (\phi_{g} | C_N ) $ is prime then  we have $ \widetilde{\ker \phi_{g}} =  \widetilde{\mathfrak{p}_{g}} \widetilde{A_N}$  
\item[(iii)]  if  $  K\phi_g(A_N) =  M_n(K) $ then   $ \mathfrak{p}_{g} $ is  a prime ideal
	\end{enumerate}	
\end{lemma}
\noindent Proof. Zariski closure of $\phi_g(F_N) $ in $ \textnormal{GL}_n(K) $ contains $  G^{0}(K) $  and is, therefore a reductive algebraic group (cf. e.g. \cite{Miln-ALG}). Hence,  as in the proof of Lemma 5.2 (i),  $ \phi_g | C_N $ is a homomorphism into a direct sum of copies of the field $K$. If, moreover, $  K \phi_g(F_N)  =  M_n(K) $ then $ \phi_g(C_N) $ is a scalar matrix and $ \phi_g | C_N $ is a homomorphism into a field. That proves (i) and (iii), The proof of (ii) is no different from the proof of the similar statement in Lemma 5.3  
\newline\newline
	Let  $ g  \in G_N/G^0_N $ be such that  
	$ R_n(\mathcal{M}_{i_1,g}, \cdots , \mathcal{M}_{i_t,g}) \neq 0 $ as in Lemma 6.3  (i) and Lemma 6.4 (ii). Denote by $ K_{g} $  the field of fractions of the integral domain  $  C_N/ \mathfrak{p}_{g} $ and let
$ \phi'_g : K_g \hookrightarrow K $ be an embedding of the  field 
		 $ K_{g} $  induced by the homomorphism $ \phi_g| C_N : C_N \rightarrow K  $.
The statement (ii) of the following lemma mimics similar arguments in \cite{KR}, \cite{R} that are due to Razmyslov (see also section 6.2)
\begin{lemma} The following statements are equivalent. 
	\begin{enumerate}
		\item[(i)] $ g  \in G_N/G^0_N $ satisfies condition (i) of Lemma 6.3
		\item[(ii)] There is a $K$-algebra isomorphism 
		$$ 
		 K \otimes_{\phi'_g(K_g)} A_N  \rightarrow  K\phi_g(A_N) = \textnormal{M}_n(K) $$ 
	\end{enumerate}
\end{lemma}
\noindent 
Proof.  It follows from (5.4) that dimension of $ \phi_g(A_N) $ over $ \phi'(K_g) $ is exactly $ n^2$  and therefore (i) implies (ii). The implication $ (ii) \implies (i) $ is trivial

\paragraph{Step 5} We can state now the following analogue of Lemma 5.5
\begin{lemma}There is a unique up-to-permutation,  irredundant set $ \mathfrak{p} = \{ 	\mathfrak{p}_1 , \cdots , \mathfrak{p}_{q} \} $ of  prime ideals of $ C_N$ such that 
	\begin{enumerate}
	\item [(i)] 
	\begin{align}
		\mathfrak{p}_1 \cap \cdots \cap \mathfrak{p}_{q} = \{0\}  
	\end{align}
	\item[(ii)]  the subset 
	\begin{align}
	  \{ 	\mathfrak{p}_{g} = \ker (\phi_{g} | C_N ) , \;   
	  g \in G_N/G^0_N
	 \} \subset \mathfrak{p}       \tag{6.12'}
	 \end{align}
	of ideals in $ \mathfrak{p} $ that satisfy equivalent conditions of Lemma 6.5 is nonempty 
	 
		\item [(iii)]  for any  $ g  \in G_N/G^0_N $ in the list (6.13') the closure of the group 
		$ \phi_g(F_N) $ in Zariski topology of $\textnormal{GL}_n(K)$
		is contained in the group of points  $ G(K) \equiv \textnormal{spec}( K \otimes_k \mathcal{O}) $ 
		
				\item [(iv)]  At least one   $ g  \in G_N/G^0_N $ satisfying the condition (ii) of Lemma 6.3 is contained in the list (6.13'); for such $ g $  the closure of the group 
		$ \phi_g(F_N) $ in Zariski topology of $\textnormal{GL}_n(K)$
		coincides with the group of points  $ G(K) \equiv \textnormal{spec}( K \otimes_k \mathcal{O}) $ 
		
	\end{enumerate}	
\end{lemma}
\noindent The proof of the statement (i)  is essentially no different from the proof of lemmas 5.4 and 5.5.  The statements (ii) and (iii) are covered by Lemma 6.3. The statement (iv) follows from Lemmas 6.2 and 6.3

\paragraph{Step 6.} 
  Turning to the irrep $ \sigma : H \rightarrow \textnormal{GL}(W) $  we can assume that $ \dim V = \dim W $.  Applying the same construction  we get  the map
\begin{align}
	\psi : kF_N \rightarrow \textnormal{M}_n(\mathcal{O}'_N), \;\;
	\psi( y_i ) = \mathcal{M}'_i, \; i = 1, \cdots, N  \nonumber
\end{align}
\noindent where $ N $ is large enough and $ \mathcal{O}'_N $ is a an appropriate coordinate ring of Cartesian power of the group $ \rho(H) $. Eventually, we will arrive at the same primary decomposition (6.12) of the central ring $ C_N$ and corresponding set of prime ideals    
	\begin{align}
	 \{ 	\mathfrak{q}_{h} = \ker (\psi_{h} | C_N ) , \;   
	h \in H_N/H^0_N
	\} \subset \mathfrak{p}       \tag{6.12''}
\end{align}
\noindent The sets of ideals (6.12') and (6.12'') coincide and there is a pair of indices $ g \in G_N/G^0_N $  and $  h \in H_N/H^0_N $   such that  
$ \mathfrak{p}_g = \mathfrak{q}_h $. Therefore, the fields of fractions  
   $ K_g $ of $ C_N/\mathfrak{p}_g $  and $ K_h$ of $ C_N/\mathfrak{q}_h =$ are the same. Since  $K$ is a universal extension there is  an automorphism $ \alpha_{g,h}$ of   $ K $ such that in notation of Lemma 6.5 $ \alpha_{g,h}( \phi_g'(K_g)) = \psi_h'(K_h) $. By Lemma 6.5 we have  induced  automorphism  of $ \textnormal{M}_n(K) $ 
(denoted by the ame letter)
 
\begin{align}
	 \alpha_{g,h} : K \otimes_{\phi'_g(K_g)} A_N  \approx K \otimes_{\psi'_h(K_h)} A_N   
\end{align} 
and, moreover
$ \widetilde{\ker \alpha_{g,h} \phi_{g} } = 
  \widetilde{\ \ker \psi_h } $, by Lemma 6.4 (ii).  We can now apply the arguments in   
  subsection 5.2.1 with the following modification.    
   It is easy to see that the automorphism $\alpha_{g,h} $ induces a homeomorphism in Zariski topology of the matrix algebra $ M_n(K) $ and, therefore it follows from Lemma 6.6 (iii)-(iv) that  $ G(K) = \overline{\phi_g(F_N)}  \subset H(K) $  for at least one of the corresponding pairs $ g, h $. For the same reason, there is a matching pair $ g', h' $ such that   $ H(K) = \overline{\psi_{h'}(F_N)}  \subset G(K) $. Hence $ G(K) = H(K) $ and since $ K $ and $ k $ are algebraically closed fields, we see that that  $  G(k) = H(k)$. This in turn implies that $ \sigma(H) $ is conjugate to $ \rho(G) $ in 
$ \textnormal{GL}_n(k) $ as stated by the Theorem 6.1.
\begin{remark}
	If $ G  $ is connected (cf. \cite{irrop}) then $C$ is an integral domain, $q = 1$ and both sets of ideals (6.12'), (6.12'') contain just one element $\{0\}$. In this case the fields $ K_g, K_h$ coincide with the field of fractions of $C_N$,  $A_N$ embeds into $ K \otimes_{\phi_g(C_N)} A_N $ and  $ \ker \phi_g = \ker \psi_h = 0 $.	If $G$ is finite then  the coordinate ring of polynomial functions on on $G$ is a ring of constant functions. In this case $K_g = K_h = k $ and the automorphism (6.13) is trivial
\end{remark}

\begin{remark}
	Thus, over algebraically closed field of characteristic zero, we have yet another proof of Theorem 1.1. Note that the central Laurent polynomial (5.1) was replaced in this proof by a direct splitting of  a coordinate ring of constant functions on a finite algebraic group   in characteristic zero (cf. \cite{Miln-ALG}). Most probably, these arguments  remain valid when $ \textnormal{char}(k) $ is relatively prime to $ | G/G^{0} |   $    
\end{remark}
\noindent Some corollaries and generalizations of Theorem 6.1 are discussed below
\subsection{Compact Lie Groups}
\begin{corollary}
Let $ G$ and $H$ be compact (real) Lie groups let $ \rho : G \rightarrow \textnormal{GL}(V) $ and $ \sigma : H \rightarrow \textnormal{GL}(W) $ be faithful finite dimensional unitary irreps with the same identities. Then representations $ \rho$ and $ \sigma$ are similar.
\end{corollary}	
\noindent Proof. As in the proof of Corollary 6.1 we can assume that $ G $ and $ H $ are subgroups of $ \textnormal{GL}_n(\mathbb{C}) $. 
In this case, it is well known (see e.g. \cite{OV}) that both $ G $ and $ 
H $ are real algebraic groups that are compact real forms  of 
$ G(\mathbb{C}) $ and  $ H(\mathbb{C}) $ respectfully. 
Now by Theorem 6.1, $ G(\mathbb{C})$ is conjugate to $  H(\mathbb{C}) $ in $ \textnormal{GL}_n(\mathbb{C}) $, since $ G $ (resp. $H$) is Zariski dense in $ G(\mathbb{C}) $ (resp. $ H(\mathbb{C}) $).
Finally, any two compact real forms of a reductive complex algebraic group are conjugate to each other (cf. \cite{OV}), and therefore irreps $\rho $ and $ \sigma$ are similar. 
\newline\newline\noindent Recalling  Remark 2.4 we can restate Corollary 6.3 as follows 
\begin{corollary}
Let $ G$ and $H$ be compact (real) Lie groups and let $ \rho : G \rightarrow \textnormal{GL}(V) $ and $ \sigma : H \rightarrow \textnormal{GL}(W) $ be faithful finite dimensional unitary irreps. Then the following conditions are equivalent
\begin{enumerate} 
	\item[(a)] 
$
 	\mathbb{E}(\rho)(tr(uu^{*})) = 0 \iff \mathbb{E}(\sigma)(tr(uu^{*})) = 0
 		\; \; \textnormal{ for } \!\!\textnormal{ any }  u \in kF  \nonumber 
$
	\item[(b)] 
$
	\mathbb{E}(\rho)(u) =  \mathbb{E}(\sigma)(u) 
	\; \; \textnormal{ for } \!\!\textnormal{ any }  u \in kF  \nonumber
$
\item[(c)] irreps $ \rho $ and $ \sigma$ are similar 
\end{enumerate}
In other words, $\rho$ and $ \sigma$ are not similar if and only if there is $ u \in kF $ such that
	$$ 	\mathbb{E}(\rho)(u) \neq  \mathbb{E}(\sigma)(u)  $$
\end{corollary}
\noindent The reason for mentioning Corollary 6.4 is the following result of Larsen and Pink (cf. \cite{LP}) that is stated below in a way convenient for this presentation 
\begin{theorema} (\cite{LP}, \cite{Larsen}).
	Let $ G$ and $H$ be compact connected semi-simple (real) Lie groups and let $ \rho : G \rightarrow \textnormal{GL}(V) $ and $ \sigma : H \rightarrow \textnormal{GL}(W) $ be faithful finite dimensional unitary irreps. Then the following conditions are equivalent
	\begin{enumerate} 	
		\item[(d)] 
		$$
	\underset{g  \in  G }{\int} \chi_{\rho}(g)^a \chi_{\rho}(g)^{*b}  d\mu_G	= 	\underset{h  \in  H }{\int}  \chi_{\sigma}(h)^a  \chi_{\sigma}(h)^{*b}   d\mu_H 
$$
 where $\mu_G , \mu_H $ are normilized Haar measures on $ G, H$ respectfully and $ a, b = 0,1, \cdots$ are arbitrary non-negative integers
	
		\item[(e)] irreps $ \rho $ and $ \sigma$ are similar 
	\end{enumerate}
\end{theorema} 
\noindent Thus for compact connected semi-simple Lie groups, conditions (a)-(e) are equivalent.  An interesting  possibility would be to find a computation that directly relates  conditions  (d), (b) and (a) without appealing to the condition (e)
%

\subsection{Algebraic Sets of Linear Operators} 
 As was already mentioned in  Remark 6.1, the proof of Theorem 6.1 can be significantly simplified  in case of connected  algebraic groups. Using connected algebraic groups as a main model, we would like to point out  that there is a  general principle that is applicable to finite dimensional irreps of irreducible (connected) algebraic sets over algebraically closed field of arbitrary characteristic. Let $ \rho_1 : S_1 \hookrightarrow \textnormal{End}(V_1), \; \rho_2 : S_2 \hookrightarrow \textnormal{End}(V_2) $ be a pair of irreducible algebraic subsets of full matrix algebras that are also irreducible as sets of linear operators. It should be clear how to define identities of a representation $ \rho : S \rightarrow \textnormal{End}(V) $ of any  (even unstructured) set $ S $. Namely, take a free associative algebra $ k(X) $ over an arbitrary field $k$ where $  X = \{x_1, x_2, \cdots \}$ is a countable set of free variables. A (non-commutative) polynomial $ p(x_1, \cdots, x_t ) \in k(X)  $ is said to be an identity of $ \rho$ if $ p(\rho(s_1), \cdots, \rho(s_t)) = 0 $ for any 
 $ s_1, \cdots , s_t \in S  $.  The following statement represents  a general pattern
 \begin{theorem}(cf. \cite{irrop}).
 Let the ground field $ k $ be algebraically closed. Suppose that $ S_1 \subset \textnormal{End}(V_1) $ and $ S_2 \subset \textnormal{End}(V_2) $ ($ \dim V_1 , \dim V_2 < \infty $) are  connected algebraic subsets that are  irreducible as sets of linear operators.  Under these conditions, representations $ S_1 \hookrightarrow \textnormal{End}(V_1), \; S_2 \hookrightarrow \textnormal{End}(V_2) $ have  same identities if and only if $ V_1 \approx V_2 \approx V $  and there is an invertible matrix $ M \in \textnormal{GL}(V) $  such that $ M^{-1}S_1 M = S_2$      
 \end{theorem}

 \noindent Theorem 6.2 is a generic statement and its generic proof  can be straightforwardly distilled from the proof of Theorem 6.1. Let's review  major steps of the proof:
 \begin{enumerate}
 	\item [1)]  Embedding of a factor-algebra (call it $A$) of a free algebra by the ideal of identities  into a full matrix algebra over coordinate ring (Lemma 6.1) 
 	\item [2)]  For connected sets of operators the  center $C$ of the algebra $A$ is an integral domain  (Remark 6.1) - see also Remark 6.7 below
 	\item [3)]  For an algebraically closed field $F$ containing $C$ the field extended algebra  $F\otimes_{C}A$ is isomorphic to the full matrix algebra over $F$ (Lemma 6.5 (iii)). This fact was established by Razmyslov (cf. \cite{KR}, \cite{R}) 
 	\item[4)] 
 	Correspondence between "schema point sets" $S_1(F), \; S_1(k) $  and $S_2(F), \; S_2(k) $  (cf. Step 6 in the proof of theoren 6.1). This is essentially equivalent to an application of Hilbert's Nullstellensatz suggested by Razmyslov in \cite{KR}, \cite{R}. A bit more can be said in a case of semi-simple Lie algebras in characteristic zero (see Remark 6.6 below)  
\end{enumerate}
We would like to make some additional (mostly unrelated) remarks on a rather general subject  of finite dimensional irreps of  "algebraic" structures with the same identities. Some results that are specific to Lie algebras and Lie groups will be discussed in section 6.3 below
\begin{remark}
	The sets $ S_1, S_2$ in Theorem 6.2 could be connected algebraic groups, finite dimensional vector spaces, finite dimensional Lie algebras (\cite{KR}, \cite{R}) or Jordan  algebras (cf. \cite{DR}). In each case one would require  embeddings involved  to  be morphisms of appropriate structures. It follows in particular, that two simple finite dimensional Lie/Jordan algebras that satisfy the same identities are isomorphic (cf. (\cite{KR}, \cite{R}, \cite{Zusmanovich}, \cite{DR})  
\end{remark}	
\begin{remark}
  Another approach to results of this kind  that relies on ultraproducts instead of central polynomials  can be found in  \cite{Zusmanovich}. 
\end{remark}
\begin{remark}
As was shown by Razmyslov (\cite{KR}, \cite{R}) the proof framework 1) - 3) of Theorem 6.2 works even for  infinite dimensional irreps of finite dimensional Lie algebras.  
A variation on this theme is outlined below in section 6.3
\end{remark}

\begin{remark}	

	Let $ k $ be an algebraically closed field of characteristic zero and let $ K $ be its algebraically closed  extension. Let $ A $ be an associative $K$-algebra and let $ L \subset A $  be 
	a Lie algebra over $ k $ with the Lie bracket inherited from $ A$, i.e. $ [x,y] = xy - yx $ for $ x,y \in L $. We will say that  $ A$ envelopes  $ L $ (or that $A$ is an envelope of $L$) if $ A$ is generated by $ KL $ as an associative algebra  
\newline\newline	 
\noindent \textbf{Lemma}.	
	Let $ A_1,A_2 $ be associative $K$-algebras that envelop  semisimple Lie $k$-algebras $  L_1, L_2 $.
	If there is a $K$-isomorphism  $ \phi : A_1 \approx A_2 $ such that $ \phi(K  L_1) = K  L_2   $ then it can be chosen in such a way that $ \phi(L_1 ) = L_2 $	    
\newline\newline
\noindent Proof. First of all, note that $ KL_1 = K \otimes_k L_1 $. Indeed, we have a $K$-algebra epimorphism 
$$ K \otimes_k L_1  \rightarrow K L_1 \rightarrow 0 $$ 
that leads to an exact sequence of Lie $K$-algebras
$$ 0 \rightarrow \mathfrak{a} \rightarrow K \otimes_k L_1 \rightarrow KL_1 \rightarrow 0  $$
The Lie algebra $ K \otimes_k L_1 $ is obviously semi-simple, hence the ideal     $ \mathfrak{a} $ contains an algebra $ K \otimes_k  L'_1 $  for one of the simple components $L'_1$ of the  $k$-algebra $ L_1 $. That is, however, impossible, unless $ L'_1 = 0 $.

 Thus we can assume that $ \phi( K \otimes_k L_1) = K \otimes_k L_2 $. Hence, $ K \otimes_k  \phi(L_1) = K \otimes_k L_2$, that is $ 
\phi(L_1)$ and $ L_2 $  are $k$-forms of 
the same semi-simple Lie algebra over  $ K$. Now it follows from  general theory of semi-simple Lie algebras that there is an inner   automorphism $ \alpha $ of $ K \otimes_k L_2 $ 
that maps $  \phi(L_1)  $ unto
$  L_2 $ (e.g. because the group of inner automorphisms acts transitively on the set of Cartan subalgebras).  Being an inner automorphism, $\alpha$  can be extended to an automorphism $ \alpha'$ of the associative envelope 
$ K  A_2 $ and the composition $ \alpha' \phi $ is an  isomorphism we are looking for
\end{remark}

\begin{remark}
Speaking informally, algebraic structures that are topologically connected do not satisfy nontrivial disjunctive identities.  One can  define   a connected algebraic system as  one with all its factors lacking nontrivial disjunctive identities. In case of linear algebras that seems to be equivalent to some kind of "primeness"  
\end{remark}

\subsection{Unitary Similar Representations} 

\begin{lemma} (cf. \cite{Sutherland}) 
	Unitary (respectfully orthogonal) finite dimensional irreps 
	$ \rho_1, \rho_2 $ of a group $ G $
	 are similar   if and only if they are unitary (respectfully orthogonally) similar
\end{lemma} 
\noindent Proof. 
Suppose  that unitary representations
 $ \rho_1, \rho_2 : G \rightarrow  \textnormal{GL}_n(\mathbb{C}) $  are similar
 (the case of orthogonal representations is no different). Then there is $ X \in \textnormal{GL}_n(\mathbb{C}) $ such that 
\begin{align}
 X  \rho_1( G) X^{-1}  =   \rho_2( G) 
\end{align}
Using the fact that  $ \rho_2$ is unitary, one gets  
$$   X \rho_1( g) (X^{*}X )^{-1}  \rho_1( g)^{-1} X^* = I_n $$ and 
\begin{align}
	   \rho_1( g) (X^{*}X )^{-1}  \rho_1( g)^{-1}  = (X^*X)^{-1}  \nonumber
\end{align}  
for all $ g \in G $. 
That means that positive definite matrix $ (X^*X) $ commutes with the representation 
$ \rho_1$. Assuming without  loss  of generality that representations  $ \rho_1 $ and $ \rho_2 $ are irreducible we see that $  (X^*X) = \lambda I_n$ for some positive real $\lambda$.  Now we can replace the matrix $X$ in (6.14) with the unitary matrix 
$ \lambda^{-1/2} X $. 
\newline\newline\noindent
Denote by $ \partial \rho $ a Lie algebra representation that is tangential to a Lie group representation $\rho$,  
\begin{theorem}
	Let $ G_1, G_2 $ be  connected simply connected real Lie groups and 
	let $$ \rho_i: G_i \rightarrow \textnormal{GL}(V_i), \; i = 1,2 $$ be 
	faithful, finite dimensional irreducible unitary representations. If corresponding Lie algebra representations 
	$ \partial \rho_i, \; i =1,2 $ have the same identities then $ \rho_1 $ and $\rho_2$ are unitary similar.  	
\end{theorem}
\noindent 
Proof. First of all, it should be clear that Lie groups $ G_1, G_2 $ and Lie algebras 
$ \mathfrak{g}_1,
\mathfrak{g}_2 $
are compact and semi-simple. (cf. \cite{Bourb9}). It is also clear that 
corresponding complex irreps of complexified Lie algebras
$$  \partial \rho_i : \mathbb{C} \otimes  \mathfrak{g}_i \equiv  \mathfrak{G}_i  \rightarrow \textnormal{End}(V_i)   , \; i =1,2 $$ 
also have the same identities. Therefore,  
complex irreps   $  \partial \rho_i : \mathfrak{G}_i  \rightarrow \textnormal{End}(V_i)   , \; i =1,2 $ are similar (cf. \cite{KR}, \cite{R})
and in particular,  complex Lie algebras $ \mathfrak{G}_1,
\mathfrak{G}_2 $ are isomorphic. It follows then that real Lie algebras $  \mathfrak{g}_1 $ and $ \mathfrak{g}_2$ are isomorphic as well, since a real compact form of a complex semi-simple  Lie algebra is unique up-to an interior automorphism (cf. \cite{Bourb9}). Thus,  we are dealing with two similar irreps
$ \partial \rho_i : \mathfrak{G}  \rightarrow \textnormal{End}(V)   , \; i =1,2 $ of the same complex Lie algebra $ \mathfrak{G} =  \mathbb{C} \otimes \mathfrak{g} $ and there is an invertible linear operator  $ X \in \textnormal{GL}(V) $ such that 
$ X \partial \rho_1(\mathfrak{G}) X^{-1}  = \partial \rho_2(\mathfrak{G}) $. 
Clearly, $  \mathfrak{g}' = X \partial \rho_1( \mathfrak{g}) X^{-1} $ is a  compact real form of 
$ \rho_2( \mathfrak{G}) $ and therefore (cf. \cite{Bourb9}) there is an inner automorphism  $\alpha$ of $ \mathfrak{G} $ such that $ \alpha(  \mathfrak{g}'  ) = \partial\rho_2( \mathfrak{g}) $. Since representation $ \partial \rho_2 $ is irreducible, the automorphism $ \alpha$ can be extended to an automorphism $ \alpha'$ of $\textnormal{End}(V)$  and combining automorphisms $ X $ and $ \alpha'$, we can assume that 
$$   X' \partial \rho_1( \mathfrak{g}) X'^{-1}  =  \partial \rho_2( \mathfrak{g})  $$
for some $ X' \in \textnormal{GL}(V)$.
Finally, for the connected simply connected Lie group $G$  this also implies
$$   X'  \rho_1( G) X'^{-1}  =   \rho_2( G)  $$
and unitary equivalence of $ \rho_1$ and $ \rho_2$   follows from Lemma 6.7
\newline\newline\noindent
Recalling Corollary 6.3 we get  
\begin{corollary}
	Let $ G_1, G_2 $ be  connected simply connected real Lie groups and 
	let $$ \rho_i: G_i \rightarrow \textnormal{GL}(V_i), \; i = 1,2 $$ be 
	faithful unirreps. 
	Then $ \rho_1, \; \rho_2 $ have the same identities if and only if Lie algebra irreps  $ \partial \rho_1, \; \partial \rho_2 $ have the same identities
\end{corollary}
\begin{remark}
	Corollary 6.5 strongly suggests that there is a one to one correspondence between identities of a unirrep of a connectet simply connected Lie group and identities of the corresponding irrep of its Lie algebra  
\end{remark}

\noindent 
\subsubsection{Infinite dimensional unitary representations}
Here we will roughly outline without proof a version of  the following result  established by Razmyslov (see \cite{KR}, \cite{R})
\begin{theorema} (\cite{KR}, \cite{R}). 
Let $ \rho_i : L_i \rightarrow \textnormal{GL}(V_i), \; i =  1,2 $  be faithful irreducible (not necessarily finite dimensional) irreps of finite dimensional Lie algebras over algebraically closed field. Let $A_i $ be an associative algebra spanned by $ \rho_i(L_i) $ in $ \textnormal{End}(V_i), \; i = 1,2$. If  identities of $ \rho_1 $ and $ \rho_2 $ are the same then there is an isomorphism $ \phi: A_1 \rightarrow A_2 $ such that $ \phi_1(L_1 ) = L_2 $
\end{theorema}	  
\noindent
Let G be a connected Lie group and let $ \rho : G \rightarrow \textbf{U}(H) $ be its unirrep (no nontrivial closed invariant subspaces) in a separable Hilbert space $H$. A vector $v \in H $ is called smooth (cf. e.g. \cite{Sch}) if the function 
$ g \rightarrow \rho(g)v  $ is in $ C^{\infty}(G) $. It is well known (cf. e.g. \cite{Sch})
that the subspace  $H^{\infty}_G \subset H $  of all $G$-smooth vectors in $H$ is everywhere dense in $H$. If $ \mathfrak{g} $ is a Lie algebra of $G$  
then a tangential action  of $ \mathfrak{g} $ on $ H^{\infty}_G $ can be defined (cf. e.g \cite{Sch}) so that there is a homomorphism 
$  \partial \rho : \mathfrak{g} \rightarrow L( H^{\infty}_G)$  of  $ \mathfrak{g} $ into an algebra of unbounded linear  operators on $  H^{\infty}_G $. 
Let  $\mathfrak{G} = \mathbb{C} \otimes \mathfrak{g} $ be a complexification of  $\mathfrak{g} $. The representation $  \partial \rho $ of $ \mathfrak{g} $ in $ H^{\infty} $ extends by linearity to that of  $\mathfrak{G} $ and we have a corresponding representation of the universal enveloping algebra $ U( \mathfrak{G}) $.    Denote the image $ \partial \rho(  U( \mathfrak{G})) $ by $ A_\mathfrak{g} $.
 
 \begin{theorem}
Let $ \rho_i  : G_i \rightarrow \textbf{U}(H) $ be faithful  unirreps of connected (real) Lie groups in a separable Hilbert space and let $\mathfrak{g_i} $ be corresponding Lie algebras. If  identities of representations $ \partial \rho_i : \mathfrak{g_i} \rightarrow L(H^{\infty}_{G_i}), \; i = 1,2  $  are the same then there is an isomorphism $ \phi : A_\mathfrak{g_1} \rightarrow  A_\mathfrak{g_2} $ such that $ \phi( \partial \rho_1(\mathfrak{g}_1 )) = \partial \rho_2(\mathfrak{g}_2 ) $   	
\end{theorem}
\begin{conjecture}
Let $ \rho_i : G_i \rightarrow \textbf{U}(H) $ be faithful  unirreps of connected simply connected Lie groups in a separable Hilbert space. Suppose that   identities of corresponding  Lie algebra representations $ \partial \rho_i : \mathfrak{g_i} \rightarrow L(H^{\infty}_{G_i}), \; i = 1,2  $  are the same and that  one of the following  conditions holds
\begin{enumerate}
	\item[(a)] $ G_1 $ and $ G_2 $ are nilpotent or
	\item[(b)] $ G_i $ are semisimple and $ \rho_i $ are discrete series representations  $ i = 1,2 $
	  
\end{enumerate}
Then $ \rho_1$ is unitary similar to $ \rho_2$
\end{conjecture}

\section{Appendix 1. Identities of the natural representation of $	\textnormal{SL}_2(k) $} 
 As an old curious example of an identity directly derived from the trace, we will reproduce here a proof of the following 
\begin{theorem} (\cite{sl2}). Over a field $k$ of characteristic zero all identities of the natural representation 
\begin{align}
	\textnormal{SL}_2(k) \hookrightarrow  \textnormal{GL}_2(k) 
\end{align}
follow from the identity 
	\begin{align}
	s_2(x,y) = ( y + y^{-1}) x - x (y + y^{-1}) 
	\end{align}  
\end{theorem}
\noindent The proof relies on the notion of \textit{trace identities} for matrix algebras that was introduced in \cite{RT} and \cite{PT1} (see also \cite{PT2}). Essentially, there is a one-to-one correspondence  between identities  of the representation (7.1) and \textit{trace identities}  of the matrix algebra $ \textnormal{M}_2(k) $. The main idea behind Theorem 7.1 can be explained as follows.
\begin{remark}
Write Hamilton Cayley identity for a matrix $ A \in M_2(k)$ (cf. (3.2)) 
\begin{align}
	A^2 - \textnormal{tr}(A) A + \det(A) I   = 0,
\end{align}
and assuming that $ \det(A) = 1 $ divide by $A^2$  to get 
\begin{align}
tr(A) I  = A + A^{-1} 
\end{align}
It was established in $\cite{RT}, \cite{PT1} $ that all the trace identities of the full matrix algebra follow from Hamilton Cayley identity. In particular, all the trace identities of $ \textnormal{M}_2(k) $ follow 
from the trace identity (7.3) and it turns out that formula (7.4) provides translation between trace identities of $ \textnormal{M}_2(k)$ and  identities of the natural representation of $\textnormal{SL}_2(k)$ (7.1). In other words, 
identity (7.2) is in some sense  equivalent to Cayley Hamilton identity (7.3). The proof presented below is a formal specification of this equivalence 
\end{remark} 

\paragraph{Proof of Theorem 7.1}
Recall following \cite{RT} the definition of trace identities. These identities live in a free trace algebra that can be defined as follows (cf. \cite{RT} or \cite{PT2}). Using countable set of free variables $X=\{x_1, \; x_2, \; \cdots \} $  and a formal functional symbol $ \text{Sp}(\cdots) $  define a semigroup $ S_0 $ by the following set of  rules
\begin{enumerate}
\item[-] the free semigroup $ S'(X)$ generated by $X$ belongs to $S_0$
	\item[-] $A,B \in S_0, \; a \in S'(X), \; a \neq \emptyset  \implies  \text{Sp}(AaB),    \text{Sp}(aB),  \text{Sp}(Aa) \in S_0 $  
	\item[-] $ A,B \in S_0 \implies AB \in S_0 $
\end{enumerate}
Let $ S \equiv S(X) $ be a maximal factor-semigroup of $ S_0$ that satisfies relations
\begin{itemize}
	\item[(a)] \text{Sp}(A)B = B \text{Sp}(A) 
	\item[(b)] \text{Sp}(AB) = \text{Sp}(BA) 
	\item[(c)] \text{Sp}(A\text{Sp}(B)) = \text{Sp}(B) \text{Sp}(A) 
\end{itemize}

\noindent It follows from relations (a)-(c) that  any element in $S $ can be uniquely written as    
\begin{align}
	 w_{0} \prod_{j=1}^{t_i}\text{Sp}(w_i), \; w_0, w_i \in S'(X),  \;   w_i \neq \oslash,  \; i = 1, \cdots, i; \;
\end{align}
Define the free algebra of trace polynomials as a semigroup algebra $ k S $ of $ S $ over the ground field $k $. 
\begin{remark}
Clearly, any element in $kS$ is uniquely represented by a linear combination of terms (7.5) up to commutat relations  (a) and (b) 
\end{remark} 
\noindent
An element 
$ P = P(x_1, \cdots, x_p ) \in kS $ is called a trace identity of the matrix algebra $ M_n(k) $ if it vanishes in $ M_n(k) $ for any substitution $ x_i \rightarrow a_i, \; a_i \in M_n(k), \; i = 1,  \cdots, p $ with the standard interpretation of the trace function $\text{Sp}$.  
\begin{remark}
	We will use 
	the linearization of the second degree trace identity $T_2(x)$ (3.2)  that, as it is easy to check, looks as follows
	\begin{align}
	kS \ni h_2(x,y) = xy + yx - x \text{Sp}(y) - y \text{Sp}(x) + 
		\text{Sp}(x) \text{Sp}(y) - \text{Sp}(xy)
	\end{align}
	Recall also that the standard  expression for the determinant of a $2 \times 2 $ matrix (cf. (3.2)) also can be viewed as an element of the algebra $kS$
	\begin{align}
	kS \ni	\det(a) = 1/2 ( \textnormal{Sp}(a)^2 -  \textnormal{Sp}(a^2) ) , \; a \in S(X)
	\end{align} 	
	We note in passing that there is a considerable amount of interest in specific trace identities even in small dimensions. See for example \cite{trace1} -\cite{trace4}.
\end{remark}
\noindent The following result is the key to what follows
\begin{theorema}(\cite{RT}, \cite{PT1})
	All the trace identities of $\textnormal{M}_n(k)$ follow from the Cayley Hamilton identity $h_2(x,y)$ (7.6).
\end{theorema}
\noindent For the set  (ideal) $ \mathcal{W}$ of the trace identities  of $ \textnormal{M}_n(k) $ in $kS$ the assertion of the theorem  is that 
 any element  $ w \in \mathcal{W}$  can be written as 
\begin{align}
	w = \sum_{i}A_i h_2(C_i,D_i)B_i, \; A_i, B_i, C_i, D_i \in kS 
\end{align}
where $ h_2 $ is a linearized Cayley Hamilton polynomial (7.6) and as explained above
\begin{align}
		C_i = c_{i0} \prod_{j=1}^{t_i}  \text{Sp}(c_{ij}), 
	\; 	D_i = d_{i0} \prod_{j=1}^{s_i}  \text{Sp}(d_{ij}) 	 \\
		h_2(C_i, D_i) = h_2(  c_{i0},d_{i0} ) 
	\prod_{j=1}^{t_i}  \text{Sp}(c_{ij}) 	\prod_{j=1}^{s_i}   \text{Sp}(d_{ij}) \\
	c_{i0}, \; d_{i0}, \; c_{ij}, \; d_{ij} \in S'(X); \;   d_{ij}, c_{ij} \neq \emptyset  \textnormal{ for }  j > 0  \nonumber
\end{align}
Let's turn now to the $ \textnormal{SL}_2(k) $ representation (7.1).
Let $ U $ be the ideal  in $kF = kF(Y) $ that is a set of all elements of the form
\begin{align}
	kF \ni u s_2(f,g) v, \; u,v \in kF, \; f,g \in F 
\end{align}
In other words (see \cite{P}),  $ U $ is defined as a minimal \textit{verbal} (cf. \cite{P} and references therein) ideal spanned by $ s_2$. We need to show that any identity of the representation (7.1) belongs to $ U $. To this end we construct a linear map
$ \phi : kS \rightarrow kF/U $ as follows:

\begin{enumerate}	
	\item[(a)] for $ w = w(x_1, \cdots, x_t) \in kS' $ set 
	$ \phi( w ) = w' \equiv w(y_1,\cdots,y_t ) \mod U  $ 
	
	\item[(b)]  and set $ \phi(   \text{Sp}(w) ) =
	w' + w'^{-1} \equiv w(y_1,\cdots,y_t ) + w(y_1,\cdots,y_t )^{-1} \mod U $
	\item[(c)]   extend thus defined map $\phi$ to the whole $kS$ by multiplicativity and linearity  
\end{enumerate}
\begin{lemma} The map $\phi $ is correctly defined by the rules (a)-(c). Moreover,  $ \mathcal{W}\subset \ker \phi $
\end{lemma}
Proof. First, use Remark 7.2,  to check that definitions (a)-(c) do not depend on the order of trace terms:

\begin{align}
\phi( Sp(ab) - Sp(ba)) \; = \; a'b' +(a'b')^{-1}  -  
	 b'a' - (b'a')^{-1}  =  \nonumber \\ 
\;\;\;\;\;\;	  \{(a'b' +(a'b')^{-1} )b' -  
	 b'(b'a' +(b'a')^{-1} ) \} b'^{-1}
	 \in U  \nonumber \\
\!\!\!\!\!\!	\phi( Sp(a)b - bSp(a)) = (a + a'^{-1})b' - 
	 b'(a + a'^{-1}) \in U   
\end{align}
Now, let's show that $ \mathcal{W}\subset \ker \phi $. In accordance with (7.8)-(7.10) and (7.12), we need
to verify
that $ 	\phi(h_2(c,d)) \in U $ for any $c,d \in S(X) $.  Indeed,
\begin{align}
	\phi(h_2(c,d)) = c'd' + d'c' -c'(d' + d'^{-1}) - 
	d'( c' + c'^{-1}) + ( c' + c'^{-1})(d' + d'^{-1}) = \nonumber \\ = - \; c'd' - (c'd')^{-1} = -( d' + d'^{-1} ) c'^{-1} +  c'^{-1}( d' + d'^{-1} ) \;  
	\in \; U  \nonumber
\end{align} 
It follows from Lemma 7.1 that to prove Theorem 7.1 it is sufficient to show that for any  identity $ v = v(y_1, \cdots, y_t ) \in kF $ of the representation (7.1) there is $ w = w(v)  \in \mathcal{W} $ such that $  \phi(w) = v  \!\!\!\mod U $. We precede a construction of $ w(v) $ by a technical lemma.

 For any $ f = f(y_1, \cdots, y_t ) \in F $ 
and any $ y_i \in Y $ let $ \deg_{y} (f) $  denote the degree of the variable $ y \in Y  $ in a reduced representation of $ f \in F(Y)$. For example 
$$ \deg_{y_1}( [y_1, y_2] ) = 0,  \;  \deg_{y_1}(y_1^{-2}y_3 y_1) = -1, \; \deg_{y_1}(y_2) = 0 $$.
\begin{lemma}
	Any identity
	$$ kF \ni v  = \sum_{i=1}^{m} \alpha_i f_i, \; \alpha_i \in k , \; f_i = f_i(y_1, \cdots, y_t) \in F $$
	of the representation (7.1) can be written as $ v_1 + \cdots + v_k , \; v_i \in kF, \; i = 1, \cdots , k  $ in such a way that 
	\begin{enumerate}
		\item[(i)] all $ v_i $ are identities of (7.1)
		\item[(ii)] if $$ v_i = \sum_{j=1}^{m_i} \alpha_{ij}f_{ij}(y_1, \cdots , y_t), \;  \alpha_{ij} \in k, \; f_{ij} \in   F(Y), \; i = 1, \cdots , m $$  
		 then  for any $ y \in Y, \;  \deg_{y} (f_{ip}) \equiv  \deg_{y} (f_{iq}) \!\!\!\mod 2, \; 1 \leq p \leq q \leq t $
	\end{enumerate}
		
\end{lemma}
\noindent Proof. Write $ v = v_1 + v_2 $ where degree of the variable $y_1$ in any summand of $ v_1 $  is even and 
 degree of $y_1$ in any summand of $ v_2 $ is odd. Suppose that there are $ g_1, \cdots, g_t \in \textnormal{SL}_2(k)$ such that in the matrix algebra $ \textnormal{M}_2(k) $
  $$ v_1( g_1, \cdots, g_t ) = \beta \neq 0 $$
Then setting $ v_2( g_1, \cdots, g_t ) = \gamma  $,  we have 
 \begin{align}
 0 = v(g_1, \cdots, g_t) =  v_1( g_1, \cdots, g_t ) +  v_2( g_1, \cdots, g_t ) = \beta + \gamma 
 \end{align}  
Obviously, $ -g_1 \in \textnormal{SL}_2(k)$, since $ \det( -g_1 ) = (-1)^2 \det(g_1) = 1 $ and we can flip the sign of $ g_1 $ in (7.13), thus getting 
\begin{align}
 0 = v( -g_1, \cdots, g_t) =  v_1( -g_1, \cdots, g_t ) +  v_2( -g_1, \cdots, g_t ) = \beta - \gamma  \nonumber
\end{align}
Hence, $ \beta = \gamma = 0 $. Therefore both  $ v_1 $ and $ v_2 $ must be identities of representation (7.1). The statement of the Lemma 7.2 can be now easily verified by recursively applying the same argument to $v_1 $ and $ v_2$ separately.    
\newline
\newline
\noindent 
Let
\begin{align}
	kF \ni v = v(y_1, \cdots, y_t ) = \sum_{i=1}^m \alpha_i f_i(y_1, \cdots, y_t ), \; \alpha_i \in k , \; f_i = f_i(y_1, \cdots, y_t) \in F 
\end{align}
\noindent be an identity of (7.1). Let $ c_{ij}  = \deg_{y_j}(f_i) $. 
 By Lemma 7.2 we can assume that for a fixed $j$ all $ c_{ij}$ are either simultaneously even or simultaneously odd. In the latter case, multiply $v$ by $y_j$ (for all such $ j $) thus obtaining an identity of representation (7.1) that is equivalent to (7.14). Thus, without loss of generality we can assume that all variable degrees $ c_{ij} $ are even.

 Let's  describe a (textual) transformation $ f \rightsquigarrow \hat{f} $ from $ kF $ to $ kS $ that will be used  to find $ w(v) \in \mathcal{W} $ such that $ \phi(w) = v$. First, replace any inverse of a free variable 
 $ y^{-1}_j \in Y $ by $  \text{Sp}(x_j) - x_j $ and replace 
 every variable $ y_j $ by $ x_j $. For example,
 $$ \widehat{[y_1,y_2]} = (\text{Sp}(x_1) - x_1)(\text{Sp}(x_2) - x_2) x_1 x_2 $$ 
 Next, using (7.7) as a definition of the expression $ \det(x) \in kS $, set
 \begin{align}
 kS \ni	d_i  \; \equiv \;  	\prod_{l=1}^t \det(x_l)^{(1/2)c_{il}} , \;\; i = 1, \cdots,  m  \nonumber
 \end{align} 
\noindent and divide the term $\hat{f}_i$ by the "square root of its determinant"  $d_i$ to get a formal expression
 \begin{align}
 	(1/d_i)\hat{f}_i  \; \equiv \;  \left(	\prod_{l=1}^t \det(x_l)^{-(1/2)c_{il}} \right)\!\!\hat{f}_i, \;\; i = 1, \cdots,  m  \nonumber
 \end{align} 
 Now formally compute
 \begin{align}
 	\sum_{i=1}^m \alpha_i \frac{1}{d_i}\hat{f}_i \; \equiv \; 
 		\frac{\overset{m}{\underset{i=1}\sum} \alpha_i (\prod_{j \neq i}d_j)\hat{f}_i }{d_1 d_2 \cdots d_m}   
 \end{align}
The numerator of this expression is an element in $kS$ that we are looking for. Set
 \begin{align}
 	 w \equiv w(v) = \overset{m}{\underset{i=1}\sum} \alpha_i (\prod_{j \neq i}d_j)\hat{f}_i    
 \end{align} 
\begin{lemma} If $ v $ is an identity of (7.1) then  
	\begin{enumerate}
		\item[(i)] $ w = w(v) $ (7.16) is a trace identity of $M_2(k)$, i.e. $ w \in \mathcal{W} $ 
		\item[(ii)] and $ \phi(w) = v  \mod U $  
	\end{enumerate}	
\end{lemma}
\noindent Proof. Let $ \bar{k} $ be an algebraic closure of the field $k$. By Zariski topology arguments  $ v $ is an identity of the representation $ 	\textnormal{SL}_2(\bar{k}) \hookrightarrow \textnormal{GL}_2(\bar{k}) \subset \textnormal{M}_2(\bar{k}) $.
Taking into account that $ \det(a)^{-1/2}\; a \in SL_2(\bar{k}) $ for any $ a \in GL_2(\bar{k}) $, we see (cf. (7.15),(7.16)) that the right hand side of the expression (7.16) vanishes in $ \textnormal{M}_2(\bar{k}) $ for any assignment of invertible matrices to its free variables. Since $ \textnormal{GL}_2(\bar{k}) $ is (Zariski) dense in $ M_2(\bar{k}) $  that proves (i). 
To prove (ii), note that by definition of the map $\phi$
\begin{align}
\!\!\!\!\!\!\!\!\!\!\!\!\!\!\!\!\!\!\!\!\!\!\!\!	\phi( \text{Sp}(x) - x ) = y + y^{-1} - y = y^{-1} \!\!\!\! \mod U \nonumber 
\end{align}
and (cf. Remark 7.2) 
\begin{align}
	\phi( \det(x)) = \frac{1}{2}( \phi(\text{Sp}(x)^2  -  \phi(\text{Sp}(x^2))  = \frac{1}{2}( y^2 + y^{-2} + 2 - y^2 - y^{-2} )  = 1 \!\!\!\! \mod U  \nonumber	
\end{align}    
\newline

\section{Appendix 2. Examples of similar representations of finite Abelian $p$-groups}
The ground field in this section is again $\mathbb{C}$. We  begin with the following, most probably well known  example.

\begin{proposition}
Faithful $n$-dimensional representations  of a finite Abelian $p$-group of rank $n$ are similar. 	
\end{proposition}
\noindent Proof. 
Faithful one-dimensional representations of a cyclic $p$-group are obviously similar. Let $ \rho  : A \hookrightarrow  \textnormal{GL}(V), \;  \sigma : A \hookrightarrow  \textnormal{GL}(V), \; \dim V = n  $ be two faithful representations of a finite Abelian $p$-group $A$ of rank $ n$. Thus, there are linear characters
$ \rho_i, \sigma_i  \in \hat{A}, \; i = 1, \cdots , n $ such that  
$ \rho(g) = \textnormal{diag}(\rho_1(g), \cdots \rho_n(g)) $ and 
$ \sigma(g) = \textnormal{diag}(\sigma_1(g), \cdots \sigma_n(g)) $ are diagonal matrices for any $ g \in A $. 
Let $ p^m $ be a maximal order  among  cyclic groups $ \rho_i(A), \; i = 1, \cdots, n $. From faithfulness of representations $ \rho, \; \sigma$  it follows that
$p^m$ is a maximal order of elements in $ A $ and there are 
indices $ 1 \leq k, \; l \leq n $ such that $ \rho_k(A) \approx \sigma_l(A) \approx  \mathbb{Z}_{p^m} $. 
 
\noindent Assume that $ A $ is not cyclic and let    $ H = \ker \rho_k $ and $ K = \ker \sigma_l  $. Then there are 
elements $ a, b \in A $, both of order $ p^m $ such that  
$$ <a> \times \; H \; \approx \; A \; \approx \;\;  <b> \times \; K $$ 
It is easy to see that $ H \approx K $ hence there is an automorphism $ \alpha : A \rightarrow A $ that maps 
$ b $ into $a$ and $ K $ into $H$.  After permuting diagonal elements  of  $ \sigma \alpha(A) $  if necessary we get a representation $ \sigma' $ similar to $ \sigma$   and corresponding decomposition of representation space $ V = V_0 \oplus V_1 $ such that 
\begin{enumerate}
	\item[(1)] $\dim V_0 = 1$
	\item[(2)] $ V_0 $ and $ V $ are  invariant under $ \rho $ and $ \sigma'$ 
	\item[(3)] $ \rho\;|<a> \; = \; \sigma' \;| <a> $ is trivial on $ V_1$
	\item[(4)] both $ \rho|H$ and $ \sigma' |H $ are trivial on $ V_0 $ and exact on $ V_1$ 
\end{enumerate}
The rank of $H$ is clearly $ n-1$ and we can now apply the same arguments to representations $ \rho|H$ and $ \sigma' |H $ to finish the proof by induction 
\begin{example}
Let $ A = \mathbb{Z}_p^m $. Any $n$-dimensional representation $\rho$ of $A$ can be parametrized (see e.g. \cite{Fourier}) as 
\begin{align}
	\rho(a) \equiv \rho_V(a)= \exp \left(\frac{2\pi i}{p} Va  \right), \;\; a \in A 
\end{align}
where $V \in \textnormal{End}(\mathbb{Z}_p^m, \mathbb{Z}_p^n) $ and  we use convention (3.17-3.18') for an exponent of the column vector $Va$ . If  $ m = n $ then 
$\rho_V$ is faithful if and only if the matrix $ V $ is invertible. In this case any faithful $A$-representations $ \rho_U $ and $ \rho_V $ are similar as we have 
$ \rho_U = \rho_V (UV^{-1}) $  
\end{example}

\begin{example}

 Let $G = \mathbb{Z}_p \times \mathbb{Z}_p, \; p > 2 $
and consider two $3$-dimensional representations $ \rho_1 $ and $ \rho_2 $ of $G$ that are parametrized (cf. (8.1)) by $ 3 \times 2 $ matrices 
\begin{gather}   
V_1 =
	\begin{bmatrix}
		1 & 0 \\
		0 & 1  \\ 
		1 & 1	 	
	\end{bmatrix} , \;
V_2 = 		\begin{bmatrix}
		1 & \;\; 0 \\
		0 & \;\; 1  \\ 
		1 & -1
	\end{bmatrix}  \;\;   \in \; \;  \textnormal{End}(\mathbb{Z}_p^2, \mathbb{Z}_p^3)  \nonumber
\end{gather}
that is for $(a,b) \in G$ we have 
$$ \rho_1((a,b)) = 
  \exp \left(\frac{2\pi i}{p}(a,b, a+b)^T   \right) $$
and
$$ \rho_2((a,b)) = 
\exp \left(\frac{2\pi i}{p}(a,b, a-b)^T   \right) $$ 
 An automorphism $ \gamma$ of the group $ G $ defined by the rule $ (a,b) \rightarrow (a+b, b  )$ satisfies the equation $$ \rho_2(\gamma(a,b)) ) = \tau  \rho_1(a,b) $$   
where $ \tau $ is an automorphism of $ \mathbb{C}^3 $ that transposes first and third coordinates. Hence $ \rho_1 $ and $ \rho_2 $ are similar. 

\end{example}


\begin{thebibliography}{99}
	
	\bibitem{P}B. I. Plotkin. Varieties of representations of groups,  Russian Math. Surveys 32:5 (1977), 1-72  

	
	\bibitem{PK} Plotkin, B. I. and Kushkuley, A. H., Identities of regular representations of groups, 
	preprint (1979) (unpublished)
	
	\bibitem{Cohn} Paul M. Cohn, Universal Algebra, Harper \& Row 1965, Springer 1981


	
	\bibitem{V} Samuel M. Vovsi, DISJUNCTIVE IDENTITIES OF FINITE GROUPS, arXiv:math/9507206v1, 1995
	
		
	\bibitem{Atyah} M. F. ATIYAH and D. O. TALL, 	GROUP REPRESENTATIONS, $\Lambda$-RINGS AND THE
	J-HOMOMORPHISM, Topology Vol. 8, pp. 253-297, 1969
	
	\bibitem{Serr}Serre, Jean-Pierre, Représentations linéaires des groupes finis, Hermann, 1971 
	
	\bibitem{Serre2} Jean-Pierre Serre, Finite Groups: An Introduction, International Press, 2016

	
	\bibitem{Wildon_L} Mark Wildon,	LABELLING THE CHARACTER TABLES OF SYMMETRIC
	AND ALTERNATING GROUPS, 	arXiv:math/0612292, 2006
	\bibitem{Dade} E.C. Dade,  Answer to a question of R. Brauer. J. Algebra 1 (1964) 1–4.
	\bibitem{Meir} Ehud Meir and Markus Szymik, Adams operations and symmetries
	of representation categories, arXiv:1704.03389v3, 2019
%
	\bibitem{GW} Roe Goodman,
	Nolan R. Wallach, 
	\url{https://sites.math.rutgers.edu/~goodman/repbook.html}, Appendix G, Group Representations and the Platonic Solids, 2009

	\bibitem{S4} \url{https://people.maths.bris.ac.uk/~matyd/GroupNames/1/S4.html}
	
	\bibitem{irrop}
	A. H. Kushkuley, On identities of  irreducible sets of linear operators,  Izv. Visch. Uch. Zaved., vol. 12, 1987, pp.  55-58
	
	\bibitem{Miln-ALG}James S. Milne,  Algebraic Groups, Lie Groups, and their Arithmetic Subgroups, 2011, available at   \url{https://www.jmilne.org/math/}
	\bibitem{Miln-CA} James S. Milne, A Primer of Commutative Algebra, 2020, available at \url{https://www.jmilne.org/math/}
	\bibitem{Miln-ANT}  James S. Milne, Algebraic Number Theory (v3.08), 2020, available at \url{https://www.jmilne.org/math/}
	\bibitem{Artin}Michael Artin, NONCOMMUTATIVE RINGS, class notes, Math 251, Berkeley,  1999

	\bibitem{Seg} Gerald H. Cliff,  Surinder K . Sehgal, ON GROUPS HAVING THE SAME CHARACTER TABLES, COMMUNICATIONS IN ALGEBRA, 9 ( 6 ) , 627-640 (1981)
	\bibitem{Nenciu}Adriana Nenciu, Character tables of p-groups with derived subgroup
	of prime order I, Journal of Algebra 319 (2008) 3960–3974
	\bibitem{Davydov} A.A.Davydov, Finite groups with the same character tables, Drinfel'd algebras and Galois algebras, arXiv:math/0001119v1, 2000


\bibitem{RT} Ju. P. Razmyslov, Identities with trace in full matrix algebras over a field of characteristic zero, Izv. Akad. Nauk SSSR Ser.
Mat. 38, 723–756 (1974) (Russian); English Transl. Math. USSR, Izv. 8, 727–760 (1975)
\bibitem{PT1} C. Procesi, The invariant theory of n × n matrices, Adv. Math. 19
(1976), 306–381
\bibitem{PT2} C. Procesi, T–ideals of Cayley Hamilton algebras,  	arXiv:2008.02222, \url{https://arxiv.org/pdf/2008.02222.pdf}, 2021

\bibitem{sl2} A. H. Kushkuley, On identities of finite-dimensional representations, Proc. Riga Algebra
Sem. 1977, no. 3, 91-115

\bibitem{trace1}Hong-Hao Zhang, Wen-Bin Yan, and Xue-Song Li, Trace Formulae of Characteristic Polynomial and Cayley-Hamilton’s Theorem, and
Applications to Chiral Perturbation Theory and General Relativity, 	arXiv:hep-th/0701116v1, \url{http://arXiv.org/abs/hep-th/0701116v1}, 2007

\bibitem{trace2} DRAGOMIR Ž. $\DJ$OKOVIĆ AND BENJAMIN H. SMITH, QUATERNIONIC MATRICES: UNITARY SIMILARITY,
SIMULTANEOUS TRIANGULARIZATION AND SOME
TRACE IDENTITIES, arXiv:0709.0513v1 ,  \url{http://arxiv.org/abs/0709.0513v1}, 2007

\bibitem{trace3}S. Humphries and C. Krattenthaler,
TRACE IDENTITIES FROM IDENTITIES FOR DETERMINANTS, arXiv:math/0411061v1,
\url{http://arxiv.org/abs/math/0411061v1}, 2004

\bibitem{trace4} T. H. Marshall, G. J. Martin,
Polynomial Trace Identities in SL(2, C),
Quaternion Algebras, and Two-generator
Kleinian Groups,  arXiv:1911.11643v1,  \url{http://arxiv.org/abs/1911.11643v1}, 2019

\bibitem{TF} \url{https://en.wikipedia.org/wiki/Characteristic_polynomial}
\bibitem{TF1} \url{https://en.wikipedia.org/wiki/Faddeev%E2%80%93LeVerrier_algorithm}

\bibitem{Adams}
\url{https://en.wikipedia.org/wiki/Adams_operation}

\bibitem{Phase} Hartmut Führ,
Vignon Oussa, Phase Retrieval for Nilpotent Groups, Journal of Fourier Analysis and Applications, 29:47, \url{https://link.springer.com/article/10.1007/s00041-023-10031-5}, 2023

\bibitem{Pointwise}JOEL SEGAL, POLYNOMIAL INVARIANT RINGS ISOMORPHIC AS
MODULES OVER THE STEENROD ALGEBRA, J. London Math. Soc. (2) 62 (2000) 729–739

\bibitem{conj}Silvio Dolfi, Gabriel Navarro, Pham Huu Tiep, FINITE GROUPS WHOSE
SAME DEGREE CHARACTERS ARE GALOIS CONJUGATE, ISRAEL JOURNAL OF MATHEMATICS 198 (2013), 283–331

\bibitem{Stone}MOHAMMAD BARDESTANI, KEIVAN MALLAHI-KARAI, AND HADI SALMASIAN, MINIMAL DIMENSION OF FAITHFUL REPRESENTATIONS FOR p-GROUPS, arXiv:1505.00626v3, 2016

\bibitem{Howe} Roger E. Howe. Some Recent Applications of Induced Representations. In Group Representations, Ergodic
Theory, and Mathematical Physics: a Tribute to George W. Mackey, volume 449 of Contemp. Math., pages 173–191. Amer. Math. Soc., Providence, RI, 2008.

\bibitem{Heisenberg} SHANE CERNELE, MASOUD KAMGARPOUR, AND ZINOVY REICHSTEIN, MAXIMAL REPRESENTATION DIMENSION OF FINITE
p-GROUPS,  arXiv:0911.0637v3, \url{http://arxiv.org/abs/0911.0637v3}, 2010

\bibitem{Bump} D. Bump. Automorphic forms and representations, Cambridge University Press, Cambridge, 1997


\bibitem{JO} 
Chung-Nim Lee and Arthur G. Wasserman,
On the groups JO(G),  Memoirs of the American Mathematical Society, Number 159, 1975

\bibitem{S} Wujie Shi, A counterexample for the conjecture of finite simple groups, arXiv:1810.03786v1, 2018

\bibitem{S1}
A. V. Vasil’ev, M. A. Grechkoseeva, and V. D. Mazurov, Characterization of the finite simple groups by spectrum and order, Algebra Logic 48 (2009), no. 6, 385–409.

\bibitem{S2} Maria A. Grechkoseeva, Victor D. Mazurov, Wujie Shi, Andrey V. Vasil’ev,
and Nanying Yang, FINITE GROUPS ISOSPECTRAL TO SIMPLE GROUPS,
arXiv:2111.15198v3 , \url{http://arxiv.org/abs/2111.15198v3},  2022


\bibitem{Mackey} G. W. Mackey, Ergodic group actions with a pure point spectrum, Illinois J.
Math. 8 (1964), 593–600

\bibitem{Sutton} Craig J. Sutton, Isospectral simply-connected homogeneous spaces
and the spectral rigidity of group actions, arXiv:math/0301376v1,  \url{http://arxiv.org/abs/math/0301376v1}, 2003

\bibitem{Todd} J. A. Todd, On a conjecture in group theory, J. London Math. Soc. 25 (1950), 246

\bibitem{Ikeda} Akira Ikeda, 
On spherical space forms which are isospectral
but not isometric, J. Math. Soc. Japan
Vol. 35, No. 3, 1983

\bibitem{Ikeda_1} Akira Ikeda, ON SPACE FORMS OF REAL GRASSMANN MANIFOLDS WHICH ARE ISOSPECTRAL BUT NOT ISOMETRIC, KODAI MATH. J.
20 (1997), 1-7

\bibitem{Vincent} G. Vincent, Les groupes lineaires finis sans points fixes, Comment. Math. Helv., 20. (1947), 117-171

\bibitem{Wolf}
J. Wolf, Spaces of constant curvature, McGraw Hill, 1967

\bibitem{Wolf2} Joseph A. Wolf, Isospectrality for Spherical Space Forms, Result. Math. 40 (2001) 321-338


\bibitem{isospectral-CR} GERSON GUTIERREZ, EMILIO A. LAURET, AND JUAN PABLO ROSSETTI, ISOSPECTRAL CR MANIFOLDS WITH RESPECT TO THE KOHN
LAPLACIAN, arXiv:2312.15214v1, \url{http://arxiv.org/abs/2312.15214v1}, 2023

\bibitem{Sunada} Toshikazu Sunada, Riemannian Coverings and Isospectral Manifolds, The Annals of Mathematics, 2nd Ser., Vol. 121, No. 1, 1985, pp. 169-186

\bibitem{Sutton} CRAIG J. SUTTON, EQUIVARIANT ISOSPECTRALITY AND SUNADA’S METHOD,  arXiv:math/0608557v2 ,
\url{http://arxiv.org/abs/math/0608557v2},2010

\bibitem{Sutherland} ANDREW V. SUTHERLAND, ARITHMETIC EQUIVALENCE AND ISOSPECTRALITY, Lecture Notes, 2018


\bibitem{Allcock}DANIEL ALLCOCK, SPHERICAL SPACE FORMS REVISITED,
arXiv:1509.00906v3, \url{http://arxiv.org/abs/1509.00906v3}, 2016 

\bibitem{Gilkey}PETER B. GILKEY
On spherical space forms with meta-cyclic
fundamental group which are isospectral but
not equivariant cobordant,
Compositio Mathematica, tome 56, no 2 (1985), p. 171-200

\bibitem{Howe1} SHAMGAR GUREVICH AND ROGER HOWE, A LOOK AT REPRESENTATIONS OF $ \textnormal{SL}_2(F_q)$ THROUGH THE LENS OF SIZE,
arXiv:2105.11008v1,  2021

\bibitem{commprob}
TIMOTHY C. BURNESS, ROBERT M. GURALNICK, ALEXANDER MORETÓ,
AND GABRIEL NAVARRO,
ON THE COMMUTING PROBABILITY OF p-ELEMENTS
IN A FINITE GROUP,
arXiv:2112.08681v2, 2022

\bibitem{Marin} I. MARIN AND J. MICHEL, AUTOMORPHISMS OF COMPLEX REFLECTION GROUPS, REPRESENTATION THEORY,
An Electronic Journal of the American Mathematical Society
Volume 14, Pages 747–788 (December 14, 2010)

\bibitem{S5} \url{https://groupprops.subwiki.org/wiki/Linear_representation_theory_of_symmetric_group:S5}
\bibitem{S6}
S. G. Ngulde, F. B. Ladan, M. M. Karagama, IOSR Journal of Mathematics,
Volume 14, Issue 1 , PP 43-48, 2018
\bibitem{M11} \url{https://math.stackexchange.com/questions/4856291/distinct-characters-with-the-same-character-values-outer-automorphisms-and-galo}
\bibitem{Hughes} Sam Hughes,
Representation and Character Theory
of the Small Mathieu Groups, Master’s thesis, University of South Wales,
\url{https://nickpgill.github.io/MMath_Sam.Hughes.pdf}

\bibitem{Stanley}
RICHARD P. STANLEY, INVARIANTS OF FINITE GROUPS AND THEIR
APPLICATIONS TO COMBINATORICS, BULLETIN (New Series) OF THE
AMERICAN MATHEMATICAL SOCIETY
Volume 1, Number 3, May 1979

\bibitem{Smith} LARRY SMITH, POLYNOMIAL INVARIANTS OF FINITE GROUPS
A SURVEY OF RECENT DEVELOPMENTS, BULLETIN (New Series) OF THE
AMERICAN MATHEMATICAL SOCIETY
Volume 34, Number 3, July 1997, Pages 211–250

\bibitem{JSegal}Joel Segal, Pointwise Conjugate Groups and Modules over
the Steenrod Algebra,  PhD Thesis, Göttingen 1999

\bibitem{Sandlinng} Sandling, R., Units in the modular group algebra of a finite abelian p-group, J. Pure Appl.
Algebra 33 (1984), 337–346.

\bibitem{Bovdi}  ADALBERT BOVDI, The group of units of a group algebra of characteristic $p$,  Publ. Math. Debrecen
52 / 1-2 (1998), 193–244

\bibitem{circulant} \url{https://en.wikipedia.org/wiki/Circulant_matrix}

\bibitem{Liebeck} H. Liebeck, Concerning nilpotent wreath products, Proc. Cambridge Philos. Soc. 58
(1962) 443–451
\bibitem{Swa} Ed Swartz, Matroids and quotients of spheres, Math. Z. 241 (2002),
no. 2, 247–269

\bibitem{Gor} Claudio Gorodski, Topics in Polar Actions, arXiv:2208.03577v3, 2022
\bibitem{Davis} Davis, Michael, Lectures on orbifolds and reflection groups, 2015  
\bibitem{Fourier} Tullio Ceccherini-Silberstein,
Fabio Scarabbotti,  
Filippo Tolli, Discrete Harmonic Analysis, Cambridge University Press, 2018
\bibitem{semidirect} Tullio Ceccherini-Silberstein,
Fabio Scarabbotti,  
Filippo Tolli, Representation Theory of Finite Groups Extensions, Springer, 2022
\bibitem{Razmyslov} Yu. P. Razmyslov, On a problem of Kaplansky, Izv. Akad. Nauk SSSR, Ser. Mat.
37 (1973), 483-501. Translation: Math. USSR, Izv. 7 (1973), 479-496.
\bibitem{binarytetrahedral} \url{https://en.wikipedia.org/wiki/Binary_tetrahedral_group}



\bibitem{Granville}
DEFECT ZERO $p$-blocks FOR FINITE SIMPLE GROUPS, 
Andrew Granville and Ken Ono
Transactions of the American Mathematical Society, 348, 1, 1996, pages 331-347.

\bibitem{Erdos} P. Erdös and P. Turan, On some problems of statistical group theory, Acta Math.
Acad. Sci. Hung. 19 (1968) 413-435
\bibitem{Webb} Peter Webb, A Course in Finite Group Representation Theory, Cambridge University Press, 2016  
\bibitem{comm}
Alexander Kushkuley,
Some Remarks on Commuting Probability,
arXiv:2510.04458, \url{https://arxiv.org/abs/2510.04458}, 2025

\bibitem{OV} A. L. Onishchik, E. B. Vinberg,
Lie Groups
and Algebraic Groups, Springer,  1990
\bibitem{LP} M. Larsen and R. Pink, Determining representations from invariant dimensions,
Invent. Math. 102 (1990), 377–398

\bibitem{Larsen} Michael Larsen, Rigidity in the Invariant Theory
of Compact Groups
arXiv:math/0212193v1,  2002

\bibitem{KR} A. H. Kushkuly, Yu. P. Razmyslov, 
Varieties generated by irreducible Representations of Lie Algebras, Vestnik  Moskovskogo Universiteta, 1983, no 5, pp. 4-7

\bibitem{R}
Yu.P. Razmyslov, Identities of Algebras and Their Representations, Nauka, Moscow, 1989, AMS, 1994

\bibitem{Zusmanovich}PASHA ZUSMANOVICH, ON THE UTILITY OF ROBINSON–AMITSUR ULTRAFILTERS,
arXiv:0911.5414v8, 2016

\bibitem{DR} V.S. Drensky and M.L. Racine, Distinguishing simple Jordan algebras by means of polynomial identities, Comm.
Algebra 20 (1992), 309–327.

\bibitem{Bourb9} N. BOURBAKI, GROUPES ET ALGEBRES DE LIE, ch. 9, MASSON, 1982
\bibitem{Sch} Konrad Schm\.{u}dgen,
Unbounded Operator
Algebras and
Representation Theory, Springer Basel, 1990 
\bibitem{Drinfeld} MITYA BOYARCHENKO AND VLADIMIR DRINFELD, A MOTIVATED INTRODUCTION TO CHARACTER SHEAVES AND
THE ORBIT METHOD FOR UNIPOTENT GROUPS IN POSITIVE
CHARACTERISTIC, arXiv:math/0609769v2, http://arxiv.org/abs/math/0609769v2, 2010

\end{thebibliography}
\end{document}